%% file: dv.tex
\documentclass[a4paper,10pt]{article}
\usepackage[english,frenchb]{babel}
\usepackage[latin1]{inputenc}
\usepackage[T1]{fontenc}
\usepackage{amsmath}
\usepackage{array}
\usepackage{amsthm}
\usepackage{amssymb}
\input xy
\xyoption{all}

\newcommand{\A}{{\mathcal{A}}}
\newcommand{\B}{{\mathcal{B}}}
\newcommand{\C}{{\mathcal{C}}}
\newcommand{\D}{{\mathcal{D}}}
\newcommand{\F}{{\mathcal{F}}}
\newcommand{\E}{{\mathcal{E}}}
\newcommand{\Eq}{{\mathcal{E}_q}}
\newcommand{\Eqd}{{\mathcal{E}_q^{deg}}}

\newcommand{\FF}{{\mathbb{F}_2}}
\newcommand{\Gr}{{\mathcal{G}r}}

\newcommand{\col}{{\rm colim}\,}

\newcommand{\kk}{{\Bbbk}}

\title{Sur l'homologie des groupes orthogonaux et symplectiques \`a coefficients tordus}

\author{Aur\'elien DJAMENT\thanks{CNRS, laboratoire de math\'ematiques Jean Leray (Nantes) ; djament@math.univ-nantes.fr.}\; et Christine VESPA\thanks{Institut de Recherche Math\'ematique Avanc\'ee, universit\'e de Strasbourg ; vespa@math.u-strasbg.fr.}\;\thanks{Les auteurs ont \'et\'e partiellement soutenus par les contrats INTAS 06-1000017-8609 et ANR BLAN08-2-338236 (HGRT : {\em nouveaux liens entre la th\'eorie de l'homotopie et la th\'eorie des groupes et des repr\'esentations}).}}

\newtheorem{thm-intro}{Théorème}

\newtheorem{thm}{Théorème}[section]
\newtheorem{pr}[thm]{Proposition}
\newtheorem{cor}[thm]{Corollaire}
\newtheorem{lm}[thm]{Lemme}

\theoremstyle{definition}
\newtheorem{defi}[thm]{Définition}
\newtheorem{nota}[thm]{Notation}
\newtheorem{hyp}[thm]{Hypothèse}
\newtheorem{conv}[thm]{Convention}

\theoremstyle{remark}
\newtheorem{rem}[thm]{Remarque}
\newtheorem{ex}[thm]{Exemple}

\begin{document}

\maketitle

\input{intro-dv}

\section{Cadre formel}\label{s-un}
Dans cette section, on donne le cadre g\'en\'eral de cet article, qui permet de traiter de l'homologie stable des groupes orthogonaux ou symplectiques. On verra, dans les appendices \ref{appBS} et \ref{apsym}, que ce cadre s'applique également aux groupes lin\'eaires et sym\'etriques. 

\subsection{Hypothèses générales}\label{p-un}

On introduit ici des axiomes que l'on supposera vérifiés dans tout l'article. Les exemples et des hypothèses supplémentaires souvent utiles seront donnés dans le paragraphe suivant --- en effet, tous les cas intéressants d'application de la situation générale, hormis un cas technique apparaissant en cours de démonstration, relèveront desdites hypothèses supplémentaires. 

On se donne une catégorie (essentiellement) petite $\C$ et deux foncteurs $S : \mathbb{N}\to\C$ et $G : \mathbb{N}\to\mathbf{Grp}$, où $\mathbb{N}$ désigne la catégorie associée à l'ensemble ordonné $\mathbb{N}$ (il y a exactement une flèche d'un entier positif $i$ vers un autre entier positif $j$ si $i\leq j$, et aucune sinon) et $\mathbf{Grp}$ la catégorie des groupes.

Avant de donner les trois hypothèses que nous ferons sur $(\C,S,G)$, signalons que l'on pourrait remplacer $\mathbb{N}$ par un autre ensemble ordonné filtrant à droite et possédant un plus petit élément sans modifier la plupart des considérations qui suivent. Cependant, une telle généralisation semble présenter un intérêt modeste dans la mesure où l'on ne connaît aucun exemple qui ne puisse se ramener au cas ici décrit (remplacer l'ensemble ordonné par une partie cofinale à droite contenant le plus petit élément ne modifie guère la situation).

Notre premier axiome est relatif au foncteur $S$, dont il exprime une sorte de propriété de cofinalité :\\
(C) {\em pour tout objet $c$ de $\C$, il existe $i\in\mathbb{N}$ et un morphisme dans $\C$ de source $c$ et de but $S(i)$.}\\

Notre second axiome renvoie à une forme stable de transitivité de l'action au but des automorphismes sur un ensemble de morphismes:\\
(W) {\em étant donnés $i\in\mathbb{N}$ et $c\in {\rm Ob}\,\C$, pour tous morphismes $u, v : c\to S(i)$ de $\C$, il existe $j\geq i$ dans $\mathbb{N}$ et $g\in {\rm Aut}_\C S(j)$ tels que le diagramme
$$\xymatrix{c\ar[r]^u\ar[rd]_v & S(i)\ar[r]^{S(i\leq j)} & S(j)\ar[d]^g \\
 & S(i)\ar[r]^{S(i\leq j)} & S(j)
}$$
commute.}

Dans de nombreux cas, l'axiome suivant, qui implique (W), qui exprime une transitivité instable, sera vérifié :\\
(W') {\em pour tous $c\in {\rm Ob}\,\C$ et $i\in\mathbb{N}$, le groupe ${\rm Aut}_\C S(i)$ opère transitivement sur l'ensemble ${\rm Hom}_\C(c,S(i))$.}
 
L'archétype de résultat fournissant ce type de propriété est le théorème de Witt (cf. paragraphe suivant).

On astreint enfin le foncteur $G$ à satisfaire la propriété de compatibilité au foncteur~$S$ suivante  :\\
(G) {\em sur les objets, $G$ est donné par $G(i)={\rm Aut}_\C (S(i))$ pour tout $i\in\mathbb{N}$. De plus, on suppose que pour tous entiers $i$ et $j$ tels que $i\leq j$, le morphisme $S(i\leq j) : S(i)\to S(j)$ est $G(i)$-équivariant, où $S(i)$ est muni de l'action à gauche tautologique de $G(i)$ et $S(j)$ de celle déduite du morphisme $G(i\leq j) : G(i)\to G(j)$}.

\begin{hyp}
Dans la suite du paragraphe~\ref{p-un}, on suppose que $(\C,S,G)$ vérifie les hypothèses (C), (W) et (G).
\end{hyp}

\begin{nota}
La colimite du foncteur $G$ est not\'ee $G_{\infty}$. Pour tout foncteur $F: \C\to\mathbf{Mod}_\kk$, on note $F_{\infty}$ la colimite du foncteur compos\'e $\mathbb{N}\xrightarrow{S}\C\xrightarrow{F}\mathbf{Mod}_\kk$.
\end{nota}

\begin{rem}\label{rqact}
Pour tout foncteur 
$F: \C\to\mathbf{Mod}_\kk$ et tous entiers $i\leq j$, l'application lin\'eaire $(F\circ S)(i\leq j): F(S(i)) \to F(S(j))$ est $G(i)$-\'equivariante, par l'axiome (G). Ainsi, $F_{\infty}$ est naturellement un $\kk[G_{\infty}]$-module.
\end{rem}

On peut donc donner la d\'efinition suivante, qui introduit l'objet d'étude de cet article.
\begin{defi}
L'{\em homologie stable} de la suite de groupes $( G(i) )_{i\in\mathbb{N}}$ \`a coefficients dans $F\in {\rm Ob}\,\C-\mathbf{Mod}$ est $H_*(G_{\infty} ; F_{\infty})$.
\end{defi}

\begin{rem}
Le caractère filtrant à droite de $\mathbb{N}$ implique que le morphisme canonique\linebreak $\col_{i\in\mathbb{N}} \ H_*(G(i) ; F(S(i))) \to H_*(G_{\infty} ; F_{\infty})$ est un isomorphisme. Dans la suite on identifiera les deux groupes via cet isomorphisme.
\end{rem}

 On rappelle que l'on note $P_c^{\C^{op}}$  le foncteur $\kk[{\rm Hom}_{\C^{op}}(c,-)]$, appel\'e projectif standard de $\mathbf{Mod}-\C$ (voir appendice \ref{appA}). On termine ce paragraphe par quelques r\'esultats généraux sur la colimite de ces projectifs standards.

\begin{lm} \label{Pinf}
Le foncteur
$$P_{\infty}^{\C^{op}}:=\underset{i\in\mathbb{N}}{\col} P_{S(i)}^{\C^{op}}$$
 est plat et est muni d'une action du groupe $G_{\infty}$.
 
 De plus, il existe un isomorphisme de $G_\infty$-modules
 $$P_{\infty}^{\C^{op}}\underset{\C}{\otimes}F\simeq F_\infty$$
 naturel en $F\in {\rm Ob}\,\C-\mathbf{Mod}$.
\end{lm}
\begin{proof}
La platitude d\'ecoule du caractère filtrant de $\mathbb{N}$. Comme $G(i)$ agit sur $P_{S(i)}^{\C^{op}}$, on en d\'eduit \'egalement une action canonique de $G_{\infty}$ sur $P_{\infty}^{\C^{op}}$.

La dernière partie se vérifie immédiatement.
\end{proof}

\begin{lm} \label{Pinf2}
On a 
$$(P_{\infty}^{\C^{op}})_{G_{\infty}}\simeq \kk$$
o\`u $(P_{\infty}^{\C^{op}})_{G_{\infty}}$ d\'esigne les coïnvariants de $P_{\infty}^{\C^{op}}$ par l'action de $G_{\infty}$ et $\kk$ est le foncteur constant.
\end{lm}
\begin{proof}
Il s'agit d'une conséquence directe de l'hypoth\`ese (W).
\end{proof}

\subsection{Cas particuliers fondamentaux}\label{p1-deux}

Dans ce paragraphe, on considère une catégorie (essentiellement) petite $\C$ munie d'une structure monoïdale symétrique $\oplus : \C\times\C\to\C$ dont l'unité sera notée $0$. On fait également l'hypothèse que $0$ est objet initial de $\C$.

Quitte à remplacer $\C$ par une cat\'egorie monoïdale symétrique \'equivalente, on pourra supposer que le foncteur $\oplus$ est strictement associatif et que $0$ en est un élément neutre strict, ce qui permet de donner un sens univoque à des expressions comme $A^{\oplus n}$, où $n\in\mathbb{N}$ et $A\in {\rm Ob}\,\C$.

Soit $A\in {\rm Ob}\,\C$. On peut définir un foncteur $S_A : \mathbb{N}\to\C$ par $S_A(n)=A^{\oplus n}$ et
$$S_A(n\leq m) : A^{\oplus n}=A^{\oplus n}\oplus 0\xrightarrow{A^{\oplus n}\oplus (0\to A^{\oplus (m-n)})}A^{\oplus n}\oplus A^{\oplus (m-n)}=A^{\oplus m}.$$

Ce choix de fonctorialité consiste, lorsque $\C$ est une catégorie de modules (cf. exemple~\ref{exfond}.\,\ref{ex-sco} ci-dessous), par exemple, à "ajouter des zéros à droite" en termes matriciels.

On définit également un foncteur $G_A : \mathbb{N}\to\mathbf{Grp}$ par $G_A(n)={\rm Aut}_\C(A^{\oplus n})$ et
$$G_A(n\leq m) : {\rm Aut}_\C(A^{\oplus n})\to {\rm Aut}_\C(A^{\oplus m})\quad u\mapsto u\oplus A^{\oplus (m-n)}.$$

On vérifie aussitôt le fait suivant :
\begin{pr}\label{prevg}
Le triplet $(\C,S_A,G_A)$ vérifie l'hypothèse (G).
\end{pr}

Nous aurons également à considérer l'hypothèse suivante :\\
(S) {\em pour tout morphisme $f : d\to c$ de $\C$ et tout $i\in\mathbb{N}$, le morphisme du stabilisateur de $f$ sous l'action de ${\rm Aut}_\C(c)$ vers le stabilisateur de $S(i)\oplus f$ sous l'action de ${\rm Aut}_\C(S(i)\oplus c)$ induit par le foncteur $S(i)\oplus -$ est un isomorphisme.}

L'hypothèse (S) permettra de donner des renseignements supplémentaires sur la deuxième page de la suite spectrale pour l'homologie de $G_\infty$ à coefficients tordus que nous obtiendrons dans la section~\ref{secss} lorsque les axiomes (C), (W) et (G) sont vérifiés.

Nous introduisons enfin une hypothèse plus forte que (C) mais moins forte que l'essentielle surjectivité du foncteur $S_A$ :\\
(C') {\em pour tout objet $c$ de $\C$, il existe un objet $b$ de $\C$ et un entier $i$ tels que $b\oplus c\simeq S(i)$.}

\smallskip

Nous terminons ce paragraphe en donnant les exemples fondamentaux qui interviendront dans cet article.

\begin{ex}\label{exfond}
\begin{enumerate}
\item\label{exquad} Soient $k$ un corps commutatif (éventuellement de caractéristique $2$)  et $\Eqd(k)$ la catégorie (qui sera notée simplement $\Eqd$ lorsqu'aucune confusion ne peut en résulter) dont les objets sont les espaces quadratiques de dimension finie sur $k$ et les morphismes les applications lin\'eaires {\em injectives} compatibles aux formes quadratiques (appelées aussi applications orthogonales). Remarquons qu'une forme quadratique sur un espace vectoriel $V$ est un polyn\^ome homog\`ene de degr\'e $2$ sur $V$ (\cite{Pfister}), c'est \`a dire un \'el\'ement de $S^2(V^*)$, o\`u $S^2$ est la deuxi\`eme puissance sym\'etrique. Comme d'habitude, on note $O(A)$ pour ${\rm Aut}_\Eqd (A)$ le groupe orthogonal associé à un objet $A$ de $\Eqd$. On notera par ailleurs $\Eq$ la sous-catégorie pleine de $\Eqd$ dont les objets sont les espaces quadratiques non dégénérés.

La somme orthogonale, notée $\perp$, définit une structure monoïdale symétrique sur $\Eqd$ dont l'unité $0$ est objet initial de $\Eqd$.
 
Soit $\mathbf{H}$ l'objet de $\Eq$, appel\'e plan hyperbolique, dont l'espace vectoriel sous-jacent est $k^2$ et la forme quadratique l'application $k^2\to k\quad (x,y)\mapsto xy$. Comme d'habitude, on notera $O_{n,n}(k)$ le groupe $G_{\mathbf{H}}(n)$ des automorphismes de $\mathbf{H}^{\perp n}$.

Le triplet $(\Eqd, S_{\mathbf{H}}, G_{\mathbf{H}})$ vérifie l'hypothèse (C), car tout espace quadratique se plonge dans un espace quadratique non dégénéré et tout espace quadratique non dégénéré se plonge dans un espace hyperbolique (i.e. un espace quadratique qui est somme directe de sous-espaces totalement isotropes), qui est isomorphe à une somme orthogonale de copies de $\mathbf{H}$ (cf. \cite{Scharlau}). Le triplet $(\Eq, S_{\mathbf{H}}, G_{\mathbf{H}})$ vérifie pour sa part (C').

Le théorème de Witt montre que l'axiome (W') est satisfait, puisque $S_{\mathbf{H}}$ prend ses valeurs dans les espaces non dégénérés.

Le triplet $(\Eqd, S_{\mathbf{H}}, G_{\mathbf{H}})$ vérifie également l'hypothèse (S) pour la même raison. En effet, si $V$ est un espace quadratique et $H$ un espace non dégénéré, tout automorphisme de $V\perp H$ qui préserve $H$ préserve également son supplémentaire orthogonal, qui n'est autre que $V$ car $H$ est non dégénéré. Par conséquent, tout élément du stabilisateur d'un morphisme $H\perp f : H\perp V\to H\perp W$ sous l'action de $O(H\perp V)$ stabilise $H\subset H\perp V$ et $W$, donc s'écrit sous la forme $H\perp u$, où $u\in O(W)$ stabilise $f : V\to W$ ; $u$ est manifestement unique, d'où la satisfaction de (S).

 (On pourrait étendre ces considérations à un corps gauche muni d'une involution et aux formes hermitiennes afférentes.)
 
 \begin{rem}\label{rq-goinf}
  On peut remplacer $\mathbf{H}$ par n'importe quel autre $k$-espace quadratique non dégénéré de dimension finie $H$ tel qu'on puisse plonger tout autre $k$-espace quadratique non dégénéré de dimension finie dans une somme de copies de $H$. Si le corps $k$ est fini (auquel cas le treillis des classes d'isomorphisme de $k$-espaces quadratiques non dégénérés de dimension finie est particulièrement simple !), n'importe quel objet non nul $H$ de $\Eq$ convient.
 \end{rem}

\item\label{exsymp}  On peut reprendre mutatis mutandis l'exemple précédent en remplaçant les formes quadratiques par les formes bilinéaires alternées. Bornons-nous à préciser nos notations : $\E_{alt}^{deg}(k)$ (ou simplement $\E_{alt}^{deg}$) désignera la catégorie des espaces symplectiques de dimension finie (avec les {\em injections} symplectiques pour morphismes) sur le corps commutatif $k$, $\E_{alt}$ la sous-catégorie pleine des espaces non dégénérés, $\perp$ la somme orthogonale.

On notera aussi $\mathbf{H}$ l'espace symplectique $k^2$ muni de la forme $((x,y),(x',y'))\mapsto xy'-yx'$. Le groupe des automorphismes de $A\in {\rm Ob}\,\E_{alt}^{deg}$ sera noté $Sp(A)$ ; $G_{\mathbf{H}}(n)=Sp(\mathbf{H}^{\perp n})$ sera noté $Sp_{2n}(k)$.

\item\label{exbs} Soit $\Theta$ la catégorie ayant pour objets les ensembles finis et pour morphismes les fonctions injectives. La structure  monoïdale sym\'etrique de $\Theta$ est donn\'ee par la r\'eunion disjointe ; son unit\'e l'ensemble vide est objet initial de $\Theta$. Soit $A$ un ensemble à un élément fixé. La condition (C') est v\'erifi\'ee puisque $S_A$ est essentiellement surjectif. On a $G_A(n)=\mathfrak{S}_n$ (groupe symétrique sur $n$ lettres) ; on voit facilement que la condition (W') est satisfaite. La condition (S) a lieu du fait qu'une bijection d'un ensemble $E$ qui pr\'eserve un sous-ensemble $F$ pr\'eserve \'egalement son compl\'ementaire.

\item\label{ex-sco} Soient $A$ un anneau et $\mathbb{M}(A)$ la catégorie des $A$-modules à gauche  libres de type fini  avec pour morphismes les injections $A$-linéaires scindées, munie de la somme directe, dont l'unité $0$ est objet initial. Le triplet $(\mathbb{M}(A),S_A, G_A)$ vérifie l'hypothèse (C'). Comme d'habitude, on note $GL_n(A)$ pour $G_A(n)$.

L'hypothèse (W) est aussi satisfaite, mais pas (S) si $A$ est non nul et pas non plus (W') en général : pour s'en convaincre, on peut considérer un anneau non nul $A$ tel que les $A$-modules à gauche $A$ et $A^2$ soient isomorphes, par exemple l'anneau des endomorphismes d'un espace vectoriel de dimension infinie, de sorte que l'injection scind\'ee évidente $A\hookrightarrow A^2$, qui n'est pas surjective, ne saurait être conjuguée sous $GL_2(A)$ à un isomorphisme $A\xrightarrow{\simeq}A^2$ ((W') est cependant vraie si $A$ est un anneau assez gentil, par exemple principal). La satisfaction de (W) d\'ecoule du r\'esultat classique de stabilisation suivant:

\begin{lm}\label{lmsthom} Soient $\A$ une catégorie additive, $M$ et $N$ deux objets de $\A$ et $u : N\to M$ un monomorphisme scindé de $\A$. 

Il existe un automorphisme $g$ de $N\oplus M$ tel que le diagramme
$$\xymatrix{N\ar[r]^u\ar[rrd] & M\ar[r] & N\oplus M \\
 & & N\oplus M\ar[u]_g^\simeq
}$$
dans lequel les flèches non spécifiées sont les inclusions canoniques, commute.
\end{lm}

\begin{proof}
Soit $p : M\to N$ une rétraction de $u$. L'endomorphisme $g$ de $N \oplus M$ donné par la matrice $\left(\begin{array}{cc} 0 & p \\
u & id
\end{array}\right)$ fait commuter le diagramme, et c'est un automorphisme d'inverse donné par $\left(\begin{array}{cc} -id & p \\
u & id-up
\end{array}\right)$.
\end{proof}
\end{enumerate}
\end{ex}

\section{Suite spectrale pour l'homologie stable d'une suite de groupes}\label{secss}

\input{spec}

\section{Cas des groupes orthogonaux et symplectiques}\label{sid}

\input{dv-s3}

\section{Quelques calculs d'homologie de groupes orthogonaux et symplectiques}\label{scal}

Nous allons maintenant illustrer les résultats de la section précédente par des calculs explicites de groupes d'homologie de groupes orthogonaux et symplectiques.

\begin{conv}
{\bf Dans toute cette section, on suppose que $\kk=k$ est un corps fini de caractéristique $p$ et de cardinal $q=p^d$.

On rappelle qu'on note simplement $\F(k)$ ou $\F$ la catégorie de foncteurs~$\E^f(k)-\mathbf{Mod}$.}

Les mentions au corps $k$ seront souvent omises.
\end{conv}

 On s'est limité au cas $k=\kk$ car tout foncteur polynomial sans terme constant de $\E^f_k$ vers les groupes abéliens prend ses valeurs dans les $\mathbb{F}_p$-espaces vectoriels ; une extension des scalaires au but ne modifie guère le comportement homologique de notre catégorie de foncteurs. De surcroît, c'est le cas où $\kk$ égale $k$ qu'il est usuel d'étudier ; des calculs d'algèbre homologique poussés y ont été effectués (cf. \cite{FLS} et \cite{FFSS}).

\subsection{Compatibilité aux (co)produits}

Si $G$ est un groupe et $M$ et $N$ sont deux $G$-modules, la projection canonique $(M\otimes N)_G\twoheadrightarrow M_G\otimes N_G$ induit un {\em coproduit externe} en homologie $H_*(G ; M\otimes N)\to H_*(G ; M)\otimes H_*(G ; N)$. Ce morphisme naturel gradué fait de $H^*(G;-)$ un foncteur comonoïdal. Il est plus usuel de considérer la situation duale, à savoir le produit externe
$$H^*(G; M)\otimes H^*(G ; N)\to H^*(G; M\otimes N)$$
qui fait de $H^*(G;-)$ un foncteur monoïdal des $G$-modules vers les espaces vectoriels gradués et se réduit en degré $0$ à l'inclusion canonique $M^G\otimes N^G\hookrightarrow (M\otimes N)^G$.

Des constructions analogues existent en (co)homologie des foncteurs : le produit tensoriel induit en particulier, pour toute petite catégorie $\C$, des produits naturels
$${\rm Ext}_{\C}^*(A,F)\otimes {\rm Ext}_{\C}^*(B,G)\to {\rm Ext}_{\C}^*(A\otimes B,F\otimes G).$$
(On a des morphismes duaux évidents dans le cas des groupes de torsion.)

Lorsque le foncteur $A$ est muni d'une structure de {\em cogèbre}, on peut utiliser les morphismes ${\rm Ext}^*(A\otimes A,F\otimes G)\to {\rm Ext}^*(A,F\otimes G)$ induits par le coproduit $A\to A\otimes A$ pour en déduire un produit
$${\rm Ext}^*(A,F)\otimes {\rm Ext}^*(A,G)\to {\rm Ext}^*(A,F\otimes G).$$

Tous les foncteurs de précomposition, donc en particulier d'évaluation, sont compatibles aux (co)produits ainsi définis. 
Cela permet d'obtenir la compatibilité de tous les isomorphismes de la section précédente aux produits ou coproduits. Nous nous bornerons ici à l'énoncé suivant :
\begin{pr}\label{prodiso} Les isomorphismes du corollaire~\ref{corco} sont compatibles aux produits externes. Plus précisément, soient $F$ et $G$ des foncteurs polynomiaux de $\F$ de degrés respectifs $d$ et $d'$ et $n$, $i$, $j$ des entiers tels que $n\geq 2(i+j)+d+d'+6$. Supposons la caractéristique $p$ de $k$ impaire. Le diagramme
$$\xymatrix{H^i(O_{n,n}(k);F(k^{2n}))\otimes H^j(O_{n,n}(k);G(k^{2n}))\ar[r]\ar[d]^\simeq & H^{i+j}(O_{n,n}(k);(F\otimes G)(k^{2n}))\ar[d]^\simeq \\
{\rm Ext}^i_{\F(k)}(k[S^2],F)\otimes {\rm Ext}^j_{\F(k)}(k[S^2],G)\ar[r] & {\rm Ext}^{i+j}_{\F(k)}(k[S^2],F\otimes G)}$$
dans lequel les flèches horizontales sont les produits et les flèches verticales les isomorphismes donnés par le corollaire~\ref{corco} commute. On dispose d'un énoncé analogue évident en termes de groupes symplectiques.
\end{pr}

\begin{proof} Compte-tenu de la remarque précédente et de la description des isomorphismes, il suffit de démontrer que l'isomorphisme d'adjonction ${\rm Tor}^{\E^f_{inj}}_*(k[S^2]^\vee,F)\simeq H_*(\Eqd ; F)$ est compatible aux coproduits (notre énoncé s'en déduit en dualisant). Celui-ci est composé du morphisme $H_*(\Eqd ; F)={\rm Tor}^\Eqd_*(k;F)\to {\rm Tor}^\Eqd_*(k[S^2]^\vee;F)$ induit par la flèche $k\to k[S^2]^\vee$ de $\mathbf{Mod}-\Eqd$ donnée, sur un espace quadratique $(V,q)$, par l'élément $[q]$ de $k[S^2(V^*)]$, et du morphisme ${\rm Tor}^\Eqd_*(k[S^2]^\vee;F)\to {\rm Tor}^{\E^f_{inj}}_*(k[S^2]^\vee;F)$ induit par le foncteur d'oubli de la forme quadratique $\Eqd\to\E^f_{inj}$. Le premier respecte les coproduits car $k\to k[S^2]^\vee$ est un morphisme de cogèbres, le second par l'observation générale précédant la démonstration, d'où la proposition.
\end{proof}

Cette proposition permet d'obtenir une propriété générale frappante des produits externes en cohomologie stable des groupes orthogonaux, à l'aide d'observations élémentaires mais efficaces dues à Touzé --- la proposition suivante est établie (dans le contexte analogue des foncteurs polynomiaux stricts) dans~\cite{Touze}.

\begin{pr}[Touzé]\label{prtouze} Soient $\A$ une petite catégorie additive et $A$ un objet de $\A - \mathbf{Mod}$ muni d'une structure de cogèbre et d'un épimorphisme naturellement scindé $A(V\oplus W)\twoheadrightarrow A(V)\otimes A(W)$ pour $V, W\in {\rm Ob}\,\A$ de sorte que le coproduit de $A$ soit donné par la composition $A(V)\to A(V\oplus V)\twoheadrightarrow A(V)\otimes A(V)$, où la première flèche est induite par la diagonale $V\to V\oplus V$. On suppose de surcroît soit que $A$ possède une résolution projective de type fini, soit que $F$ et $G$ possèdent une résolution injective de type fini et que $A$ prend des valeurs de dimension finie. Alors le produit externe
$${\rm Ext}^*(A,F)\otimes {\rm Ext}^*(A,G)\to {\rm Ext}^*(A,F\otimes G)$$
est une injection naturellement scindée pour tous objets $F$ et $G$ de $\A-\mathbf{Mod}$.
\end{pr}

\begin{proof}
À l'aide de l'épimorphisme scindé de l'hypothèse, des foncteurs adjoints $\oplus : \A\times\A\to\A$ et $\delta : \A\to\A\times\A$ (diagonale) et de la formule de Künneth (cf. par exemple \cite{FFSS}, propriété (1.7.2)\,\footnote{Dans cet énoncé, il n'est fait mention de la formule de Künneth que dans le premier cas de finitude envisagé ; l'autre se traite de façon analogue.}), on obtient un monomorphisme naturellement scindé
$${\rm Ext}^*_\A(A,F)\otimes {\rm Ext}^*_\A(A,G)\simeq {\rm Ext}^*_{\A\times\A}(A\boxtimes A,F\boxtimes G)\hookrightarrow {\rm Ext}^*_{\A\times\A}(\oplus^*A,F\boxtimes G)$$
$$\simeq {\rm Ext}^*_\A(A,\delta^*(F\boxtimes G))={\rm Ext}^*_\A(A,F\otimes G).$$

La description du coproduit sur $A$ à partir de l'épimorphisme scindé $\oplus^*A\to A\boxtimes A$ contenue dans l'hypothèse montre que la flèche précédente n'est autre que le produit de l'énoncé.
\end{proof}

On rappelle au lecteur que les g\'en\'eralit\'es concernant les foncteurs exponentiels, dont on fait un usage fr\'equent dans la suite, sont donn\'ees dans l'appendice~\ref{apex}.

L'hypothèse sur le foncteur $A$ est en particulier vérifiée, pour $\A=\E^f$, lorsque $A$ est la composition $E\circ T$ d'un foncteur exponentiel $E$ et d'un foncteur $T$ tel que $T(0)=0$ (puisqu'alors $T(U\oplus V)$ contient $T(U)\oplus T(V)$ comme facteur direct naturel). Par conséquent, les propositions~\ref{prodiso} et~\ref{prtouze} procurent le résultat d'injectivité suivant (où l'on utilise que les foncteurs polynomiaux à valeurs de dimension finie possèdent des résolutions injectives de type fini --- cf. \cite{FLS}, §\,10, où l'énoncé est donné sous forme duale).

\begin{pr}\label{inj-ext} Soient $F$ et $G$ deux foncteurs polynomiaux de $\F(k)$ prenant des valeurs de dimension finie, de degrés respectifs $d$ et $d'$, $i$, $j$ et $n$ des entiers tels que $n\geq 2(i+j)+ d+d'+6$. Alors le produit externe 
$$H^i(O_{n,n}(k);F(k^{2n}))\otimes H^j(O_{n,n}(k);G(k^{2n}))\to H^{i+j}(O_{n,n}(k);(F\otimes G)(k^{2n}))$$
est injectif si la caract\'eristique de $k$ est impaire.
\end{pr}

\begin{rem}
On a en fait mieux (cf. \cite{Touze}, §\,6.1) : il existe en cohomologie stable des groupes orthogonaux ou symplectiques un {\em co}produit externe qui fournit une rétraction naturelle au produit externe, comme on s'en convainc aisément en reprenant la démonstration. Il est impossible de décrire directement ce coproduit, qui n'existe pas en cohomologie instable.
\end{rem}

\begin{rem}
 Il est en fait inutile d'invoquer les propriétés de finitude des foncteurs polynomiaux à valeurs de dimension finie. En effet, la variante duale (en termes de groupes de torsion) de la proposition~\ref{prtouze} est valide sans aucune hypothèse de finitude (comme la formule de Künneth pour les Tor) ; les groupes d'extensions ici considérés viennent tous de la dualisation de groupes de torsion. Nous avons néanmoins privilégié les énoncés en termes de groupes d'extensions, plus usuels et plus intuitifs.
\end{rem}

\smallskip

Nous donnons maintenant quelques résultats préliminaires aux calculs explicites de (co)homologie stabilisée de groupes orthogonaux et symplectiques tordus par des foncteurs polynomiaux classiques. Nos calculs concerneront les foncteurs exponentiels gradués usuels : puissances extérieures, divisées, symétriques.

\begin{lm} \label{iso-Hoch}
Soit $F$ un objet de $\F$. Il existe un isomorphisme gradu\'e naturel
$${\rm Tor}_*^{\E^f}(k[T^2]^\vee,F)\simeq HH_*(\E^f;(V,W)\mapsto F(V^*\oplus W)).$$
\end{lm}

\begin{proof}
Utilisant l'adjonction entre les foncteurs exacts induits par la somme directe et la diagonale, l'équivalence de catégories $(-)^* : (\E^f)^{op}\to\E^f$ et l'isomorphisme naturel $V^*\otimes W\simeq {\rm Hom}_{\E^f}(V,W)$, puis l'isomorphisme~(\ref{hochtor}) de l'appendice~\ref{appA}, on obtient des isomorphismes gradués naturels
$${\rm Tor}_*^{\E^f}(k[T^2]^\vee,F)\simeq{\rm Tor}_*^{(\E^f)^{op}\times\E^f}(k[{\rm Hom}_{(\E^f)^{op}}],(V,W)\mapsto F(V^*\oplus W))$$
$$\simeq HH_*(\E^f;(V,W)\mapsto F(V^*\oplus W)).$$
\end{proof}

\begin{rem}
Pour un foncteur polynomial $F$, par le th\'eor\`eme de Betley-Suslin \cite{Bet} \cite[appendice]{FFSS}, on a un isomorphisme naturel 
$$HH_*(\E^f;(V,W)\mapsto F(V^*\oplus W)) \simeq \underset{n \in \mathbb{N}}{\col}H_*(GL_n(k); F(k^n \oplus k^n))$$
o\`u $g \in GL_n(k)$ 	agit sur $k^n \oplus k^n$ par $(^t g^{-1},g)$. 

On v\'erifie que le diagramme suivant commute
$$\xymatrix{
{\rm Tor}_*^{\E^f}(k[T^2]^\vee,F) \ar[r]^-\simeq \ar[d]_{\alpha} & HH_*(\E^f;(V,W)\mapsto F(V^*\oplus W)) \ar[r]^-\simeq&\underset{n \in \mathbb{N}}{\col}H_*(GL_n(k); F(k^n \oplus k^n)) \ar[d]^-{\beta}\\
{\rm Tor}_*^{\E^f}(k[S^2]^\vee,F) \ar[rr]^-\simeq &&  \underset{n \in \mathbb{N}}{\col}H_*(O_{n,n}(k); F(k^n \oplus k^n)) 
}$$
o\`u $\alpha$ est induit par le morphisme canonique $T^2 \to S^2$, $\beta$ par le morphisme de groupes 
$$GL_n(k) \to O_{n,n}(k) \qquad g \mapsto \left(\begin{array}{cc} 
 ^t g^{-1}  & 0 \\
0 & g 
\end{array}\right)$$
et l'isomorphisme du bas est induit par le corollaire \ref{crf-o}.
\end{rem}

On note $D : \F^{op}\to\F$ le foncteur de dualité de $\F$ composé commutatif de $(-)^\vee$ et $(-)^*$ (cf. notation~\ref{not-gi}), c'est-à-dire donné par $(DF)(V)=F(V^*)^*$. On remarque que ce foncteur vérifie la propriété d'auto-adjonction ${\rm Hom}_\F(F,DG)\simeq {\rm Hom}_\F(G,DF)$, qui s'étend aux groupes d'extensions par exactitude de~$D$.

\begin{pr}\label{preg}
Soit $E^\bullet$ un foncteur exponentiel gradué de $\F$. Pour tout entier $n$, il existe un isomorphisme gradué naturel
$${\rm Tor}_*^{\E^f}(k[T^2]^\vee,E^n)\simeq\bigoplus_{i+j=n}{\rm Tor}_*^{\E^f}((E^i)^\vee,E^j).$$
Le dual de cet espace vectoriel s'identifie à
\begin{equation}\label{eq-ee}
{\rm Ext}^*_\F(E^n,I\circ T^2)\simeq\bigoplus_{i+j=n}{\rm Ext}^*_\F (E^j,DE^i).
\end{equation}
\end{pr}

\begin{proof} Le premier isomorphisme s'obtient par le lemme  \ref{iso-Hoch} en développant $E^n(V^*\oplus W)$ à l'aide de la propriété exponentielle et en utilisant l'expression de l'homologie de Hochschild d'un produit tensoriel extérieur comme groupe de torsion.

La deuxième assertion s'obtient à partir de la première et de l'isomorphisme~(\ref{dualfct}) de l'appendice~\ref{appA}.
\end{proof}

\begin{rem}
Le substitut de propriété exponentielle indiqué dans la remarque~\ref{tensexp} permet d'obtenir d'une manière analogue, pour les puissances tensorielles, un isomorphisme gradué
$${\rm Tor}_*^{\E^f}(k[T^2]^\vee,T^n)\simeq\bigoplus_{i+j=n}{\rm Tor}_*^{\E^f}((T^i)^\vee,T^j)\underset{\mathfrak{S}_i\times\mathfrak{S}_j}{\otimes}k[\mathfrak{S}_n].$$ 
Ces groupes peuvent être entièrement calculés ; on laisse au lecteur le soin d'écrire les analogues des résultats de la suite de cette section en termes de puissances tensorielles.
\end{rem}

Par la suite, nous utiliserons systématiquement la forme duale en termes d'Ext, car tous les calculs effectués dans la catégorie $\F$ ont été donnés en termes de groupes d'extensions et non de torsion (cf. \cite{FFPS} par exemple).

Dans l'énoncé suivant, les produits et coproduits sont {\em internes} ; ils s'obtiennent à partir des structures externes en utilisant la structure de foncteur de Hopf.

\begin{pr}\label{lmhc}
Soit $E^\bullet$ un foncteur exponentiel de Hopf {\em commutatif ou anticommutatif} de $\F$. L'isomorphisme bigradué
$${\rm Ext}^*_\F(E^\bullet,I\circ T^2)\simeq {\rm Ext}^*_\F(E^\bullet,DE^\bullet)$$
de la proposition~\ref{preg} est un isomorphisme d'algèbres de Hopf.
\end{pr}

\begin{proof}
La compatibilité aux unités et coünités est évidente.

L'isomorphisme
$$\gamma_{E^\bullet} : {\rm Ext}^*_\F(E^\bullet,I\circ T^2)\simeq {\rm Ext}^*_\F(E^\bullet,DE^\bullet)$$
de la proposition~\ref{preg} est naturel en $E^\bullet$, la naturalité étant relative aux morphismes respectant la structure exponentielle.
Nous obtenons la compatibilité aux produits et coproduits en montrant que $\gamma_{E^\bullet}$
est monoïdal. (On rappelle qu'un foncteur monoïdal $F$ entre deux catégories monoïdales $(\A,\otimes_A)$ et $(\B,\otimes_B)$ est un foncteur muni de morphismes naturels $F(A)\otimes_\B F(A')\to F(A\otimes_\A A')$ ; une transformation naturelle entre foncteurs monoïdaux est monoïdale si elle est compatible à ces morphismes naturels en un sens évident.) En effet, comme le foncteur exponentiel gradué $E^\bullet$ est supposé (anti)commutatif, les morphismes de multiplication $E^\bullet\otimes E^\bullet\to E^\bullet$ et de comultiplication $E^\bullet\to E^\bullet\otimes E^\bullet$ sont des morphismes de foncteurs exponentiels, ce qui permet de déduire du caractère monoïdal de $\gamma_{E^\bullet}$ sa compatibilité aux produits et coproduits en revenant à leur définition. 

La démonstration de la proposition~\ref{preg} montre que $\gamma_{E^\bullet}$ s'obtient par composition des trois isomorphismes suivants. Pour alléger, on note dans la suite de la démonstration $\sigma : \E^f\times\E^f\to\E^f$ le foncteur somme directe, $\pi : \E^f\times\E^f\to\E^f$ le foncteur produit tensoriel, $\delta : \E^f\to\E^f\times\E^f$ le foncteur diagonal et $bi-\F$ pour $\E^f\times\E^f-\mathbf{Mod}$.
\begin{enumerate}
\item l'isomorphisme
$${\rm Ext}_{bi-\F}(\sigma^*F,\pi^*I)\xrightarrow{\simeq} {\rm Ext}_\F(F,I\circ T^2)$$
naturel en $F\in {\rm Ob}\,\F$ déduit de l'adjonction entre $\sigma$ et $\delta$. C'est un isomorphisme monoïdal (où la structure monoïdale des deux foncteurs $\F^{op}\to\E$ considérés se déduit de la structure d'algèbre sur le foncteur $I$) car il est composé du morphisme naturel monoïdal ${\rm Ext}_{bi-\F}(\sigma^*F,\pi^*I)\to {\rm Ext}_\F((\sigma\delta)^*F,I\circ T^2)$ induit par le foncteur $\delta$ et du morphisme naturel monoïdal ${\rm Ext}_\F((\sigma\delta)^*F,I\circ T^2)\to {\rm Ext}_\F(F,I\circ T^2)$ induit par l'unité $id_{\E^f}\to\sigma\delta$ ;
\item  l'isomorphisme
$${\rm Ext}_{bi-\F}(\sigma^* E^\bullet,\pi^*I)\simeq {\rm Ext}_{bi-\F}(E^\bullet\boxtimes E^\bullet,\pi^*I)$$
déduit de la structure exponentielle de $E^\bullet$, qui est un morphisme naturel monoïdal de foncteurs contravariants depuis les foncteurs exponentiels vers les espaces vectoriels ;
\item l'isomorphisme
$${\rm Ext}^*_{bi-\F}(F\boxtimes G,\pi^*I)\xrightarrow{\simeq}{\rm Ext}^*_\F(F,DG)$$
naturel en $(F,G)\in {\rm Ob}\,\F\times\F$ est également monoïdal. Il suffit de le montrer en degré cohomologique nul, en lequel il prend la forme explicite suivante : à une transformation naturelle $(u_{V,W} : F(V)\otimes G(W)\mapsto I(V\otimes W))_{V,W}$ on associe la transformation naturelle donnée par la collection d'applications $F(V)\to G(V^*)^*$ adjointes aux applications
$$F(V)\otimes G(V^*)\xrightarrow{u_{V,V^*}}I(V\otimes V^*)\simeq k^{{\rm End}\,V}\to k$$
où la dernière flèche est l'évaluation en le morphisme identique de $V$. Ce morphisme est donc compatible aux produits, ce qui achève d'établir la proposition.
\end{enumerate}
\end{proof}

\smallskip

La structure du foncteur $T^2$ dépendant de la parité de la caractéristique $p$ de $k$, nous distinguons le cas $p$ impair du cas $p=2$ dans nos investigations ultérieures. Tandis que nous obtiendrons des calculs complets sur les foncteurs usuels dans le premier cas, nous ne pourrons donner que des résultats très partiels dans le second.

\subsection{Calculs en caractéristique impaire}\label{pci}

\begin{conv}
{\bf Dans tout ce paragraphe, on suppose $p$ impair.

On suppose que $E^\bullet$ est un foncteur exponentiel de Hopf gradué commutatif ou anticommutatif.}

D'une manière générale, si $u : A\to B$ est une flèche de $\F$ entre foncteurs sans terme constant (i.e. $A(0)=B(0)=0$), on note $h(u)=(I\circ u)_* : {\rm Ext}^*_\F(E^\bullet,I\circ A)\to {\rm Ext}^*_\F(E^\bullet,I\circ B)$ le morphisme d'algèbres de Hopf considéré dans la proposition~\ref{convexp}.
\end{conv}

Comme $2$ est inversible dans $k$, le foncteur $T^2$ se scinde en somme des deux foncteurs simples $\Gamma^2(\simeq S^2)$ et $\Lambda^2$. La propriété exponentielle du foncteur $I$ procure donc un isomorphisme $I\circ T^2\simeq (I\circ\Gamma^2)\otimes (I\circ\Lambda^2)$. Comme le foncteur constant en $k$ est facteur direct de $I$, $I\circ\Gamma^2$ et $I\circ\Lambda^2$ sont en particulier facteurs directs de $I\circ T^2$. Pour mener à bien nos calculs, nous avons besoin d'exprimer précisément l'effet des idempotents correspondant à cette décomposition sur les isomorphismes de la proposition~\ref{preg}.

À cette fin, on note $\tau$ l'involution du foncteur $T^2$ échangeant les deux facteurs du produit tensoriel puis $e_\Gamma=h\big(\frac{1+\tau}{2}\big)$ et $e_\Lambda=h\big(\frac{1-\tau}{2}\big)$ les deux idempotents de ${\rm Ext}^*_\F(E^\bullet,I\circ T^2)$ dont les images respectives sont ${\rm Ext}^*_\F(E^\bullet,I\circ\Gamma^2)$ et ${\rm Ext}^*_\F(E^\bullet,I\circ\Lambda^2)$. 

Le résultat suivant permettra, dans les cas usuels, de décrire l'involution $h(\tau)$ en termes explicites :

\begin{lm}\label{echt}
Dans l'isomorphisme (\ref{eq-ee}) de la proposition~\ref{preg}, l'involution $h(\tau)$ est donnée par les isomorphismes ${\rm Ext}^*_\F (E^j,DE^i)\xrightarrow{\simeq} {\rm Ext}^*_\F (E^i,DE^j)$ fournis par l'auto-adjonction du foncteur~$D$ au signe $\epsilon^{ij}$ près, où $\epsilon$ vaut $1$ si $E^\bullet$ est commutatif et $-1$ si $E^\bullet$ est anticommutatif.
\end{lm}

\begin{proof} Si $F$ est un objet de $\F=\E^f-\mathbf{Mod}$ et $G$ un objet de $\mathbf{Mod}-\E^f$, le diagramme d'isomorphismes suivant commute
$$\xymatrix{{\rm Tor}^{\E^f}_*(G,F)\ar[d] & HH_*(\E^f;G\boxtimes F)\ar[l]\ar[r]\ar[d] & {\rm Tor}_*^{(\E^f)^e}(k[{\rm Hom}_{(\E^f)^{op}}],G\boxtimes F)\ar[d] \\
{\rm Tor}^{\E^f}_*(F^\vee,G^\vee) & HH_*(\E^f;F^\vee\boxtimes G^\vee)\ar[l]\ar[r] &  {\rm Tor}_*^{(\E^f)^e}(k[{\rm Hom}_{(\E^f)^{op}}],F^\vee\boxtimes G^\vee)
}$$
où les flèches verticales sont induites par l'échange des facteurs et l'anti-équivalence de catégories $(-)^*$ et l'on a noté $(\E^f)^e=(\E^f)^{op}\times\E^f$.

Supposons maintenant $F=E^\bullet$ et $G=(E^\bullet)^\vee$ et formons le diagramme d'isomorphismes suivant
$$\xymatrix{{\rm Tor}_*^{(\E^f)^e}(k[{\rm Hom}_{(\E^f)^{op}}],(E^\bullet)^\vee\boxtimes E^\bullet)\ar[d]\ar[r] & {\rm Tor}_*^{(\E^f)^e}(k[{\rm Hom}_{(\E^f)^{op}}],\sigma^* E^\bullet)\ar[d]\ar[r] & {\rm Tor}^{\E^f}_*(k[T^2]^\vee,E^\bullet)\ar[d] \\
{\rm Tor}_*^{(\E^f)^e}(k[{\rm Hom}_{(\E^f)^{op}}],(E^\bullet)^\vee\boxtimes E^\bullet)\ar[r] & {\rm Tor}_*^{(\E^f)^e}(k[{\rm Hom}_{(\E^f)^{op}}],\sigma^* E^\bullet)\ar[r] & {\rm Tor}^{\E^f}_*(k[T^2]^\vee,E^\bullet)
}$$
dans lequel les flèches verticales sont construites comme précédemment et où l'on désigne par $\sigma$ le foncteur $(\E^f)^e\to\E^f\quad (W,V)\mapsto W^*\oplus V$. Le carré de droite commute car les flèches horizontales sont des isomorphismes d'adjonction entre deux foncteurs ($\sigma$ et $V\mapsto (V^*,V)$) invariants par l'auto-équivalence $(W,V)\mapsto (V^*,W^*)$ de $(\E^f)^e$. Le carré de gauche commute au signe $\epsilon^{ij}$ près lorsque l'on se restreint en degrés covariant $i$ et contravariant $j$, par définition de l'(anti)commutativité de $\E^\bullet$.

Il suffit alors de reprendre la démonstration de la proposition~\ref{preg} pour obtenir la conclusion.
\end{proof}

Le lemme suivant permet de déterminer explicitement, via les résultats de \cite{FFSS}, l'involution décrite dans l'énoncé précédent, lorsque $E^\bullet$ est un foncteur exponentiel gradué usuel. On y note par abus $Id\in {\rm Ob}\,\F$ le foncteur d'inclusion $\E^f\hookrightarrow\E$.

\begin{lm}\label{lmdfls}
L'involution sur ${\rm Ext}^*_\F(Id,Id)$ induite par l'auto-adjonction du foncteur $D$ et l'auto-dualité du foncteur $Id$ est triviale.
\end{lm}

\begin{proof} Franjou, Lannes et Schwartz ont déterminé dans l'article \cite{FLS} ${\rm Ext}^*_\F(Id,Id)$ : cet espace vectoriel gradué est de dimension $1$ en degré pair et nul en degré impair ; comme algèbre (pour le produit de composition), ${\rm Ext}^*_\F(Id,Id)$ est une algèbre symétrique (donc en particulier commutative) sur des générateurs $e_i$ de degré $2p^{i-1}$, pour $i>0$, quotientée par l'idéal des puissances $p$-ièmes. L'involution donnée par la dualité étant un (anti)morphisme d'algèbre graduée, il suffit de vérifier qu'elle préserve les $e_i$.

Pour cela, nous aurons besoin d'utiliser la catégorie $\mathcal{P}$ des foncteurs "polynomiaux stricts"  sur $k$ de Friedlander-Suslin (cf. \cite{FS}). Il existe un foncteur exact $\mathcal{P}\to\F$ compatible à la dualité ($\mathcal{P}$ est munie d'une dualité analogue à celle de $\F$) tel que $e_i$ est l'image d'un élément, traditionnellement encore noté de la même façon, appartenant à ${\rm Ext}_\mathcal{P}^{2p^{i-1}}(Id^{(i)},Id^{(i)})$, où $Id^{(i)}$ désigne le foncteur identité $Id$ tordu $i$ fois par le morphisme de Frobenius (cf. \cite{FS}, §\,$4$). De plus, le lemme~$4.2$ de \cite{FS} montre que $e_1$ est auto-dual. On utilise maintenant le corollaire~$5.9$ de \cite{FFSS} : il montre que l'image de $e_i$ par le morphisme ${\rm Ext}_\mathcal{P}^{2p^{i-1}}(Id^{(i)},Id^{(i)})\to {\rm Ext}_\mathcal{P}^{2p^{i-1}}(\Gamma^{p^{i-1}(1)},S^{p^{i-1}(1)})$ donné par la postcomposition par le morphisme de Frobenius (itéré et tordu) $Id^{(i)}\to S^{p^{i-1}(1)}$ et la précomposition par son dual (de sorte que ce morphisme est compatible à la dualité) envoie $e_i$ sur $e_1^{p^{i-1}}$, où le produit est cette fois relatif à l'algèbre (trigraduée) ${\rm Ext}_\mathcal{P}^*(\Gamma^{*(1)},S^{*(1)})$, dont la structure multiplicative, déduite des structures exponentielles duales de $\Gamma^*$ et $S^*$, est compatible à la dualité. Il s'ensuit que $e_i$ est auto-dual, d'où le lemme. 
\end{proof}

Nous rappelons maintenant le premier point du théorème~$6.3$ de \cite{FFSS}, dans la version bigraduée (moins précise que la version trigraduée originelle) qui nous intéresse, où l'assertion relative à la dualité se déduit du lemme~\ref{lmdfls} :

\begin{thm}[Franjou-Friedlander-Scorichenko-Suslin]\label{ffss1}
L'algèbre de Hopf bigraduée ${\rm Ext}^*_\F(\Gamma^\bullet,S^\bullet)$ est une algèbre symétrique $S(U)$, où $U$ est un espace vectoriel bigradué de dimension $2$ en bidegré $(2q^s m,q^s+1)$ pour tous entiers $m\geq 0$ et $s>0$, de dimension $1$ en bidegré $(2m,2)$ pour tout entier $m\geq 0$ et nul dans les autres bidegrés. Le premier degr\'e est le degr\'e cohomologique et le second le degr\'e interne.
 
 De surcroît, la dualité est égale à $S(t)$, où $t$ est une involution graduée de $U$ dont $1$ est valeur propre simple en chaque bidegré où $U$ est non nul.
\end{thm}

On en d\'eduit le th\'eor\`eme suivant :

\begin{thm}\label{thfc1}
L'algèbre de Hopf bigraduée ${\rm Ext}^*_\F(\Gamma^*,I\circ \Gamma^2)$ (resp. ${\rm Ext}^*_\F(\Gamma^*,I\circ\Lambda^2)$) est une algèbre symétrique $S(V_\Gamma)$ (resp. $S(W_\Gamma)$), où  $V_\Gamma$ (resp. $W_\Gamma$) est un espace vectoriel bigradué de dimension $1$ en bidegré $(2q^s m,q^s+1)$ pour tous entiers $m\geq 0$ et $s\geq 0$ (resp. $s>0$) et nul dans les autres bidegrés.
\end{thm}

\begin{proof}
Le théorème précédent (dont on conserve les notations) décrit l'algèbre de Hopf bigraduée ${\rm Ext}^*_\F(\Gamma^\bullet,S^\bullet)\simeq S(U)$, on a donc aussi ${\rm Ext}^*_\F(\Gamma^\bullet,I\circ T^2)\simeq S(U)$ par la proposition~\ref{lmhc}, ainsi que $h(\tau)\simeq S(t)$ par le lemme~\ref{echt}. On en déduit par les propositions~\ref{convef} et~\ref{convexp}.\,\ref{ece3}
$$h(1+\tau)\simeq S(1)*S(t)\simeq S(1+t)\qquad\qquad\text{puis}$$
$$h\Big(\frac{1+\tau}{2}\Big)=h\Big(\frac{p+1}{2}.(1+\tau)\Big)=h(1+\tau)^{*\frac{p+1}{2}}\simeq S(1+t)^{*\frac{p+1}{2}}=S\Big(\frac{p+1}{2}.(1+t)\Big)=S\Big(\frac{1+t}{2}\Big) ;$$
de même $h(\frac{1-\tau}{2})\simeq S(\frac{1-t}{2})$, ce qui donne le résultat.
\end{proof}

En procédant de façon similaire, utilisant cette fois la dernière assertion du théorème~$6.3$ de \cite{FFSS}, on obtient le résultat suivant :

\begin{thm}\label{thfc2}
L'algèbre de Hopf bigraduée ${\rm Ext}^*_\F(\Lambda^*,I\circ\Gamma^2)$ (resp. ${\rm Ext}^*_\F(\Lambda^*,I\circ\Lambda^2)$) est une algèbre à puissances divisées $\Gamma(V_\Lambda)$ (resp. $\Gamma(W_\Lambda)$), où $V_\Lambda$ (resp. $W_\Lambda$) est un espace vectoriel bigradué de dimension $1$ en bidegré $(2q^s m +q^s-1,q^s+1)$ pour tous entiers $m\geq 0$ et $s>0$ (resp. $s\geq 0$) et nul dans les autres bidegrés.
\end{thm}

Le corollaire \ref{crf-o} et les th\'eor\`emes \ref{thfc1} et \ref{thfc2} impliquent :
\begin{thm}
\begin{enumerate} 
\item 
L'alg\`ebre bigradu\'ee  
$$H^*(O_\infty;S^\bullet_\infty)$$
 est une algèbre symétrique $S(V_S)$, où  $V_S$ est un espace vectoriel bigradué de dimension $1$ en bidegré $(2q^s m,q^s+1)$ pour tous entiers $m\geq 0$ et $s\geq 0$ et nul dans les autres bidegrés.
 \item
 L'algèbre bigraduée 
 $$H^*(O_\infty;\Lambda^\bullet_\infty)$$
 est une algèbre à puissances divisées $\Gamma(V_\Lambda)$, où $V_\Lambda$ est un espace vectoriel bigradué de dimension $1$ en bidegré $(2q^s m +q^s-1,q^s+1)$ pour tous entiers $m\geq 0$ et $s>0$ et nul dans les autres bidegrés.
 \item 
L'alg\`ebre bigradu\'ee  
$$H^*(Sp_\infty;S^\bullet_\infty)$$
 est une algèbre symétrique $S(W_S)$, où  $W_S$ est un espace vectoriel bigradué de dimension $1$ en bidegré $(2q^s m,q^s+1)$ pour tous entiers $m\geq 0$ et $s> 0$ et nul dans les autres bidegrés.
 \item
 L'algèbre bigraduée 
 $$H^*(Sp_\infty;\Lambda^\bullet_\infty)$$
 est une algèbre à puissances divisées $\Gamma(W_\Lambda)$, où $W_\Lambda$ est un espace vectoriel bigradué de dimension $1$ en bidegré $(2q^s m +q^s-1,q^s+1)$ pour tous entiers $m\geq 0$ et $s \geq 0$ et nul dans les autres bidegrés.

\end{enumerate}
\end{thm}

\begin{thm}
\begin{enumerate} 
\item 
La cog\`ebre bigradu\'ee  
$$H_*(O_\infty;\Gamma^\bullet_\infty)$$
 est une cogèbre \`a puissances divis\'ees  $\Gamma(V_\Gamma)$, où  $V_\Gamma$ est un espace vectoriel bigradué de dimension $1$ en bidegré $(2q^s m,q^s+1)$ pour tous entiers $m\geq 0$ et $s\geq 0$ et nul dans les autres bidegrés.
 \item
 La cogèbre bigraduée 
 $$H_*(O_\infty;\Lambda^\bullet_\infty)$$
 est une cogèbre sym\'etrique $S(V_\Lambda)$, où $V_\Lambda$ est un espace vectoriel bigradué de dimension $1$ en bidegré $(2q^s m +q^s-1,q^s+1)$ pour tous entiers $m\geq 0$ et $s>0$ et nul dans les autres bidegrés.
 \item 
La cog\`ebre bigradu\'ee  
$$H_*(Sp_\infty;\Gamma^\bullet_\infty)$$
 est une cogèbre \`a puissances divis\'ees  $\Gamma(W_\Gamma)$, où  $W_\Gamma$ est un espace vectoriel bigradué de dimension $1$ en bidegré $(2q^s m,q^s+1)$ pour tous entiers $m\geq 0$ et $s> 0$ et nul dans les autres bidegrés.
 \item
 La cogèbre bigraduée 
 $$H_*(Sp_\infty;\Lambda^\bullet_\infty)$$
 est une cogèbre sym\'etrique $S(W_\Lambda)$, où $W_\Lambda$ est un espace vectoriel bigradué de dimension $1$ en bidegré $(2q^s m +q^s-1,q^s+1)$ pour tous entiers $m\geq 0$ et $s \geq 0$ et nul dans les autres bidegrés.

\end{enumerate}
\end{thm}

Par le corollaire~\ref{crf-os}, le théorème pr\'ec\'edent  calcule les espaces vectoriels $H_i(O_{n,n},\Gamma^j(\mathbb{F}_q^{2n}))$, $H_i(O_{n,n},\Lambda^j(\mathbb{F}_q^{2n}))$, $H_i(Sp_{2n},\Gamma^j(\mathbb{F}_q^{2n}))$ et $H_i(Sp_{2n},\Lambda^j(\mathbb{F}_q^{2n}))$ pour $n\geq 2i+j+6$.

\smallskip

Pour déterminer par la même méthode $H_*(O_\infty,S^\bullet_\infty)$ et $H_*(Sp_\infty,S^\bullet_\infty)$, on a besoin de connaître ${\rm Ext}^*_\F(S^\bullet,\Gamma^\bullet)$. Ce calcul, manquant dans \cite{FFSS}, a été effectué récemment par Cha{\l}upnik dans \cite{Chal}, corollaire~$4.6$ (toujours en transitant par la catégorie $\mathcal{P}$). On laisse au lecteur le soin d'écrire le résultat, un petit peu plus technique, obtenu.

Pour conclure ce paragraphe, mentionnons un corollaire frappant de ce théorème, dont les auteurs ignorent s'il peut être établi de manière plus directe.

\begin{cor}\label{cr-radj}
 Soient $n$ et $i$ deux entiers tels que $n\geq 2i+8$. Alors le $i$-ème groupe d'homologie du groupe $O_{n,n}$ (resp. $Sp_{2n}$) à coefficients dans sa représentation adjointe est nul.
\end{cor}

\begin{proof} Soit $(V,q)$ un espace quadratique de dimension finie non dégénéré. L'isomorphisme d'espaces vectoriels $\varphi : T^2(V)\xrightarrow{\simeq} {\rm End}\,V$ composé de l'isomorphisme canonique $T^2(V)\simeq {\rm Hom}_k (V^*,V)$ et de ${\rm Hom}_k(\phi,V)$, où $\phi : V\xrightarrow{\simeq}V^*$ est l'isomorphisme déduit de $q$, est $O(V,q)$-équivariant, où l'action est la restriction de l'action de $GL(V)$ donnée par la fonctorialité de $T^2$ à la source et la conjugaison au but (cf. remarque~\ref{rq-bif-AT}). De plus, $\varphi$ se restreint en un isomorphisme $O(V,q)$-équivariant entre $\Lambda^2(V)$ et la représentation adjointe de $O(V,q)$.
 
Par conséquent, l'annulation dans le cas orthogonal provient de celle de $H_*(O_\infty;\Lambda^2)$ contenue dans le théorème précédent et du corollaire~\ref{crf-os} (cf. remarque suivant le théorème précédent).

Le cas symplectique s'obtient de la même manière à partir de l'annulation de $H_*(Sp_\infty;S^2)$.
\end{proof}

\subsection{Un calcul en caractéristique $2$}

\begin{conv}
{\bf Dans ce paragraphe, $k$ est le corps à $2$ éléments $\FF$.}
\end{conv}

On pourrait procéder de manière analogue sur tout corps fini de caractéristique $2$, mais nous nous restreignons à ce cas car nous nous appuyons sur des résultats de Troesch qui n'ont été énoncés que sur des corps premiers.

Sur $\FF$, le foncteur $T^2$ n'est pas semi-simple, et $\Gamma^2$ n'est pas simple (il n'est pas non plus isomorphe à $S^2$) : on a des suites exactes courtes non scindées
$$0\to \Lambda^2\to\Gamma^2\to Id\to 0\qquad\text{et}\qquad 0\to\Gamma^2\to T^2\to\Lambda^2\to 0,$$
où la flèche $\Gamma^2\to Id$ ("Verschiebung") est duale du morphisme de Frobenius $Id\to S^2$.

Nous calculerons seulement les groupes ${\rm Ext}^*_\F(Id,I\circ\Gamma^2)$ et ${\rm Ext}^*_\F(Id,I\circ\Lambda^2)$, la détermination complète de ${\rm Ext}^*_\F(F,I\circ\Gamma^2)$ ou ${\rm Ext}^*_\F(F,I\circ\Lambda^2)$ semblant hors de portée lorsque $F$ est un foncteur polynomial de degré strictement supérieur à~$1$.

De fait, les foncteurs ${\rm Ext}^*_\F(Id,-)$ jouissent de la remarquable propriété suivante (qui est vraie sur tout corps fini premier), établie dans \cite{T1} (théorème~$3.1$) :

\begin{thm}[Troesch]\label{th-T}
Soient $C^\bullet$ un complexe exact d'objets de $\F$ sans terme constant et $F$ un foncteur analytique de $\F$. On a 
$${\rm Ext}^*_\F(Id,H^*(F\circ C^\bullet))=0.$$
\end{thm}
(Troesch énonce ce résultat seulement pour $F$ polynomial, mais le résultat général s'en déduit aussitôt par passage à la colimite filtrante, puisque $Id$ possède une résolution projective de type fini, d'après un résultat classique de Schwartz, établi par exemple dans \cite{FLS}, proposition~$10.1$.)

Par conséquent, les deux suites exactes précédentes induisent après postcomposition par $I$ (le caractère analytique de ce foncteur est un fait classique --- cf. par exemple \cite{FFPS}) des suites exactes longues :
$$\cdots\to {\rm Ext}^{i-1}(Id,I)\to {\rm Ext}^i(Id,I\circ\Lambda^2)\to {\rm Ext}^i(Id,I\circ\Gamma^2)\to {\rm Ext}^{i}(Id,I)\to\cdots$$
et
$$\cdots\to {\rm Ext}^i(Id,I\circ T^2)\to {\rm Ext}^i(Id,I\circ\Lambda^2)\to {\rm Ext}^{i+1}(Id,I\circ\Gamma^2)\to {\rm Ext}^{i+1}(Id,I\circ T^2)\to\cdots$$

Mais ${\rm Ext}^i(Id,I)$ est nul pour $i\neq 0$ et isomorphe à $\FF$ pour $i=0$, tandis que ${\rm Ext}^i(Id,I\circ T^2)=0$ pour tout $i$ d'après la proposition~\ref{preg}. On en déduit aussitôt :

\begin{thm}\label{thcal2}
Le $\FF$-espace vectoriel ${\rm Ext}^i_\F(Id,I\circ\Gamma^2)$ est de dimension $1$ pour $i\geq 2$ et nul sinon.

Le $\FF$-espace vectoriel ${\rm Ext}^i_\F(Id,I\circ\Lambda^2)$ est de dimension $1$ pour $i\geq 1$ et nul pour $i=0$.
\end{thm}

\begin{rem}
De façon similaire on obtient que ${\rm Ext}^i(Id,I\circ S^2)$ est de dimension $1$ pour tout $i\in\mathbb{N}$.
\end{rem}

On en déduit, par le corollaire~\ref{crf-os}, le résultat suivant. 

\begin{cor}\label{crc22}
Soient $n$ et $i$ des entiers naturels tels que $n\geq 2i+7$.
\begin{enumerate}
 \item L'espace vectoriel $H_i(O_{n,n}(\FF),\FF^{2n})$ est de dimension $i-1$ pour $i\geq 2$ et nul sinon.
 \item L'espace vectoriel $H_i(Sp_{2n}(\FF),\FF^{2n})$ est de dimension $1$ pour $i\geq 1$ et nul sinon.
\end{enumerate}
\end{cor}

\appendix

\section{Algèbre homologique dans les catégories de foncteurs}\label{appA}

Soit $\C$ une petite catégorie ; on rappelle que l'on note $\C-\mathbf{Mod}$ la catégorie des foncteurs de source $\C$ et de but la catégorie $\mathbf{Mod}_\kk$ et $\mathbf{Mod}-\C$ pour $\C^{op}-\mathbf{Mod}$. Nous rappelons dans cet appendice quelques résultats homologiques élémentaires et bien connus de cette catégorie, qui se comporte à plusieurs égards comme une catégorie de modules.

Cette catégorie est abélienne, l'exactitude se testant au but. Elle possède toutes limites et colimites, qui se calculent au but. (On renvoie, par exemple, à \cite{Gab} pour les généralités sur les catégories abéliennes.) C'est même une catégorie de Grothendieck : elle possède un générateur et les colimites filtrantes sont exactes. On peut préciser le premier point de la façon suivante : pour tout objet $i$ de $\C$, on note $P^\C_i$ le foncteur $\kk[{\rm Hom}_\C(i,-)]$ (où $\kk[E]$ désigne le $\kk$-module libre de base $E$). Le lemme de Yoneda procure un isomorphisme canonique ${\rm Hom}_\C(P^\C_i,F)\simeq F(i)$ (on note ici ${\rm Hom}_\C$ pour ${\rm Hom}_{\C-\mathbf{Mod}}$ ; des conventions analogues pour les groupes d'extensions seront utilisées), d'où l'on déduit que les $P^\C_i$ sont projectifs et engendrent $\C-\mathbf{Mod}$. On les appelle parfois {\em générateurs projectifs standard} de $\C-\mathbf{Mod}$.

\smallskip

Le {\em produit tensoriel au-dessus de $\C$} (cf. par exemple l'appendice C de \cite{Loday}) est le foncteur $\underset{\C}{\otimes} : \mathbf{Mod}-\C\times\C-\mathbf{Mod}\to\mathbf{Mod}_\kk$ donné ainsi : pour $F\in {\rm Ob}\,\mathbf{Mod}-\C$ et $G\in {\rm Ob}\,\C-\mathbf{Mod}$, $F\underset{\C}{\otimes}G$ est le quotient de $\underset{i\in {\rm Ob}\,\C}{\bigoplus} F(i)\otimes G(i)$ (on rappelle que les produits tensoriels non spécifiés sont pris sur $\kk$) par le sous-module engendré par les éléments du type $F(f)(x')\otimes y-x'\otimes G(f)(y)$ pour $i\xrightarrow{f}i'$ flèche de $\C$, $y\in G(i)$ et $x'\in F(i')$.

Le bifoncteur $\underset{\C}{\otimes}$ commute aux colimites en chaque variable. Il possède des propriétés de commutativité et d'associativité évidentes. De plus, on a des isomorphismes naturels
$$P^{\C^{op}}_i\underset{\C}{\otimes}G\simeq G(i)\qquad\text{et}\qquad F\underset{\C}{\otimes} P^\C_i\simeq F(i).$$
Une autre manière de caractériser le produit tensoriel est l'adjonction suivante :
\begin{equation}\label{homtens}
{\rm Hom}_\kk(F\underset{\C}{\otimes}G,V)\simeq {\rm Hom}_\C(G,\mathcal{H}om_\kk(F,V))\,,
\end{equation}
où, pour $V\in {\rm Ob}\,\mathbf{Mod}_\kk$ et $F\in {\rm Ob}\,\mathbf{Mod}-\C$, $\mathcal{H}om_\kk(F,V)$ désigne l'objet de $\C-\mathbf{Mod}$ donné par $C\mapsto {\rm Hom}_\kk(F(C),V)$.

Pour des raisons standard d'algèbre homologique, dériver à gauche le bifoncteur $\underset{\C}{\otimes}$ relativement à l'une ou l'autre des variables donne des résultats canoniquement isomorphes, notés ${\rm Tor}^\C_*$. (On peut aussi définir ces groupes comme l'homologie d'un complexe simplicial explicite --- cf. \cite{Loday}, appendice C, par exemple.)

Un foncteur $G\in {\rm Ob}\,\C-\mathbf{Mod}$ est dit {\em plat} si $-\underset{\C}{\otimes}G$ est exact. Tout foncteur projectif est plat, et les foncteurs plats sont stables par colimite filtrante. Comme dans le cas des modules, on peut remplacer, pour calculer les groupes de torsion sur $\C$, les résolutions projectives par des résolutions plates.
\smallskip

L'isomorphisme~(\ref{homtens}) montre que le produit tensoriel au-dessus de $\C$ est en quelque sorte dual du foncteur Hom sur $\C-\mathbf{Mod}$. Cela fonctionne particulièrement bien lorsque $\kk$ est un corps ; on note alors $F^*$ pour $\mathcal{H}om_\kk(F,\kk)$, la post-composition de $F$ par le foncteur de dualité des espaces vectoriels, noté $V\mapsto V^*$. L'exactitude de cette dualité donne alors un isomorphisme naturel gradué 
\begin{equation}\label{dualfct}
{\rm Tor}^\C(F,G)^*\simeq {\rm Ext}_\C(G,F^*).
\end{equation}

\smallskip

L'{\em homologie} de la catégorie $\C$ à coefficients dans un foncteur $F\in {\rm Ob}\,\C-\mathbf{Mod}$ est par définition le $\kk$-module gradué ${\rm Tor}^\C_*(\kk,F)$ (où $\kk$ désigne le foncteur de $\mathbf{Mod}-\C$ constant en $\kk$), qui est noté $H_*(\C;F)$. On remarque que $H_0(\C;F)$ s'identifie canoniquement à la colimite du foncteur $F$. Lorsque $\C$ est la catégorie à un objet associée à un groupe (ou plus généralement un monoïde) $G$, $F$ se réduit à un $\kk[G]$-module et l'on retrouve la notion habituelle d'homologie du groupe $G$ (ou plus généralement de groupes de torsion sur $G$). Par ailleurs, le module gradué $H_*(\C;\kk)$ à coefficients dans le foncteur constant s'identifie à l'homologie singulière du nerf de~$\C$.

\smallskip

Nous aurons besoin du résultat suivant déduit de la th\'eorie de l'obstruction des complexes de cha\^ines exposée par Dold dans \cite{Dold} ; ce n'est qu'une transposition au cas des cat\'egories de foncteurs de type $\C-\mathbf{Mod}$ de la proposition $1.6$ de l'article \cite{P-hodge} de Pirashvili sur les $\Gamma$-modules.

\begin{pr}[Pirashvili]\label{obstruction}
Soient $F\in {\rm Ob}\,\C-\mathbf{Mod}$ et $C_\bullet$ un complexe de cha\^ines $\mathbb{N}$-gradu\'e d'objets projectifs de $\mathbf{Mod}-\C$.
\begin{enumerate}
\item Il existe une suite spectrale du premier quadrant
$$E^2_{p,q}={\rm Tor}_p^{\C}(H_q(C_\bullet),F) \Rightarrow H_{p+q}(C_\bullet \underset{\C}{\otimes}  F).$$
\item
Supposons que 
\begin{equation}\label{obstruction1}
{\rm Ext}^{m-n+1}_{\C^{op}}(H_n(C_\bullet), H_m(C_\bullet))=0 \quad\text{pour}\quad n<m.
\end{equation}
Alors la suite spectrale s'effondre au terme $E^2$; de plus, il existe une d\'ecomposition :
\begin{equation}\label{obstruction2}
H_n(C_\bullet \underset{\C}{\otimes} F) \simeq \bigoplus_{p+q=n} {\rm Tor}^\C_p(H_q(C_\bullet),F)
\end{equation}
naturelle en $F$.
\end{enumerate}
\end{pr}
(La suite spectrale en question s'obtient en considérant le complexe double produit tensoriel au-dessus de $\C$ de $C_{\bullet}$ et d'une résolution projective de $F$.)

\smallskip

Une généralisation de l'homologie et des groupes de torsion (modulo une hypothèse de projectivité sur les coefficients) sur la catégorie $\C$ est celle d'homologie des {\em bi}foncteurs, i.e. des objets de $\C^{op}\times\C_\mathbf{Mod}$. Comme cette notion correspond, dans le cas d'une catégorie à un objet, à celle d'homologie de Hochschild de l'algèbre du monoïde associé, nous noterons $HH_*$ cette théorie homologique. Si $B$ est un tel bifoncteur, on définit $HH_0(\C;B)$ comme le quotient de $\underset{i\in {\rm Ob}\,\C}{\bigoplus}B(i,i)$ par le sous-module engendré par les éléments du type $B(f,id_i)(x)-B(id_j,f)(x)$ pour $i$, $j$ objets de $\C$, $x\in B(j,i)$ et $f\in {\rm Hom}_\C(i,j)$. Le foncteur $HH_0(\C,-)$ est exact à droite ; les $\kk$-modules $HH_i(\C;B)$ sont par définition l'évaluation en $B$ de ses foncteurs dérivés. Si $B$ est un produit tensoriel extérieur $F\boxtimes G$ avec $F\in {\rm Ob}\,\mathbf{Mod}-\C$ et $G\in {\rm Ob}\,\C-\mathbf{Mod}$, i.e. $B(i,j)=F(i)\otimes G(j)$, et que l'un des foncteurs $F$ et $G$ a pour valeurs des $\kk$-modules projectifs, on a $HH_*(\C;B)={\rm Tor}^\C_*(F,G)$ (c'est vrai sans hypothèse en degré $0$).

Une autre description de $HH_*$ est donnée par l'isomorphisme gradué naturel suivant :
\begin{equation}\label{hochtor}
HH_*(\C;B)\simeq {\rm Tor}_*^{\C^{op}\times\C}(\kk[{\rm Hom}_{\C^{op}}(-,-)],B)
\end{equation}
(cf. \cite{Loday}, appendice C, par exemple).

\smallskip

L'homologie de $\C$ à coefficients dans un bifoncteur (et de même les groupes de torsion sur $\C$) possède une fonctorialité en $\C$ : si $Q : \D\to\C$ est un foncteur entre petites catégories et $B$ un objet de ${\rm Ob}\,\C^{op}\times\C-\mathbf{Mod}$, on dispose d'un morphisme canonique $HH_*(\D;Q^*(B))\to HH_*(\C;B)$, où l'étoile indique la précomposition et où l'on note par abus encore $Q$ pour $Q^{op}\times Q : \D^{op}\times\D\to\C^{op}\times\C$. En degré $0$, il se définit en constatant que le morphisme évident
$$\bigoplus_{j\in {\rm Ob}\,\D}B(Q(j),Q(j))\to\bigoplus_{i\in {\rm Ob}\,\C}B(i,i)$$
induit un morphisme $HH_0(\D;Q^*(B))\to HH_0(\C;B)$, qu'on prolonge ensuite en un morphisme de foncteurs homologiques (c'est possible par exactitude de $Q^*$). C'est un isomorphisme si $Q$ est une équivalence de catégories.
\smallskip

Il existe un analogue de la notion de foncteurs adjoints pour le produit tensoriel au-dessus d'une petite catégorie. Le cas qui nous intéresse est celui correspondant aux extensions de Kan, qu'on peut traiter de façon analogue :

\begin{pr}\label{kan-t}
Soient $\C$ et $\D$ des petites catégories et $Q : \D\to\C$ un foncteur. Il existe un foncteur $Q_! : \mathbf{Mod}-\D\to\mathbf{Mod}-\C$, unique à isomorphisme canonique près, tel qu'existe un isomorphisme naturel
$$F\underset{\D}{\otimes} Q^*(G)\simeq Q_!(F)\underset{\C}{\otimes} G$$
pour $F\in {\rm Ob}\,\mathbf{Mod}-\D$ et $G\in {\rm Ob}\,\C-\mathbf{Mod}$ (où $Q^*$ désigne la précomposition par $Q$). On peut définir $Q_!$ par 
\begin{equation}\label{kanev}
Q_!(F)(a)=F\underset{\D}{\otimes} Q^*(P^\C_a).
\end{equation}

Le foncteur $Q_!$ est exact à droite ; il commute même à toutes les colimites. Ses foncteurs dérivés à gauche sont donnés par
\begin{equation}\label{kander}
\mathbb{L}_iQ_!(F)(a)={\rm Tor}^\D_i(F,Q^*(P^\C_a)).
\end{equation}

De plus, le foncteur $Q_!$ est adjoint à gauche au foncteur $(Q^{op})^* : \mathbf{Mod}-\C\to\mathbf{Mod}-\D$ ; son effet sur les projectifs standard est donné par un isomorphisme naturel
$$Q_!(P^{\D^{op}}_b)\simeq P^{\C^{op}}_{Q(b)}.$$
\end{pr}

La démonstration, qui consiste en une variation sur le lemme de Yoneda, est laissée au lecteur.

\smallskip

Nous terminons ces rappels par un résultat plus particulier qui intervient dans la section~\ref{sid}.

\begin{pr}\label{torcomod}
Soient $X$ un foncteur de $\C^{op}$ vers la catégorie des ensembles et $\C_X$ la catégorie dont les objets sont les couples $(i,x)$ formés d'un objet $i$ de $\C$ et d'un élément $x$ de $X(i)$, les morphismes $(i,x)\to (j,y)$ étant les morphismes $f : i\to j$ de $\C$ tels que $X(f)(y)=x$. On note $U : \C_X\to\C$ le foncteur d'oubli (donné sur les objets par $(i,x)\mapsto i$), et $\Omega_X : \mathbf{Mod}-\C_X\to\mathbf{Mod}-\C$ le foncteur donné par
$$\Omega_X(F)(i)=\bigoplus_{x\in X(i)}F(i,x),$$
le morphisme $\Omega_X(F)(f) : \Omega_X(F)(i)\to\Omega_X(F)(j)$ induit par un morphisme $f : i\to j$ de $\C^{op}$ ayant pour composante $F(i,x)\to F(j,y)$ l'application induite par $F$ si $X(f)(i)=j$ et $0$ sinon.

Il existe un isomorphisme gradué
$${\rm Tor}^{\C_X}(F,U^*(G))\simeq {\rm Tor}^\C(\Omega_X(F),G)$$
naturel en les objets $F$ de $\mathbf{Mod}-\C_X$ et $G$ de $\C-\mathbf{Mod}$.
\end{pr}

Ce résultat élémentaire est la traduction en terme de foncteur Tor des adjonctions étudiées dans le §\,3.1 de \cite{Dja}. Le cas du degré~$0$ se déduit de la formule~(\ref{kanev}) ; l'exactitude du foncteur $\Omega_X$ permet de propager l'isomorphisme en tout degré homologique.

\smallskip

Le cas des coefficients constants est ce que Loday nomme {\em lemme de Shapiro pour l'homologie des petites catégories} (\cite{Loday}, appendice C.12).

\section{Foncteurs exponentiels}\label{apex}

Dans cet appendice, on s'intéresse uniquement à la catégorie $\F(k)=\E^f(k)-\mathbf{Mod}$ des foncteurs de source $\E^f(k)$ et de but $\E(k)$. On notera également $\F_2(k)=\E^f(k)\times\E^f(k)-\mathbf{Mod}$.

Les notions que nous rappelons sont classiques : on pourra se référer à \cite{FFSS}.

\smallskip

Notons $\pi : \E^f(k)\times\E^f(k)\to\E^f(k)$ le foncteur de somme directe et $\delta : \E^f(k)\to\E^f(k)\times\E^f(k)$ le foncteur d'inclusion diagonale. Chacun d'entre eux est adjoint à droite et à gauche à l'autre. Il en résulte que la même propriété vaut pour les foncteurs de précomposition $\pi^* : \F(k)\to\F_2(k)$ et $\delta^* : \F_2(k)\to\F(k)$. Cette observation simple s'avère tout-à-fait efficace pour traiter des foncteurs possédant la propriété suivante.

\begin{defi}\begin{enumerate}
\item Un {\em foncteur exponentiel} de $\F(k)$ est un objet $E$ de $\F(k)$ muni d'un isomorphisme $E(0)\simeq k$ et d'un isomorphisme entre les objets $\pi^*(E)$ et $E\boxtimes E$ de $\F_2(k)$ (i.e. d'un isomorphisme $E(U\oplus V)\simeq E(U)\otimes E(V)$ naturel en les objets $U$ et $V$ de $\E^f_k$). 
\item Un {\em foncteur exponentiel gradué} de $\F(k)$ est une suite $E^\bullet=(E^n)_{n\in\mathbb{N}}$ d'objets de $\F(k)$ prenant des valeurs de dimension finie muni d'un isomorphisme $E^0\simeq k$ (foncteur constant) et d'un isomorphisme gradué $\pi^*(E^\bullet)\simeq E^\bullet\boxtimes E^\bullet$ (i.e. pour tout $n\in\mathbb{N}$, d'un isomorphisme $\pi^*(E^n)\simeq\underset{i+j=n}{\bigoplus}E^i\boxtimes E^j$).
\item Soit $E$ un foncteur exponentiel, éventuellement gradué, de $\F(k)$. Si la structure exponentielle de $E$ est compatible aux isomorphismes d'associativité de la somme directe et du produit tensoriel\,\footnote{Cela signifie que le diagramme évident d'isomorphismes
$$\xymatrix{E((U\oplus V)\oplus W)\ar[r]\ar[d] & E(U\oplus (V\oplus W))\ar[d] \\
E(U\oplus V)\otimes E(W)\ar[d] & E(U)\otimes E(V\oplus W)\ar[d] \\
(E(U)\otimes E(V))\otimes E(W)\ar[r] & E(U)\otimes (E(V)\otimes E(W))
}$$
commute pour tous $U, V, W\in {\rm Ob}\,\E^f(k)$.}, on dit que $E$ est un {\em foncteur exponentiel de Hopf}.
\end{enumerate}
\end{defi}

Par exemple, les foncteurs projectifs standard de $\F(k)$ ou le foncteur injectif $I$ (cf. notation~\ref{not-inji}) sont des foncteurs exponentiels de Hopf. Les foncteurs puissance symétrique $S^*$, puissance divisée $\Gamma^*$ (dual gradué du précédent) et puissance extérieure $\Lambda^*$ sont exponentiels gradués de Hopf. (Cela provient directement des propriétés universelles qui les caractérisent.) 
De plus, $S^*$ et $\Gamma^*$ sont commutatifs et $\Lambda^*$ est anticommutatif (au sens gradué) en un sens évident, qui est précisé dans \cite{FFSS}, page~$675$.

Noter que les composantes d'un foncteur exponentiel gradué sont des foncteurs polynomiaux.

\begin{rem}\label{tensexp} Le foncteur gradué puissance tensorielle $T^*$ n'est pas exponentiel, mais il possède une propriété analogue dont on peut se servir dans les calculs d'une façon similaire à la propriété exponentielle : la formule du binôme se traduit par un isomorphisme
$$\pi^*(T^n)\simeq\bigoplus_{i+j=n}(T^i\boxtimes T^j)\underset{\mathfrak{S}_i\times\mathfrak{S}_j}{\otimes} k[\mathfrak{S}_n],$$
où le groupe symétrique $\mathfrak{S}_l$ agit par permutation des facteurs du produit tensoriel et le groupe $\mathfrak{S}_i\times\mathfrak{S}_j$ est plongé de manière usuelle dans $\mathfrak{S}_{i+j}$.
\end{rem}

La proposition suivante est une partie du théorème~$1.7$ de \cite{FFSS} ; sa démonstration repose essentiellement sur l'adjonction entre les foncteurs $\pi^*$ et $\delta^*$ évoquée plus haut et son prolongement aux groupes d'extensions (qui découle de leur exactitude).

\begin{pr}\label{expext} Soient $E^*$ un foncteur exponentiel gradué et $F$, $G$ deux foncteurs de $\F(k)$. Il existe, pour tout $n\in\mathbb{N}$, un isomorphisme gradué
$${\rm Ext}^*_{\F(k)}(E^n,F\otimes G)\simeq\bigoplus_{i+j=n} {\rm Ext}^*_{\F(k)}(E^i,F)\otimes {\rm Ext}^*_{\F(k)}(E^j,G)$$
naturel en $F$ et $G$.

(Autrement dit, on dispose d'un isomorphisme bigradué ${\rm Ext}^*_{\F(k)}(E^\bullet,F\otimes G)\simeq {\rm Ext}^*_{\F(k)}(E^\bullet,F)\otimes {\rm Ext}^*_{\F(k)}(E^\bullet,G)$.)
\end{pr}

Si $E$ est un foncteur exponentiel, éventuellement gradué, on dispose de morphismes
$k\simeq E(0)\to E$, appelé {\em unité}, $E\to E(0)\simeq k$, appelé {\em coünité}, $E\otimes E=\delta^*(E\boxtimes E)\simeq\delta^*\pi^*(E)\to E$ (induit par la coünité de l'adjonction), appelé {\em produit}, $E\to E\otimes E=\delta^*(E\boxtimes E)\simeq\delta^*\pi^*(E)$, appelé {\em coproduit}, et $E(-id) : E\to E$, appelé {\em antipode}. On vérifie facilement (cf. \cite{FFSS}) que ces applications font de $E$ une {\em algèbre de Hopf} (graduée connexe si $E$ est gradué) dans la catégorie monoïdale symétrique $(\F(k),\otimes,k)$ lorsque $E$ est un foncteur exponentiel de Hopf (c'est en fait une caractérisation des foncteurs exponentiels de Hopf).

Ce type de structure est utile pour définir des {\em produits de convolution}. Supposons en effet que $E^\bullet$ est un foncteur exponentiel de Hopf et $F$ un foncteur de $\F(k)$ muni d'une structure d'algèbre. Alors ${\rm Ext}^*_{\F(k)}(E^\bullet,F)$ est une algèbre bigraduée pour le produit, dit de convolution,
$$* : {\rm Ext}^*_{\F(k)}(E^\bullet,F)\otimes {\rm Ext}^*_{\F(k)}(E^\bullet,F)\simeq {\rm Ext}^*_{\F(k)}(E^\bullet,F\otimes F)\to {\rm Ext}^*_{\F(k)}(E^\bullet,F),$$
où l'isomorphisme est donné par la proposition~\ref{expext} et la seconde flèche est induite par le produit de $F$ ; l'unité de ce produit est donnée par l'élément de ${\rm Hom}_{\F(k)}(E^0,F)\simeq F(0)$ unité de l'algèbre $F$.

Si $F$ est un foncteur muni d'une structure de cogèbre, on définit de même un coproduit sur ${\rm Ext}^*_{\F(k)}(E^\bullet,F)$ comme le morphisme
$${\rm Ext}^*_{\F(k)}(E^\bullet,F)\to {\rm Ext}^*_{\F(k)}(E^\bullet,F\otimes F)\simeq {\rm Ext}^*_{\F(k)}(E^\bullet,F)\otimes {\rm Ext}^*_{\F(k)}(E^\bullet,F)$$
où la première flèche provient du coproduit de $F$. Lorsque $F$ est muni d'une structure d'algèbre de Hopf dans $\F(k)$, on obtient sur ${\rm Ext}^*_{\F(k)}(E^\bullet,F)$ une structure de $k$-algèbre de Hopf bigraduée (l'antipode étant induite par celle de $F$), connexe si $F(0)=k$.

On utilise, au paragraphe~\ref{pci}, les deux propositions suivantes où intervient la convolution.

\begin{pr}\label{convef} Soit $E^\bullet$ un foncteur exponentiel de Hopf gradué, commutatif ou anticommutatif.
\begin{enumerate}
\item Pour tout objet $V$ de $\E^f(k)$, $E^\bullet(V)$ est naturellement une $k$-algèbre de Hopf graduée.
\item Si $f : V\to W$ est une flèche de $\E^f(k)$, alors $E^\bullet(f) : E^\bullet(V)\to E^\bullet(W)$ est un morphisme de $k$-algèbres de Hopf graduées.

Si $g : V\to W$ est une autre flèche de $\E^f(k)$, on a $E^\bullet(f+g)=E^\bullet(f) * E^\bullet(g)$.
\end{enumerate}
\end{pr}

Cette conséquence facile des définitions est laissée au lecteur.

\begin{rem}
 La nécessité d'une hypothèse d'(anti)commutativité est omise dans l'article \cite{FFSS} ; elle a été remarquée par Touzé dans \cite{Touze} (cf. sa remarque~5.6). On renvoie d'une manière générale le lecteur à la section 5 de \cite{Touze} où toutes les questions de compatibilité (et parfois de regraduation) nécessaires dans les énoncés relatifs aux structures exponentielles sont traitées avec soin.
\end{rem}

\begin{pr}\label{convexp}
Soient $E^\bullet$ un foncteur exponentiel gradué de Hopf commutatif ou anticommutatif de $\F(k)$ et $F$ un foncteur muni d'une structure d'algèbre (resp. d'algèbre de Hopf) dans la catégorie monoïdale symétrique $(\F(k),\otimes,k)$.
\begin{enumerate}
\item Si $A$ est un foncteur sans terme constant (i.e. $A(0)=0$) de $\F(k)$, la structure d'algèbre (resp. d'algèbre de Hopf) de $F$ induit une structure d'algèbre (resp. d'algèbre de Hopf) sur $F\circ A$.
\item Si $u : A\to B$ est un morphisme entre foncteurs sans terme constant de $\F(k)$, ${\rm Ext}^*(E^\bullet,F(u)) : {\rm Ext}^*(E^\bullet,F\circ A)\to {\rm Ext}^*(E^\bullet,F\circ B)$ est un morphisme de $k$-algèbres (resp. de $k$-algèbres de Hopf) bigraduées. Ce morphisme sera noté $h(u)$.
\item\label{ece3} Si $v : A\to B$ est un autre tel morphisme et que $F$ est un foncteur exponentiel de Hopf commutatif ou anticommutatif, on a $h(u+v)=h(u)*h(v)$.
\end{enumerate}
\end{pr}

\begin{proof}
Les deux premières assertions sont immédiates. Pour la dernière, on écrit $u+v$ comme la composée $A\hookrightarrow A\oplus A\xrightarrow{u\oplus v}B\oplus B\twoheadrightarrow B$, où le premier morphisme est la diagonale et le dernier la somme, puis on considère le diagramme commutatif suivant
$$\xymatrix{{\rm Ext}^*_{\F(k)}(E^\bullet,F_* A)\ar@<-5ex>[d]\ar[dr] & \\
{\rm Ext}^*_{\F(k)}(E^\bullet,F_* (A\oplus A))\simeq {\rm Ext}^*_{\F(k)}(E^\bullet,F_* A\otimes F_* A)\ar@<-5ex>[d]_{h(u\oplus v)}\ar[r]_-\simeq & {\rm Ext}^*_{\F(k)}(E^\bullet,F_* A)\otimes {\rm Ext}^*_{\F(k)}(E^\bullet,F_*A)\ar[d]_{h(u)\otimes h(v)} \\
{\rm Ext}^*_{\F(k)}(E^\bullet,F_*(B\oplus B))\simeq {\rm Ext}^*_{\F(k)}(E^\bullet,F_* B\otimes F_* B)\ar@<-5ex>[d]\ar[r]^-\simeq & {\rm Ext}^*_{\F(k)}(E^\bullet,F_* B)\otimes {\rm Ext}^*_{\F(k)}(E^\bullet,F_* B)\ar[dl] \\
{\rm Ext}^*_{\F(k)}(E^\bullet,F_*B)  & & 
}$$
où l'on note pour alléger $F_*$ la postcomposition par $F$ et dans lequel les flèches verticales ou obliques non spécifiées sont induites par coproduit ou produit : le chemin vertical gauche fournit $h(u+v)$ et celui de droite $h(u)*h(v)$.
\end{proof}

\section{Inversion de Möbius et équivalences de Morita}\label{amm}

\input{appc}

\section{Quelques résultats d'annulation en homologie des foncteurs}\label{sect-cof}

\input{apcof}

\section{Les résultats de Betley sur les groupes symétriques revisités}\label{apsym}

On conserve les notations de l'appendice précédent ; notre anneau de base est $\kk=\mathbb{Z}$. On note également $U : \Theta\to\Gamma$ le foncteur composé de l'inclusion de $\Theta$ dans la catégorie des ensembles finis et de $(-)_+$. On dispose donc de $U_! : \mathbf{Mod}-\Theta\to\mathbf{Mod}-\Gamma$ par la notation introduite \`a la proposition~\ref{kan-t}.

On considère aussi la catégorie $\Sigma$ des ensembles finis avec bijections ; les foncteurs d'inclusion $J^\Omega : \Sigma\to\Omega$ et $J^\Theta : \Sigma^{op}\to\Theta^{op}$ donnent lieu à $J^\Omega_! : \mathbf{Mod}-\Sigma\to\mathbf{Mod}-\Omega$ et $J^\Theta_! : \Sigma-\mathbf{Mod}\to\Theta-\mathbf{Mod}$.

\begin{lm}
Pour $E$ un ensemble fini, on a
$$J^\Theta_!(F)(E)=\bigoplus_{E'\subset E} F(E').$$
\end{lm}
(On laisse au lecteur, dans cet énoncé comme dans la suite de l'appendice, le soin de préciser la fonctorialité en $E$, qui est analogue à plusieurs cas déjà rencontrés --- cf. la proposition~\ref{torcomod}.)

\begin{proof}
Cela découle de l'isomorphisme canonique
$$(J^\Theta)^*(P^{\Theta^{op}}_E)=\mathbb{Z}[{\rm Hom}_\Theta(-,E)]\simeq\bigoplus_{E'\subset E}\mathbb{Z}[{\rm Hom}_\Sigma(-,E')]=\bigoplus_{E'\subset E}P^{\Sigma^{op}}_{E'}.$$
\end{proof}

\begin{pr}\label{prcdk}
La composée $\mathbf{Mod}-\Theta\xrightarrow{U_!}\mathbf{Mod}-\Gamma\xrightarrow{cr}\mathbf{Mod}-\Omega$ est isomorphe à $\mathbf{Mod}-\Theta\xrightarrow{(J^\Theta)^*}\mathbf{Mod}-\Sigma\xrightarrow{J^\Omega_!}\mathbf{Mod}-\Omega$.

De surcroît, le foncteur $U_!$ est exact.
\end{pr}

\begin{proof}
Par un argument d'adjonction (ou plutôt sa variante en terme de produit tensoriel), la première assertion équivaut à l'existence d'un isomorphisme $U^*\circ i_!\simeq J^\Theta_!\circ (J^\Omega)^*$, qui résulte du lemme précédent.

Pour la deuxième assertion, il suffit de voir que $J^\Omega_!$ est exact, puisque $cr$ est une équivalence de catégories ; cela provient d'une description explicite analogue à celle du lemme précédent (faisant intervenir les quotients au lieu des sous-ensembles d'un ensemble fini).
\end{proof}

La suite spectrale du théorème suivant est équivalente à la conjonction des théorèmes~$1.23$ et~$3.2$ de \cite{Bet-sym} (on rappelle que la naturalité du scindement suggérée par le théorème~$1.23$ semble incorrecte). On y note $F_\infty$ pour $(U^*F)_\infty$, et $cr_i(F)$ pour $cr(F)(\{1,\dots,i\})$.

\begin{thm}[Betley]
Pour tout $F\in {\rm Ob}\,\Gamma-\mathbf{Mod}$, il existe des isomorphismes naturels
$$H_n(\mathfrak{S}_\infty;F_\infty)\simeq\bigoplus_{i\in\mathbb{N}}H_n(\mathfrak{S}_\infty\times\mathfrak{S}_i;cr_i(F))\simeq\underset{i\in\mathbb{N}}{\bigoplus_{p+q=n}}{\rm Tor}_p^{\mathfrak{S}_i}(H_q(\mathfrak{S}_\infty;\mathbb{Z}),cr_i(F))$$
(où l'action de $\mathfrak{S}_i$ sur $H_q(\mathfrak{S}_\infty;\mathbb{Z})$ est triviale), ainsi qu'une suite spectrale naturelle
$$E^2_{p,q}=H_p(\mathfrak{S}_\infty;H_q(\mathfrak{S}_i;cr_i(F)))\Rightarrow H_{p+q}(\mathfrak{S}_\infty;F_\infty)$$
(où l'action de $\mathfrak{S}_\infty$ est triviale) qui s'effondre à la deuxième page et procure un isomorphisme non naturel
$$H_n(\mathfrak{S}_\infty;F_\infty)\simeq\underset{i\in\mathbb{N}}{\bigoplus_{p+q=n}}H_p(\mathfrak{S}_\infty;H_q(\mathfrak{S}_i;cr_i(F))).$$
\end{thm}

\begin{proof}Les propositions~\ref{prcdk} et~\ref{prdk} procurent des isomorphismes : 
$${\rm Tor}^\Theta_*(G,U^*(F))\simeq {\rm Tor}^\Gamma(U_!(G),F)\simeq {\rm Tor}^\Omega_*(cr\circ U_!(G),cr(F))\simeq {\rm Tor}^\Omega_*(J_!^\Omega\circ (J^\Theta)^*(G),cr(F))$$
$$\simeq {\rm Tor}^\Sigma_*((J^\Theta)^*(G),(J^\Omega)^*\circ cr(F))\simeq\bigoplus_{i\in\mathbb{N}}{\rm Tor}_*^{\mathfrak{S}_i}((J^\Theta)^*(G),cr_i(F)).$$

On en déduit les isomorphismes annoncés à l'aide des propositions~\ref{fdo} et~\ref{arret2}. La suite spectrale et son effondrement proviennent de la proposition~\ref{eff-kun}.
\end{proof}

\begin{rem}
Notre méthode ne diffère pas fondamentalement de celle employée par Betley pour établir ce théorème, si l'on excepte le remplacement de la $K$-théorie stable généralisée par l'homologie de la catégorie~$\Theta$.
\end{rem}

\section{Un aperçu du cas des modules d'après Betley et Scorichenko}\label{appBS}

Dans cet appendice, on se donne un anneau (unitaire) $A$, non nécessairement commutatif ; on note $\F(A)$ pour $\mathbb{P}(A)-\mathbf{Mod}$, où $\mathbb{P}(A)$ est la catégorie des $A$-modules à gauche projectifs de type fini.

Comme dans l'exemple~\ref{exfond}\,.\ref{ex-sco}, $\mathbb{M}(A)$ désigne la sous-catégorie de $\mathbb{P}(A)$ avec les mêmes objets et les injections scindées comme morphismes. Nous aurons aussi besoin de la sous-catégorie $\mathbb{E}(A)$ ayant les mêmes objets que $\mathbb{P}(A)$ et les épimorphismes pour flèches. La catégorie $\mathbb{E}(A)^{op}$ est équivalente à $\mathbb{M}(A^{op})$ via le foncteur de dualité ${\rm Hom}_A(-,A)$.

 On utilise également les notations $L_q$ et $St$ introduites respectivement  en~\ref{not-fdl} et en~\ref{nost2}, le triplet $(\C,S,G)$ sous-jacent étant celui de l'exemple~\ref{exfond}.\,\ref{ex-sco}.

\begin{lm}\label{lmsts}
Pour tout objet $V$ de $\mathbb{P}(A)$, il existe une suite spectrale du premier quadrant
$$E^2_{p,q}=H_p(GL_\infty(A);H_q({\rm Hom}_A(-,V))_\infty)\Rightarrow L_{p+q}(V)$$
où ${\rm Hom}_A(-,V)$ est vu comme foncteur de $\mathbb{E}(A)^{op}$ vers $\mathbf{Ab}$.
\end{lm}

\begin{proof}
On note d'abord qu'on a $H_q({\rm Hom}_A(-,V))_\infty\simeq H_q({\rm Hom}_A(-,V)_\infty)$ canoniquement puisqu'homologie et colimites filtrantes commutent.

Si $E$ est un $A$-module projectif de type fini, le stabilisateur de $V\hookrightarrow V\oplus E$ sous l'action de $GL(V\oplus E)$ s'identifie au produit semi-direct ${\rm Hom}_A(E,V)\rtimes GL(E)$. On en déduit aisément l'énoncé, en utilisant la suite spectrale de Lyndon-Hochschild-Serre.
\end{proof}

Noter qu'il n'y a a priori {\em aucune fonctorialité en $V$} sur cette suite spectrale : l'aboutissement est "naturellement" fonctoriel contravariant sur les monomorphismes scindés, tandis que le terme $E^2$ est "naturellement" fonctoriel covariant sur toutes les applications $A$-linéaires ! En particulier, on ne semble pas disposer de description simple du morphisme canonique $L_n(V)\simeq H_n(St(V))\to L_n(0)\simeq H_n(GL_\infty(A))$ (qui induit par l'inclusion de $St(V)$ dans $GL_\infty(A)$).

On a néanmoins le résultat suivant : 
\begin{lm}\label{rqfonct}
Soit $n\in\mathbb{N}$. Supposons que, dans la suite spectrale précédente, le terme $E^2_{p,q}$ soit nul pour $p<n$ et $q>0$. Alors le foncteur $L_n$ est constant. 
\end{lm}
 
 \begin{proof}
 Cette assertion provient des deux observations suivantes :
\begin{itemize}
\item l'épimorphisme de groupes $\pi : St(V)\twoheadrightarrow GL_\infty(A)$ donné par la description précédente comme produit semi-direct induit, sous l'hypothèse d'annulation de l'énoncé, un isomorphisme en homologie de degré au plus $n$, puisque $H_*(\pi)$ est le "coin" $L_*(V)\to E^2_{*,0}$ (propriété générale de la suite spectrale de Lyndon-Hochschild-Serre) ;
\item la composée du monomorphisme $GL_\infty(A)\hookrightarrow St(V)$ (toujours donné par la description de $St(V)$ de la démonstration du lemme précédent) avec l'inclusion $St(V)\hookrightarrow GL_\infty(A)$ induit un isomorphisme en homologie (cf. la démonstration de la proposition~\ref{fdo}), et sa composée avec l'épimorphisme $\pi$ est l'identité.
\end{itemize}
\end{proof}

\begin{thm}[Betley]\label{thfb}
Si $F$ est un foncteur polynomial de $\F(A)$ tel que $F(0)=0$, alors $H_*(GL_\infty(A),F_\infty)=0$.
\end{thm}

Cette assertion est le théorème~4.2 de \cite{Bet2}, où l'anneau $A$ est supposé commutatif, mais les arguments de Betley ne semblent en fait pas réellement requérir cette hypothèse (Betley suppose aussi $\kk=\mathbb{Z}$, mais ce n'est pas restrictif).

\begin{proof}
On démontre  l'assertion par récurrence sur le degré homologique. Si elle est vraie en degré $<n$, on a $E^2_{p,q}=0$ dans la suite spectrale du lemme~\ref{lmsts} pour $p<n$ et $q>0$, puisqu'alors $H_q({\rm Hom}_A(-,V))$ est un foncteur polynomial (car il est en fait exponentiel gradué par la formule de Künneth) nul en $0$ (on applique ici l'hypothèse de récurrence à $A^{op}$). On en déduit que $L_i$ est un foncteur constant pour $i\leq n$, en utilisant le lemme~\ref{rqfonct}. La suite spectrale du §\,\ref{ssf} et le théorème~\ref{thscogl} (dans le cas particulier où le foncteur contravariant est constant) donnent alors $H_n(GL_\infty(A) ; F_\infty)=0$ pour $F$ polynomial sans terme constant, d'où le théorème.
\end{proof}

\begin{rem}
Cette démonstration diffère profondément de celle de Betley. En effet, dans \cite{Bet2}, il établit le résultat à partir de  faits généraux sur la structure des foncteurs polynomiaux et surtout de théorèmes d'annulation démontrés dans son article antérieur \cite{Bet1}. Celui-ci se fondait sur des considérations explicites sur les groupes linéaires utilisant des arguments arithmétiques (différents selon qu'il s'agit de prouver l'annulation de la composante sans torsion de l'homologie ou de sa composante $p$-primaire pour un nombre premier $p$), qui permettent d'obtenir les annulations souhaitées, d'abord pour le cas crucial de $GL_\infty(\mathbb{Z})$. À l'inverse, la méthode ici suivie ne nécessite aucune considération d'ordre arithmétique, mais repose lourdement sur le théorème~\ref{thscogl}, résultat d'annulation tout-à-fait non trivial en homologie des foncteurs (même dans le cas particulier où le foncteur contravariant dans le groupe de torsion est constant).
\end{rem}

\bibliographystyle{smfalpha}
\bibliography{bibli-dv}
\end{document}

%% file: intro-dv.tex
\begin{abstract}
On calcule dans cet article l'homologie stable des groupes orthogonaux et symplectiques sur un corps fini $k$ à coefficients tordus par un endofoncteur usuel $F$ des $k$-espaces vectoriels (puissance extérieure, symétrique, divisée...). Par homologie stable, on entend, pour tout entier naturel $i$, les colimites des espaces vectoriels $H_i(O_{n,n}(k) ; F(k^{2n}))$ et $H_i(Sp_{2n}(k) ; F(k^{2n}))$ --- dans cette situation, la stabilisation (avec une borne explicite en fonction de $i$ et $F$) est un résultat classique de Charney. 

Tout d'abord, nous donnons un cadre formel pour relier l'homologie stable de certaines suites de groupes à l'homologie de petites catégories convenables, à l'aide d'une suite spectrale, qui dégénère dans de nombreux cas favorables. Cela nous permet d'ailleurs de retrouver des résultats de Betley sur l'homologie stable des groupes linéaires et des groupes symétriques, par des méthodes purement algébriques (sans recours à la $K$-théorie stable).

Pour une application exploitable de ce formalisme aux groupes orthogonaux ou symplectiques sur un corps fini, nous réinterprétons la deuxième page de notre suite spectrale en termes de foncteurs de Mackey non additifs et utilisons leurs propriétés d'acyclicité. Cela permet d'obtenir une simplification spectaculaire de la deuxième page de la suite spectrale en employant de puissants résultats d'annulation connus en homologie des foncteurs.

Dans le cas où les groupes orthogonaux ou symplectiques sont pris sur un corps fini et les coefficients à valeurs dans les espaces vectoriels sur ce même corps, nous pouvons mener le calcul de cette deuxième page grâce à des résultats classiques : annulation homologique à coefficients triviaux (Quillen, Fiedorowicz-Priddy), et calcul des groupes de torsion entre foncteurs usuels (Franjou-Friedlander-Scorichenko-Suslin, Cha{\l}upnik). 
Ceci permet de nombreux calculs d'homologie stable à coefficients.

\smallskip

\begin{center}
\textbf{Abstract}
\end{center}

We compute the stable homology of orthogonal and symplectic groups, over a finite field $k$, when the coefficients module is twisted by a usual endofunctor $F$ of $k$-vector spaces (e.g. an exterior, a symmetric, or a divided power) --- that is, for each natural integer $i$, we compute the colimit of the vector spaces $H_i(O_{n,n}(k) ; F(k^{2n}))$ and $H_i(Sp_{2n}(k) ; F(k^{2n}))$. Stabilization in this situation is a classical result of Charney.

We first set a formal framework, within which the stable homology of some families of groups, relates through a spectral sequence, to the homology of suitable small categories. The spectral sequence collapses in many cases. We illustrate this purely algebraic method to retrieve results of Betley for the stable homology of the general linear groups and of the symmetric groups.

We then apply our approach to orthogonal and symplectic groups over a finite field. To this end, we reinterpret the second page of our spectral sequence with Mackey functors and use their acyclicity properties. It allows us to simplify the second page of the spectral sequence, by using powerful cancellation results for functor homology.

For the orthogonal as for the symplectic groups over a finite field, and for coefficients modules over the same field, we compute the second page of the spectral sequence. Classical results prove useful at this point: homological cancellation with trivial coefficients (Quillen, Fiedorowicz-Priddy), and calculation of the torsion groups between usual functors (Franjou-Friedlander-Scorichenko-Suslin, Cha{\l}upnik). This provides extensive computations of stable homology with coefficients.
\end{abstract}

\smallskip

\noindent
\begin{small}{\em Mots-clés : } homologie stable, groupes orthogonaux, groupes symplectiques, homologie des foncteurs, foncteurs de Mackey non additifs.

\noindent
{\em Keywords : } stable homology, orthogonal groups, symplectic groups, functor homology, non-additive Mackey functors.

\noindent
{\em Classification math. : } 20J06, 20J05, 20G10, 18G40.
\end{small}

\section*{Introduction}
Cet article a pour objet l'homologie stable \`a coefficients tordus de familles de groupes classiques, c'est-à-dire la colimite de leur homologie ; celle-ci, dans de nombreux cas, est atteinte en temps fini pour chaque degré. Alors que cette homologie stable poss\`ede un comportement plus r\'egulier que l'homologie instable, elle s'av\`ere gén\'eralement inaccessible au calcul direct. Depuis les travaux de Betley  (\cite{Bet}) et Suslin  (\cite{FFSS}), on dispose cependant d'une interpr\'etation de l'homologie stable des groupes lin\'eaires en terme d'homologie des foncteurs, qui permet de mener \`a bien de nombreux calculs. N\'eanmoins aucun analogue n'\'etait jusqu'alors connu dans les cas des groupes orthogonaux et symplectiques. 
\par
Le présent travail établit un isomorphisme naturel entre l'homologie stable des groupes orthogonaux (ou symplectiques) sur un corps fini, à coefficients tordus par un foncteur polynomial, et des groupes de torsion entre endofoncteurs des espaces vectoriels. Il est remarquable que cet isomorphisme ne fasse pas intervenir de cat\'egorie d'espaces quadratiques (ou symplectiques), et qu'il ne s'exprime que par la catégorie déjà bien étudiée des endofoncteurs entre espaces vectoriels. Cela rend calculables les groupes d'homologie stable des groupes orthogonaux  pour les foncteurs polynomiaux usuels : puissances symétriques, extérieures, tensorielles, etc.

L'homologie stable à coefficients constants des groupes orthogonaux sur un corps fini $k$ a \'et\'e calcul\'ee dans les années $1970$ par Fiedorowicz et Priddy \cite{FPH}, en généralisant les méthodes initiées par Quillen \cite{Qui} pour les groupes lin\'eaires. L'homologie stable à coefficients dans un corps de même caractéristique que $k$ est triviale pour les groupes orthogonaux (si la caractéristique de $k$ est impaire), symplectiques ou linéaires sur $k$. De plus, on dispose de résultats de stabilité homologique pour les familles de groupes classiques sur les corps finis,  à coefficients constants ou tordus par un foncteur polynomial --- le cas des groupes orthogonaux étant dû à Charney \cite{Charney}. Pour autant, la d\'etermination de cette valeur stable semblait jusqu'\`a pr\'esent inabordable, y compris pour des coefficients tordus par un foncteur polynomial non constant \'el\'ementaire.

\smallskip 

L'annulation de l'homologie stable du groupe linéaire sur un anneau, à coefficients tordus par un foncteur polynomial sans terme constant, a \'et\'e obtenue par Betley \cite{Bet2} par des m\'ethodes compl\`etement diff\'erentes de celles utilisées pour les coefficients constants. En $1999$, Betley \cite{Bet} et Suslin \cite[appendice]{FFSS} ont démontré indépendamment une généralisation du résultat précédent, pour des coefficients tordus par un {\em bi}foncteur, polynomial en chaque variable. L'homologie stable n'est alors plus généralement nulle, mais naturellement isomorphe à l'homologie de Hochschild de la catégorie des $k$-espaces vectoriels de dimension finie, à coefficients dans le bifoncteur. Ces groupes d'homologie sont accessibles pour les foncteurs usuels, comme l'a notamment montré l'article \cite{FFSS}. La démonstration de Betley repose sur un analogue, en termes de groupes {\em algébriques}, du lien entre l'homologie des groupes linéaires et l'homologie de Hochschild des bifoncteurs polynomiaux, résultat établi un peu plus tôt par Friedlander et Suslin \cite{FS}. Suslin s'appuie entièrement, pour sa part, sur des considérations internes aux foncteurs entre espaces vectoriels. Sa démarche a été étendue peu après à l'homologie stable des groupes linéaires sur un anneau arbitraire par Scorichenko \cite{Sco}, qui a obtenu un isomorphisme entre $K$-théorie stable et homologie de Hochschild d'un bifoncteur polynomial.

Ces r\'esultats constituent une illustration de la richesse du point de vue des cat\'egories de foncteurs, dont on trouvera une synth\`ese dans \cite{FFPS}. Cette approche a d\'ej\`a permis de r\'esoudre de difficiles conjectures de finitude, comme en t\'emoignent les travaux de Friedlander-Suslin \cite{FS}, et tout r\'ecemment de Touz\'e et van der Kallen \cite{TvdK}. L'\'etude fine des cat\'egories de foncteurs s'est \'egalement poursuivie avec, par exemple, les r\'esultats cohomologiques de Cha\l upnik \cite{Chal}  et Touzé \cite{Tcu},
compl\'etant ceux de \cite{FFSS}, ou l'introduction de nouvelles cat\'egories de foncteurs reli\'ees aux groupes orthogonaux dans \cite{CV}.

\medskip

Nous présentons maintenant le contenu de l'article.
Une première section dégage un cadre formel pour étudier l'homologie stable d'une suite convenable de groupes, à partir de l'homologie d'une catégorie adaptée à la situation. 
Plus pr\'ecis\'ement, on consid\`ere une petite cat\'egorie mono\"idale sym\'etrique $(\C, \oplus, 0)$ dont l'unit\'e $0$ est objet initial, et $A$ un objet de $\C$. Pour chaque entier naturel $i$, on note $G(i)$ le groupe d'automorphismes ${\rm Aut}_\C(A^{\oplus i})$.
La cat\'egorie $\C$ sera par exemple celle des modules  projectifs de type fini sur un anneau $A$, avec les {\em monomorphismes scindés} comme flèches, et la somme directe pour structure mono\"idale. 
On peut prendre aussi pour $\C$ la cat\'egorie des ensembles finis avec injections, pour $A$ un ensemble \`a $1$ \'el\'ement et la r\'eunion disjointe pour $\oplus$, ce qui nous permettra de retrouver des résultats de Betley (\cite{Bet-sym}) sur l'homologie stable des groupes symétriques. 
Pour le sujet principal de cet article, on considère une catégorie $\C$ d'espaces quadratiques, ou symplectiques, de dimension finie sur un corps commutatif, pour $A$ un espace hyperbolique ou symplectique de dimension~$2$ et pour $\oplus$ la somme orthogonale. Ces exemples sont détaillés au paragraphe~\ref{p1-deux}.

Dans le cas général, la suite de morphismes :
$$0\to A \to\cdots \to A^{\oplus n} \to A^{\oplus(n+1)}\to\cdots$$
donn\'es par :
$$A^{\oplus n} \simeq A^{\oplus n} \oplus 0 \xrightarrow{Id \oplus (0 \to A)} A^{\oplus n} \oplus A \simeq  A^{\oplus(n+1)}$$
est compatible aux actions des groupes d'automorphismes $G(n)={\rm Aut}(A^{\oplus n})$, o\`u $G(n)$  agit sur $A^{\oplus(n+1)}$ via le morphisme : 
$$G(n) \xrightarrow{g \mapsto g \oplus A} G(n+1).$$ 
Ceci induit une suite naturelle de morphismes :
$$\cdots \to H_*(G(n); F(A^{\oplus n})) \to H_*(G(n+1); F(A^{\oplus (n+1)}))\to\cdots ,$$
o\`u $F$ est un foncteur de $\C$ vers les groupes ab\'eliens. On appelle {\it homologie stable} des groupes $G(n)$ \`a coefficients dans $F$ la colimite de cette suite.

Il existe toujours un morphisme naturel de l'homologie stable des groupes $G(n)$ \`a coefficients dans $F$ vers l'homologie de $\C$ à coefficients dans $F$. Sous de bonnes hypoth\`eses sur la catégorie $\C$, ce morphisme est  un isomorphisme en degré $0$, mais ce n'est plus en général le cas au-delà, où l'on obtient, au paragraphe~\ref{ssf}, une suite spectrale convergeant vers cette homologie stable et dont la deuxième page s'exprime par des groupes de torsion sur la catégorie~$\C$.

Cette suite spectrale offre une alternative purement algébrique à la suite spectrale de la $K$-théorie stable convergeant vers l'homologie stable du groupe linéaire à coefficients tordus, et à ses avatars la généralisant aux familles de groupes usuelles. En effet, dans les cas connus où l'on sait exprimer la $K$-théorie stable en termes algébriques plus simples, et où la suite spectrale correspondante s'arrête dès la deuxième page, on observe un phénomène analogue dans notre formalisme --- voir notamment nos propositions~\ref{fdo} et~\ref{arret2}. Ceci n'est guère surprenant dans la mesure où les travaux de Scorichenko (cas classique de la $K$-théorie stable) ou Betley (cas des groupes symétriques) n'utilisent en fait nullement la définition de la $K$-théorie stable, mais seulement l'existence d'une suite spectrale naturelle d'un certain type, et où l'article \cite{BP} de Betley et Pirashvili identifie dans de nombreux cas la $K$-théorie stable à un foncteur dérivé défini de façon purement algébrique.

\medskip
Notre r\'esultat principal relatif \`a l'homologie stable des groupes orthogonaux, d\'emontr\'e au paragraphe~\ref{thm-fond}, est le suivant :
\begin{thm-intro} \label{thm-1}
Soit $k$ un corps fini de caractéristique impaire. Pour tout foncteur \textbf{polynomial} (voir d\'efinition~\ref{poly}) $F$ entre $k$-espaces vectoriels,
il existe un isomorphisme naturel :
$$\underset{n \in \mathbb{N}}{\col}H_*(O_{n,n}(k);F(k^{2n})) \simeq {\rm Tor}^{\E^f_k}_*(V \mapsto k[S^2(V^*)],F)$$
o\`u $\E^f_k$ désigne la cat\'egorie des $k$-espaces vectoriels de dimension finie, et $S^2$ la seconde puissance symétrique.
\end{thm-intro}
Rappelons que les foncteurs $n$-i\`eme puissance tensorielle $T^n$ et $n$-i\`eme puissance sym\'etrique $S^n$ sont polynomiaux de degr\'e $n$.
Par les r\'esultats de \cite{Charney}, la colimite du th\'eor\`eme est donc atteinte pour $n$ fini, en chaque degr\'e homologique. Hormis pour cette considération de stabilisation, le choix des groupes orthogonaux $O_{n,n}$ plutôt que d'autres n'a pas d'importance : toute colimite analogue construite à partir de l'homologie d'autres groupes orthogonaux (associés à des formes quadratiques non dégénérées) sur $k$ est canoniquement isomorphe à celle qu'on considère (cf. remarque~\ref{rqsco}.\,\ref{qand}).

Le théorème~\ref{thm-1} s'obtient par la suite d'isomorphismes suivante :
\begin{enumerate}
\item \label{iso1}
 $$\underset{n \in \mathbb{N}}{\col}H_*(O_{n,n}(k);F(k^{2n})) \simeq H_*(\Eq; F)$$ 
 o\`u $\Eq$ est la cat\'egorie des $k$-espaces quadratiques non d\'eg\'en\'er\'es de dimension finie. Cet isomorphisme vient du cadre g\'en\'eral que nous avons \'evoqu\'e pr\'ec\'edemment et de la trivialité de $\underset{n \in \mathbb{N}}{\col}H_*(O_{n,n}(k);k)$ (pour $k$ de caractéristique impaire). 
 
 \item \label{iso2}
$$H_*(\Eq; F) \simeq H_*(\Eqd; F)$$ 
o\`u $\Eqd$ est la cat\'egorie des $k$-espaces quadratiques (\'eventuellement d\'eg\'en\'er\'es) de dimension finie avec pour morphismes les injections quadratiques. Cette \'etape constitue le c\oe ur de la d\'emonstration de ce th\'eor\`eme et en est la partie la plus d\'elicate.

\item   \label{iso3}
$$H_*(\Eqd; F) \simeq {\rm Tor}^{\E^f_{inj}}_*(V \mapsto k[S^2(V^*)],F)$$ o\`u $\E^f_{inj}$ est la sous-cat\'egorie des injections de $\E^f_k$. Cet isomorphisme s'obtient par adjonction en observant qu'une forme quadratique sur $V$ est un \'el\'ement de $S^2(V^*)$. 

\item \label{iso4}
$${\rm Tor}^{\E^f_{inj}}_*(V \mapsto k[S^2(V^*)],F) \simeq {\rm Tor}^{\E^f_k}_*(V \mapsto k[S^2(V^*)],F).$$
Cet isomorphisme est un cas particulier d'un r\'esultat de Suslin, essentiel pour exprimer l'homologie stable des groupes lin\'eaires en terme d'homologie des foncteurs. 
\end{enumerate}
Notons que les isomorphismes~\ref{iso1} et~\ref{iso3} valent pour un foncteur $F$ arbitraire alors que ceux des points~\ref{iso2} et~\ref{iso4} utilisent de mani\`ere fondamentale le caract\`ere polynomial de $F$. L'hypoth\`ese de finitude de $k$ n'intervient qu'\`a l'\'etape~\ref{iso2}.

Revenons plus en d\'etail sur celle-ci. Pour des raisons formelles, il existe une suite spectrale de Grothendieck convergente de la forme :
$$E^2_{p,q}={\rm Tor}_p^{\Eqd}( L_q, F) \Rightarrow  H_{p+q}(\Eq; F)$$
o\`u les foncteurs $L_q$ s'obtiennent \`a partir des foncteurs d\'eriv\'es d'un adjoint \`a la pr\'ecomposition par l'inclusion $\Eq \to \Eqd$. Les valeurs des  foncteurs $L_q$ sont donn\'ees par le $q$-i\`eme groupe d'homologie de sous-groupes explicites du groupe orthogonal. Elles sont l'aboutissement de suites spectrales de Serre dont on peut calculer la deuxi\`eme page, mais pas les diff\'erentielles. Ces valeurs sont donc inaccessibles pour $q > 0$.

Nous contournons ce problème en observant que ces foncteurs $L_q$ transforment l'inclusion d'un espace quadratique  dans sa somme orthogonale avec un espace {\em non dégénéré}, en un isomorphisme. On peut donc les définir depuis la catégorie de fractions où l'on inverse ces inclusions.
 On montre  au théorème~\ref{thm-fractions} l'équivalence  de cette cat\'egorie avec la {\em catégorie de Burnside} sur les espaces vectoriels avec injections. Autrement dit, nos foncteurs $L_q$ peuvent être vus comme des {\em foncteurs de Mackey non additifs} sur les espaces vectoriels avec injections, ou encore comme des familles de représentations des différents groupes linéaires, par un résultat général d'équivalence de Morita (cf. \cite{CV}). Gr\^ace au théorème d'annulation homologique principal de \cite{Dja}, on en d\'eduit un isomorphisme :
 $$E^2_{p,q}={\rm Tor}_p^{\Eqd}( L_q, F) \simeq {\rm Tor}_p^{\Eqd}( L_q(0), F).$$
 Or $L_q(0) \simeq H_q(O_{\infty}(k) ; k)$ est nul pour $q>0$ par le r\'esultat de Fiedorowicz-Priddy cit\'e plus haut. L'isomorphisme 2 en découle.
 
\medskip

La cons\'equence suivante du th\'eor\`eme~\ref{thm-1} illustre en quoi la (co)homologie stable est plus r\'eguli\`ere que la (co)homologie instable.
\begin{thm-intro} \label{thm-2}
Soient $k$ un corps fini de caractéristique impaire, $F$ et $G$ deux foncteurs polynomiaux entre $k$-espaces vectoriels et $i$, $j$ des entiers. Pour tout entier $n$ assez grand, le produit externe
$$H^i(O_{n,n}(k);F(k^{2n}))\otimes H^j(O_{n,n}(k);G(k^{2n}))\to H^{i+j}(O_{n,n}(k);F(k^{2n})\otimes G(k^{2n}))$$
est injectif.
\end{thm-intro}
Ce r\'esultat est  \'enonc\'e en terme de cohomologie pour traiter de produits plut\^ot que de coproduits.
Il s'obtient \`a partir de notre th\'eor\`eme principal par un raisonnement formel d\^u \`a Touz\'e (cf.  \cite{Touze}).

Le th\'eor\`eme~\ref{thm-1} permet aussi des calculs explicites, donn\'es aux théorèmes~\ref{thfc1} et~\ref{thfc2}. Nous obtenons entre autres, à l'aide du calcul de \cite{FFSS} des groupes d'extensions entre puissances divisées  sur un corps fini, que la cohomologie stable des groupes orthogonaux ou symplectiques à coefficients dans une algèbre polynomiale est elle-même une algèbre de polynômes. Plus précisément :
\begin{thm-intro} \label{thm-3}
Soit $k$ un corps fini de cardinal $q$ impair. L'algèbre de cohomologie stable des groupes orthogonaux (resp. symplectiques) sur $k$ à coefficients dans les puissances symétriques est polynomiale sur des générateurs $\alpha_{m,s}$ (resp. $\beta_{m,s}$) de bidegré $(2q^s m,q^s+1)$ indexés par des entiers $m\geq 0$ et $s\geq 0$ (resp. $s>0$), où le premier degré est le degré homologique.
\end{thm-intro}

\medskip

La partie formelle du présent travail (tout comme l'article \cite{Bet-sym} sur les groupes symétriques) ne traite que d'homologie stable à coefficients dans un foncteur, pas à coefficients dans un {\em bi}foncteur, comme le fait Scorichenko pour les groupes linéaires (à la suite de Betley et Suslin). Néanmoins, contrairement à ce qui advient lorsqu'on étudie l'homologie stable des groupes linéaires (où l'annulation à valeurs dans un foncteur polynomial sans terme constant contraste avec le résultat général pour un bifoncteur polynomial), l'homologie stable des groupes orthogonaux ou symplectiques à coefficients tordus par un bifoncteur ne s'avère pas plus générale que le cas particulier des foncteurs. De fait, toute forme quadratique ou symplectique non dégénérée sur un espace vectoriel $V$ déterminant un isomorphisme entre $V$ et son dual, l'homologie à coefficients tordus par un bifoncteur $B$ (avec une première variable contravariante et la seconde covariante) s'identifie à l'homologie à coefficients tordus par le {\em foncteur} $V\mapsto B(V^*,V)$ (cf. remarque~\ref{rq-bif-AT}). Comme application, on obtient au corollaire \ref{cr-radj} l'annulation stable de l'homologie du groupe orthogonal (ou symplectique) sur un corps fini de caractéristique impaire à coefficients dans sa représentation adjointe.
 
\smallskip

Enfin, signalons que Touzé a tout récemment obtenu (cf. \cite{Touze}), par des méthodes différentes, des résultats analogues à ceux de cet article pour la cohomologie {\em rationnelle} (i.e. comme groupes algébriques) stable des groupes orthogonaux et symplectiques.

\paragraph*{Organisation de l'article} La section~\ref{s-un} présente en détail le cadre g\'en\'eral adapt\'e  à notre approche de l'homologie stable. On y discute également une classe d'exemples fondamentale, qui contient tous nos cas d'application, et une condition supplémentaire qui intervient dans l'étude de la suite spectrale de la section~\ref{secss}. Celle-ci étudie, dans le formalisme de la section~\ref{s-un}, le morphisme puis la suite spectrale naturels qui relient l'homologie d'une catégorie convenable à l'homologie stable de la suite de groupes correspondante. On y discute notamment de la simplification de la deuxième page de cette suite spectrale, de son arrêt et de sa comparaison avec les suites spectrales classiques dérivées de constructions du type de la $K$-théorie stable.

La section~\ref{sid} constitue le c\oe ur de ce travail : elle donne les arguments non formels qui rendront accessible au calcul l'homologie stable des groupes orthogonaux ou symplectiques sur les corps finis, à coefficients tordus raisonnables. Elle identifie d'abord la catégorie de fractions des espaces éventuellement dégénérés où l'inclusion dans la somme orthogonale avec un espace non dégénéré est inversée. Ensuite, elle combine ce résultat avec les théorèmes d'annulation idoines connus en homologie des foncteurs pour obtenir le th\'eor\`eme~\ref{thm-1}.

La section~\ref{scal} donne les applications de ce th\'eor\`eme. On y traite de compatibilit\'e aux (co)produits pour obtenir notamment le th\'eor\`eme~\ref{thm-2}.  On effectue ensuite des calculs d'homologie stable de groupes orthogonaux ou symplectiques. Il faut distinguer la caractéristique impaire, qui se prête à des calculs complets (th\'eor\`eme~\ref{thm-3}), de la caractéristique~$2$ où les mêmes méthodes ne suffisent plus, mais où quelques résultats partiels sont déduits des travaux de Troesch (\cite{T1}).

Les trois premiers appendices donnent des rappels et des notations sur les catégories de foncteurs utilisés dans le corps de l'article : généralités d'algèbre homologique, foncteurs exponentiels puis quelques équivalences de catégories de foncteurs.

L'appendice~\ref{sect-cof} expose les deux résultats sur l'homologie des foncteurs (dus à Djament et Suslin) utilisés dans la démonstration du théorème principal de l'article et en rappelle les arguments.

Dans les deux derniers appendices, on montre comment retrouver rapidement à l'aide de notre formalisme des théorèmes dus à Betley sur l'homologie stable des groupes symétriques et linéaires.

\paragraph*{Quelques notations utilisées dans tout l'article}
On se donne un anneau commutatif (unitaire) $\kk$ de "base", au-dessus duquel tous les produits tensoriels non spécifiés seront pris. On désigne par $\mathbf{Mod}_\kk$ la catégorie des $\kk$-modules.

 Si $\C$ est une catégorie (essentiellement) petite, on note $\C-\mathbf{Mod}$ la catégorie des foncteurs de $\C$ vers $\mathbf{Mod}_\kk$. Quelques généralités sur cette catégorie abélienne sont rappelées dans l'appendice~\ref{appA}. On pose également $\mathbf{Mod}-\C=\C^{op}-\mathbf{Mod}$.
 
 Si $k$ est un corps commutatif, on note $\E(k)$ la catégorie des espaces vectoriels sur $k$ et $\E^f(k)$ (ou simplement $\E^f$) la sous-catégorie pleine des espaces de dimension finie.
 
 La précomposition par un foncteur $Q$ est notée $Q^*$.

On note $\mathbb{N}$ l'ensemble des entiers positifs ou nuls.

\paragraph*{Remerciements} Les auteurs témoignent leur gratitude à Vincent Franjou pour ses nombreuses discussions utiles à la réalisation de cet article. Ils remercient Serge Bouc pour leur avoir indiqué l'utilité de l'inversion de Möbius pour obtenir des équivalences de Morita, en particulier le théorème de Pirashvili à la Dold-Kan. Ils sont reconnaissants envers Antoine Touzé pour de fructueuses conversations qui ont permis d'améliorer la présentation ou les résultats de cet article en plusieurs occurrences. Ils le remercient notamment pour leur avoir signalé la remarque \ref{rq-bif-AT} permettant de traiter le cas des bifoncteurs, ainsi qu'une erreur présente dans une version préliminaire de ce travail, relative à la caractéristique~$2$.

Le second auteur remercie pour son hospitalité le laboratoire de mathématiques Jean Leray (Nantes), où une grande partie de ce travail a été élaboré, ainsi que le programme MATPYL de la Fédération de Recherche 2962 ``Math\'ematiques des Pays de Loire'' pour son soutien.

\tableofcontents

%% file: spec.tex
Dans cette section, on introduit une suite spectrale obtenue de mani\`ere purement alg\'ebrique,  convergeant vers l'homologie stable, \`a coefficients tordus, des familles de groupes qui nous int\'eressent. Cette suite spectrale est l'ingr\'edient original de la partie formelle de ce travail ; elle fournit une alternative \`a la suite spectrale de la $K$-th\'eorie stable consid\'er\'ee dans les travaux de Betley et Scorichenko.

\begin{conv}
{\bf Dans toute cette section, on se donne un triplet $(\C,S,G)$ vérifiant les hypothèses (C), (W) et (G) du paragraphe~\ref{p-un}.}
\end{conv}

\subsection{Morphisme de l'homologie stable vers l'homologie de catégorie} \label{section-2.1}

Le but de ce paragraphe est de construire des morphismes naturels :

\begin{equation} \label{eq1}
H_*(G_{\infty} ; F_{\infty}) \to H_*(\C ; F)
\end{equation}
o\`u $H_*(\C ; F)$ est l'homologie de la cat\'egorie $\C$ \`a coefficients dans un foncteur $F \in {\rm Ob}\,\C-\mathbf{Mod}$, dont on rappelle la d\'efinition dans l'appendice \ref{appA}, et
\begin{equation} \label{eq2}
H_*(G_{\infty} ; F_{\infty}) \to H_*(\C\times G_\infty ; \Pi^*F)
\end{equation}
où l'on voit le groupe $G_\infty$ comme une catégorie à un objet $\Pi : \C\times G_\infty\to\C$ désigne le foncteur de projection. (On rappelle que $\Pi^*: \C-\mathbf{Mod}\to\C\times G_\infty-\mathbf{Mod}$ désigne la pr\'ecomposition par $\Pi$.)

Ces morphismes naturels peuvent \^etre d\'efinis \`a partir des morphismes d'\'evaluation. On en donne ici une pr\'esentation utilisant une cat\'egorie auxiliaire not\'ee $\widetilde{\C}$, qui interviendra uniquement dans cette section. 

L'int\'er\^et de cette cat\'egorie provient de ce qu'on peut obtenir des r\'esolutions plates du foncteur constant $\kk\in {\rm Ob}\, \mathbf{Mod}-\widetilde{\C}$ \`a partir de r\'esolutions projectives du $G_{\infty}$-module trivial $\kk$. Ceci permet, entre autre, d'exprimer $H_*(G_{\infty} ; F_{\infty})$ en terme d'homologie de la cat\'egorie~$\widetilde{\C}$.

\begin{defi}
Soit $\widetilde{\C}$ la catégorie ayant $\mathbb{N}$ pour ensemble d'objets, telle que ${\rm Hom}_{\widetilde{\C}}(i,j)=G(j)$ pour $i\leq j$ et $\emptyset$ sinon. La composition est définie par $t'\circ t:=t'\cdot G(j\leq l)(t)$ pour $i\leq j\leq l$ et $t\in{\rm Hom}_{\widetilde{\C}}(i,j)=G(j)$ et $t'\in{\rm Hom}_{\widetilde{\C}}(j,l)=G(l)$. 
\end{defi}

\begin{nota}
On désigne par $\widetilde{S} : \mathbb{N}\to\widetilde{\C}$ le foncteur égal à l'identité sur les objets et associant à chaque relation $i\leq j$ de $\mathbb{N}$ le morphisme de $\widetilde{\C}$ correspondant à $1\in G(j)$. 

Pour la cohérence des notations, on note également $\widetilde{G}:=G : \mathbb{N}\to\mathbf{Grp}$.
\end{nota}

La propriété suivante est immédiate.

\begin{pr}\label{prtrax}
Le triplet $(\widetilde{\C},\widetilde{S},\widetilde{G})$ vérifie les hypothèses (C), (W') et (G).
\end{pr}

En revanche, même si $(\C,S,G)$ provient d'une structure monoïdale sur $\C$ comme au paragraphe~\ref{p1-deux}, il n'en est pas nécessairement de même pour $(\widetilde{\C},\widetilde{S},\widetilde{G})$ (la définition évidente que l'on est tenté de donner à partir de l'addition sur $\mathbb{N}$ et de la structure monoïdale sur $\C$ n'est pas toujours fonctorielle).

\begin{defi}
On note $Q: \widetilde{\C}\to \C$ le foncteur donn\'e par $Q(i)=S(i)$ sur les objets et $Q(f)=f\cdot S(i\leq j)$ pour $i\leq j$ dans $\mathbb{N}$ et $f \in {\rm Hom}_{\widetilde{\C}}(i,j)=G(j)$ sur les morphismes.
\end{defi}

On v\'erifie aussit\^ot la compatibilit\'e de $Q$ \`a la composition.
\smallskip

Pour $F \in {\rm Ob}\,\C-\mathbf{Mod}$, le foncteur $Q^*F$ fournit un morphisme naturel $H_*(\widetilde{\C}, Q^*F) \to H_*(\C, F)$.

Afin d'obtenir le morphisme naturel~(\ref{eq1}), on identifie dans la suite $H_*(\widetilde{\C} ; Q^*F)$ et $H_*(G_{\infty} ; F_{\infty})$. Pour cela nous aurons besoin du r\'esultat suivant qui explique
l'avantage de la cat\'egorie $\widetilde{\C}$ sur $\C$.

\begin{lm} \label{Pinf3}
Le foncteur $P_{\infty}^{\widetilde{\C}^{op}}(=\underset{i \in\mathbb{N}}{\col} P_i^{\widetilde{\C}^{op}}$, cf. lemme~\ref{Pinf}) est tel que pour tout $i\in\mathbb{N}$ on a un isomorphisme de $G_{\infty}$-modules :
$$P_{\infty}^{\widetilde{\C}^{op}}(i)\simeq\kk[G_{\infty}].$$
\end{lm}

\begin{proof} Pour $i,j\in\mathbb{N}$ tels que $j \leq i$ on a:  $P_i^{\widetilde{\C}^{op}}(j)=\kk[{\rm Hom}_{\widetilde{\C}}(j,i)]\simeq \kk[G(i)]$. Le r\'esultat s'en d\'eduit par passage \`a la colimite.
\end{proof}

\begin{rem}
Le foncteur $P_{\infty}^{\widetilde{\C}^{op}}$ n'est pas pour autant un foncteur constant : via l'isomorphisme du lemme précédent, son action sur les morphismes n'est pas donnée par l'identité. En effet, pour $l \leq j \leq i$ si $f \in {\rm Hom}_{\widetilde{\C}^{op}}(j,l)=G(j)$ on a $P_i^{\widetilde{\C}^{op}}(f)([g])=[g\cdot G(j\leq i)(f)]$, pour $g\in G(i)$. N\'eanmoins nous avons le r\'esultat suivant.
\end{rem}

\begin{lm}\label{lmdcn}
 Le foncteur $(P_{\infty}^{\widetilde{\C}^{op}})_{G_{\infty}}: \widetilde{\C}^{op} \to \mathbf{Mod}_\kk$ est constant en~$\kk$.
\end{lm}

\begin{proof}
Notons $\alpha_j : G(j)\to G_\infty$, pour $j\in\mathbb{N}$, le morphisme canonique. La remarque précédente montre que, dans l'isomorphisme du lemme~\ref{Pinf3}, le morphisme $\kk[G_\infty]\simeq P^{\widetilde{\C}^{op}}_\infty(j)\to  P^{\widetilde{\C}^{op}}_\infty(l)\simeq\kk[G_\infty]$ induit par un morphisme $l\to j$ de $\widetilde{\C}$ correspondant à un élément $f$ de $G(j)$ est la multiplication par $\alpha_j(f)$. Par passage aux coïnvariants il induit donc l'identité de~$\kk$.
\end{proof}

\begin{pr}\label{hotilde}
Il existe  un isomorphisme gradué
$$H_*(\widetilde{\C} ; Q^*F) \simeq H_*(G_{\infty} ; F_{\infty})$$
naturel en $F\in {\rm Ob}\,\C-\mathbf{Mod}$.
\end{pr}

\begin{proof}
Soit $R_\bullet\to\kk$ une r\'esolution projective de $\kk$ en tant que $G_{\infty}$-module. Par le lemme~\ref{Pinf3} le foncteur $-\underset{G_{\infty}}{\otimes} P_{\infty}^{\widetilde{\C}^{op}} : \mathbf{Mod}_{\kk[G_\infty]}\to\mathbf{Mod}-\widetilde{\C}$ est exact. On en d\'eduit un complexe exact $R_\bullet \otimes_{G_{\infty}} P_{\infty}^{\widetilde{\C}^{op}} \to \kk \otimes_{G_{\infty}} P_{\infty}^{\widetilde{\C}^{op}}$ dans $\mathbf{Mod}-\widetilde{\C}$. Or $\kk \otimes_{G_{\infty}} P_{\infty}^{\widetilde{\C}^{op}} \simeq (P_{\infty}^{\widetilde{\C}^{op}} )_{G_{\infty}} \simeq \kk$ d'apr\`es le lemme~\ref{lmdcn}. De plus les foncteurs $R_i \otimes_{G_{\infty}} P_{\infty}^{\widetilde{\C}^{op}}$ de $\mathbf{Mod}-\widetilde{\C}$ sont plats comme $P_{\infty}^{\widetilde{\C}^{op}}$ puisque les $G_\infty$-modules $R_i$ sont projectifs. On en d\'eduit que $R_\bullet \otimes_{G_{\infty}} P_{\infty}^{\widetilde{\C}^{op}} \to \kk \otimes_{G_{\infty}} P_{\infty}^{\widetilde{\C}^{op}} \simeq \kk$ est une r\'esolution plate de $\kk\in {\rm Ob}\,\mathbf{Mod}-\widetilde{\C}$.

On utilise maintenant l'isomorphisme canonique
$$(M\underset{G_\infty}{\otimes}P_{\infty}^{\widetilde{\C}^{op}})\underset{\widetilde{\C}}{\otimes} Q^*F\simeq M\underset{G_\infty}{\otimes}(P_{\infty}^{\widetilde{\C}^{op}}\underset{\widetilde{\C}}{\otimes} Q^*F)$$
naturel en $F$ et en le $G_\infty$-module $M$, associé aux isomorphismes naturels de $G_\infty$-modules
$$P_{\infty}^{\widetilde{\C}^{op}}\underset{\widetilde{\C}}{\otimes} Q^*F\simeq (Q^* F)_\infty\simeq F_\infty$$
dont le premier est déduit du lemme~\ref{Pinf} et de la proposition~\ref{prtrax} et le second s'obtient par inspection.
\end{proof}

On peut donc définir le morphisme canonique~(\ref{eq1}) par la composition
$$H_*(G_{\infty} ; F_{\infty})\xrightarrow{\simeq}H_*(\widetilde{\C} ; Q^*F)\to H_*(\C;F)$$
de l'isomorphisme de la proposition précédente et du morphisme induit par $Q$. On peut vérifier aisément que ce morphisme est toujours un isomorphisme en degré~$0$, ce qui découlera des résultats du paragraphe suivant, qui permettent d'étudier son comportement en tout degré.

\begin{rem}\label{rqdz}
Pour formel et élémentaire qu'il soit, ce résultat en degré $0$ peut déjà procurer un point de vue efficace sur des calculs de coïnvariants stables. Par exemple (anticipant  sur les résultats que nous donnerons par la suite en tout degré homologique, résultats qui sont eux non formels et considérablement plus difficiles qu'en degré $0$) si $\kk=k$ est un corps fini, on peut en déduire sans trop de peine le fait que
\begin{equation}\label{eqcoin}
\underset{n\in\mathbb{N}}{\col}H_0(O_{n,n}(k);\Gamma^*(k^{2n}))\quad\text{(resp. } \underset{n\in\mathbb{N}}{\col}H_0(Sp_{2n}(k);\Gamma^*(k^{2n})))
\end{equation}
est isomorphe à l'espace vectoriel des transformations naturelles de $\Gamma^*$ (foncteur (gradué) puissance divisée sur les $k$-espaces vectoriels) vers le foncteur $V\mapsto k^{S^2(V^*)}$ (resp. $V\mapsto k^{\Lambda^2(V^*)}$), où l'étoile indique cette fois la dualité.

L'article de Kuhn \cite{K-comp} (cf. son théorème~1.6) montre comment calculer ces espaces vectoriels gradués, à l'aide du lien fondamental établi entre les foncteurs entre $k$-espaces vectoriels et algèbre de Steenrod sur $k$ (au moins lorsque le corps fini $k$ est premier), établi par Henn, Lannes et Schwartz (on pourra consulter le premier article de Schwartz dans l'ouvrage \cite{FFPS} à ce sujet). Ce résultat est déjà remarquable car le calcul direct des coïnvariants~(\ref{eqcoin}) n'est pas du tout immédiat !
\end{rem}

Pour définir le morphisme~(\ref{eq2}), nous utiliserons le foncteur donné par la proposition immédiate suivante :
\begin{pr} Il existe un foncteur $J : \widetilde{\C}\to G_\infty$ envoyant chaque flèche $u\in {\rm Hom}_{\widetilde{\C}}(i,j)=G(j)$ sur son image canonique dans $G_\infty$.
\end{pr}

Le morphisme naturel~(\ref{eq2}) est défini par la composition
$$H_*(G_{\infty} ; F_{\infty})\xrightarrow{\simeq}H_*(\widetilde{\C} ; Q^*F)=H_*(\widetilde{\C} ; (Q,J)^*(\Pi^*F))\to H_*(\C\times G_\infty;\Pi^*F)$$
(on rappelle que $\Pi$ désigne la projection $\C\times G_\infty\to\C$) composé de l'isomorphisme de la proposition~\ref{hotilde} et du morphisme induit par $(Q,J)$.

On notera que le morphisme~(\ref{eq1}) n'est autre que la composée du morphisme~(\ref{eq2}) avec le morphisme naturel $H_*(\C\times G_\infty;\Pi^*F)\to H_*(\C;F)$ induit par $\Pi$.

\smallskip

La proposition~\ref{hotilde} admet la généralisation suivante :

\begin{pr}\label{pr-hst2} Il existe  un isomorphisme gradué
$$H_*(\widetilde{\C} ; (Q,J)^*X) \simeq H_*(G_{\infty} ; X_{\infty})$$
naturel en $X\in {\rm Ob}\,(\C\times G_\infty)-\mathbf{Mod}$, où $X_\infty=\col_{n\in\mathbb{N}} X(S(n))$ est muni de l'action diagonale de $G_\infty$ (cet espace vectoriel est naturellement muni d'une action de $G_\infty\times G_\infty$, dont un facteur agit comme dans la remarque~\ref{rqact} et l'autre par l'action tautologique de $G_\infty$ sur $\C\times G_\infty$).
\end{pr}

\begin{proof} Elle est complètement analogue à celle de la proposition~\ref{hotilde}, en remplaçant $Q^*F$ par $(Q,J)^*X$.
\end{proof}

\begin{rem} \label{rq-bif-AT}
 Pour donner une généralisation des considérations précédentes en termes de bifoncteurs, on a besoin de données supplémentaires. Par exemple, si $B$ est un bifoncteur sur la catégorie $\mathbb{P}(A)$ des $A$-modules à gauche projectifs de type fini sur un anneau $A$, i.e. un objet de $\mathbb{P}(A)^{op}\times\mathbb{P}(A)-\mathbf{Mod}$, on définit l'homologie stable des groupes linéaires sur $A$ à coefficients dans $B$ comme étant $H_*(GL_\infty(A);B_\infty)$, où $B_\infty$ est la colimite des $B(A^n,A^n)$ construite à partir des projections $A^{n+1}\twoheadrightarrow A^n$ sur les $n$ premiers facteurs (pour la première variable, contravariante) et des inclusions $A^n\hookrightarrow A^{n+1}$ des $n$ premiers facteurs (pour la seconde variable, covariante).
 
 Néanmoins, dans le cas des groupes orthogonaux ou symplectiques, on n'obtient rien de plus général par une telle procédure. En effet, toute forme quadratique non dégénérée sur un espace vectoriel déterminant un isomorphisme de celui-ci sur son dual, la catégorie $\Eq$ est équivalente à sa catégorie opposée, de sorte que l'homologie stable des groupes orthogonaux à coefficients dans un bifoncteur sur $\Eq$ (définie comme dans la situation précédente) n'est autre que l'homologie stable du foncteur obtenu en précomposant avec 
 $$\Eq\xrightarrow{{\rm diag}}\Eq\times\Eq\xrightarrow{\simeq}(\Eq)^{op}\times\Eq.$$
 On peut procéder de manière analogue pour les formes symplectiques.
\end{rem}

\subsection{Suite spectrale fondamentale} \label{ssf}

Soit $(Q,J)_! : \mathbf{Mod}-\widetilde{\C}\to\mathbf{Mod}-(\C\times G_\infty)$ le foncteur défini (à isomorphisme canonique près) par l'isomorphisme naturel
$$\forall X\in {\rm Ob}\,(\C\times G_\infty)-\mathbf{Mod}\;\;\forall Y\in {\rm Ob}\,\mathbf{Mod}-\widetilde{\C}\qquad Y\underset{\widetilde{\C}}{\otimes}(Q,J)^*(X)\simeq (Q,J)_!(Y)\underset{\C\times G_\infty}{\otimes}X$$
(cf. proposition~\ref{kan-t}).

Cet isomorphisme se dérive en une suite spectrale de Grothendieck (homologique) de terme $E^2$ donné par
$$E^2_{p,q}={\rm Tor}^{\C\times G_\infty}_p(\mathbb{L}_q (Q,J)_!(Y),X)$$
et d'aboutissement ${\rm Tor}^{\widetilde{\C}}_*(Y,(Q,J)^*(X))$. Les foncteurs $\mathbb{L}_q (Q,J)_!(Y)$ désignent les dérivés à gauche du foncteur exact à droite $(Q,J)_!$ ; cette suite spectrale est naturelle en $X$ et concentrée dans le premier quadrant, donc en particulier convergente.

On s'intéresse maintenant au cas où $Y$ est le foncteur constant $\kk$.

\begin{pr}\label{annulss1}
Pour tout entier $q>0$, on a $\mathbb{L}_q(Q,J)_!(\kk)=0$.

De plus, $(Q,J)_!(\kk)$ est donné par un isomorphisme naturel
$(Q,J)_!(\kk)(c)\simeq (P^\C_c)_\infty$ pour $c\in {\rm Ob}\,\C$.
\end{pr}

\begin{proof} La formule~(\ref{kander}) (appendice~\ref{appA}, proposition~\ref{kan-t}) et la proposition~\ref{pr-hst2} procurent des isomorphismes
$$\mathbb{L}_q(Q,J)_!(\kk)(c)=H_q(\widetilde{\C};(Q,J)^*(P^{\C\times G_\infty}_c))\simeq H_q(G_\infty ; (P^{\C\times G_\infty}_c)_\infty).$$

Comme $(P^{\C\times G_\infty}_c)_\infty\simeq (P^\C_c)_\infty\otimes\kk[G_\infty]$ comme $G_\infty$-modules, cela démontre la proposition.
\end{proof}

On déduit donc de notre suite spectrale et de la proposition~\ref{pr-hst2} :
\begin{cor}\label{hsexp1}
Il existe un isomorphisme gradué naturel
$$H_*(G_\infty ; X_\infty)\simeq {\rm Tor}_*^{\C\times G_\infty}((Q,J)_!(\kk),X)$$
pour $X\in {\rm Ob}\,(\C\times G_\infty)-\mathbf{Mod}$.

En particulier, il existe un isomorphisme naturel
$$H_*(G_\infty ; F_\infty)\simeq {\rm Tor}_*^{\C\times G_\infty}((Q,J)_!(\kk),\Pi^* F)$$
pour $F\in {\rm Ob}\,\C-\mathbf{Mod}$.
\end{cor}

Pour étudier plus avant ces groupes, nous aurons besoin du résultat classique suivant sur l'homologie d'un produit de deux catégories :
\begin{pr}\label{homprod}
Soient $\A$ et $\B$ deux petites catégories et $\Pi : \A\times\B\to\A$ le foncteur de projection.

Il existe une suite spectrale (du premier quadrant)
$$E^2_{p,q}={\rm Tor}^A_p(A\mapsto H_q(\B;Y(A,-)),F)\Rightarrow {\rm Tor}^{\A\times\B}_{p+q}(Y,\Pi^*F)$$
fonctorielle en $F\in {\rm Ob}\,\A-\mathbf{Mod}$ et $Y\in {\rm Ob}\,\mathbf{Mod}-(\A\times\B)$.

En particulier, il existe deux suites spectrales de Künneth
$${\rm I}^2_{p,q}=H_p(\A;A\mapsto H_q(\B;X(A,-)))\Rightarrow H_{p+q}(\A\times\B ; X),$$
$${\rm II}^2_{p,q}=H_p(\B;B\mapsto H_q(\A;X(-,B)))\Rightarrow H_{p+q}(\A\times\B ; X)$$
fonctorielles en $X\in {\rm Ob}\,(\A\times\B)-\mathbf{Mod}$.
\end{pr}

\begin{proof}
Il existe un isomorphisme naturel
$$Y\underset{\A\times\B}{\otimes}\Pi^*(F)\simeq (A\mapsto H_0(\B ; Y(A,-)))\underset{\A}{\otimes}F.$$
La suite spectrale recherchée est la suite spectrale de Grothendieck correspondante.
\end{proof}

Afin d'examiner la forme que prend cette suite spectrale dans le cas qui nous intéresse, nous introduisons les définitions suivantes.

\begin{nota}\label{not-st}
Soient $i$ et $j$ deux éléments de $\mathbb{N}$ tels que $i\leq j$. Nous désignerons par $St_\C(i,j)$ le stabilisateur de $S(i\leq j)\in {\rm Hom}_\C(S(i),S(j))$ sous l'action à gauche canonique de $G(j)$.

Nous noterons $St_\C(i)$ la colimite sur $j\geq i$ des groupes $St_\C(i,j)$. Nous identifierons $St_\C(i)$ avec son image dans le groupe $G_\infty$. 
\end{nota}

On remarque que l'on a $St_\C(j)\subset St_\C(i)$ pour $i\leq j$. Cela permet de considérer $St_\C$ comme un foncteur de but $\mathbf{Grp}$ et de source $\mathbb{N}^{op}$.

Nous aurons besoin, dans les cas où le foncteur $S$ n'est pas essentiellement surjectif, de la généralisation suivante de ces stabilisateurs, qui s'effectue au prix d'une perte de fonctorialité. 

\begin{nota}\label{nost2}
Choisissons, conformément à l'axiome (C), un morphisme $c\xrightarrow{u_c}S(i_c)$ de $\C$ pour tout objet $c$ de $\C$. Pour tout entier naturel $j\geq i_c$, on note $St(c,j)$ le stabilisateur de $c\xrightarrow{u_c}S(i_c)\xrightarrow{S(i_c\leq j)}S(j)$ sous l'action à gauche canonique de $G(j)$.

Nous noterons $St(c)$ la colimite sur $j\geq i_c$ des groupes $St(c,j)$. Nous identifierons $St(c)$ avec son image dans le groupe $G_\infty$.
\end{nota}

L'hypothèse (W) montre que changer le choix des $i_c$ et des $u_c$ ne modifie pas, à conjugaison près, les groupes $St(c)$ obtenus. Pour la même raison, tout morphisme $f : b\to c$ de $\C$ fait de $St(c)$ un sous-groupe d'un conjugué de $St(b)$ dans $G_\infty$.

Malgré le manque de fonctorialité sur les groupes $St(c)$, la remarque précédente et la trivialité en homologie de l'action des automorphismes intérieurs d'un groupe font de $c\mapsto H_*(St(c);\kk)$ un foncteur de $\C^{op}$ vers les $\kk$-modules gradués. 

\begin{pr}\label{fctder} Il existe une suite spectrale 
$$E^2_{p,q}={\rm Tor}_p^\C(c\mapsto H_q(St(c);\kk),F)\Rightarrow H_{p+q}(G_\infty ; F_\infty)$$
naturelle en $F\in {\rm Ob}\,\C-\mathbf{Mod}$.
\end{pr}

\begin{proof} Le $G_\infty$-module $(P^{\C}_c)_\infty$ s'identifie à $\kk[G_\infty/St(c)]$ grâce à l'axiome (W). La conclusion découle donc des propositions~\ref{annulss1} et~\ref{homprod} et du lemme de Shapiro.
\end{proof}

\begin{rem}\begin{enumerate}
\item En utilisant la proposition~\ref{hotilde} plutôt que la proposition~\ref{pr-hst2}, on obtient une suite spectrale naturelle
$$E^2_{p,q}={\rm Tor}_p^\C(\mathbb{L}_*Q_!(\kk),F)\Rightarrow H_{p+q}(G_\infty ; F_\infty).$$
On vérifie aisément qu'existent des isomorphismes naturels $\mathbb{L}_qQ_!(\kk)(c)\simeq H_q(St(c);\kk)$ et que la suite spectrale précédente est isomorphe à celle de la proposition~\ref{fctder}.

L'intérêt de la présentation qu'on donne réside dans la possibilité d'utiliser des arguments généraux d'effondrement pour les suites spectrales de la proposition~\ref{homprod} comme on le verra dans le paragraphe suivant.
\item Les suites spectrales de ce paragraphe possèdent également des propriétés de fonctorialité relativement à $\C$ qu'on laisse au lecteur le soin d'énoncer.
\end{enumerate}
\end{rem}

\begin{nota}\label{not-fdl} Lorsqu'aucune ambiguïté n'est possible sur $(\C,S,G)$ ou $\kk$, on notera $L_q$ le foncteur $c\mapsto H_q(St(c);\kk)$.
\end{nota}

C'est dans la proposition suivante que l'intérêt de l'hypothèse (S) apparaît.

\begin{pr}\label{fdo} Supposons que les foncteurs $S$ et $G$ proviennent d'une structure monoïdale symétrique sur $\C$ comme au paragraphe~\ref{p1-deux} et que l'hypothèse (S) est satisfaite. 

Alors pour tous entiers $n$, $i$ et tout objet $c$ de $\C$, le foncteur  $L_n$ transforme le morphisme canonique $c\to S(i)\oplus c$ en un isomorphisme.

\smallskip

Si de plus $S$ vérifie l'hypothèse (C'), alors le foncteur $L_n$ est constant en $H_n(G_\infty;\kk)$. La suite spectrale de la proposition~\ref{fctder} prend donc naturellement la forme
$$E^2_{p,q}\simeq {\rm Tor}^\C_p(H_q(G_\infty ; \kk),F)\Rightarrow H_{p+q}(G_\infty ; F_\infty).$$
\end{pr}

\begin{proof} On peut supposer que $i_{S(i)\oplus c}=S(i)\oplus i_c$ et $u_{S(i)\oplus c}=S(i)\oplus u_c$ (cf. supra). Pour tout entier $j\geq i_c$, le foncteur $S(i)\oplus -$ induit donc un isomorphisme $St(c,j)\xrightarrow{\simeq} St(S(i)\oplus c,i+j)$, par l'axiome (S).

Identifiant ces deux groupes via cet isomorphisme, on voit que la flèche induite par le morphisme canonique $c\to S(i)\oplus c$ s'identifie au morphisme de groupes $St(c,j)\hookrightarrow St(c,i+j)$ induit par 
$$G(j)={\rm Aut}_\C(A^{\oplus j})\to G(i+j)={\rm Aut}_\C(A^{\oplus (i+j)})\quad u\mapsto A^{\oplus i}\oplus u.$$

Celui-ci est conjugué au morphisme $u\mapsto u\oplus A^{\oplus i}$ qui intervient dans la colimite définissant $G_\infty$ par l'automorphisme de $A^{\oplus (i+j)}$ 
$$A^{\oplus (i+j)}=A^{\oplus i}\oplus A^{\oplus j}\simeq A^{\oplus j}\oplus A^{\oplus i}=A^{\oplus (i+j)}$$
donné par l'échange des facteurs.

Comme les automorphismes intérieurs n'agissent pas en homologie, cela fournit la première partie du résultat, grâce à la proposition~\ref{fctder}.

\smallskip

Supposons maintenant que $S$ vérifie l'hypothèse (C'). Soient $c$ un objet de $\C$ et $b$ un objet de $\C$ tel qu'existent $i\in\mathbb{N}$ et un isomorphisme $b\oplus c\simeq S(i)$. Ce qui précède montre que, pour tout $a\in {\rm Ob}\,\C$ et tout $n\in\mathbb{N}$, l'application linéaire $L_n(a\to a\oplus c)$ a un inverse à droite, donné à isomorphisme près par $L_n(a\oplus  c\to a\oplus c\oplus b)$ (qui a donc lui un inverse à gauche !) ; pour la même raison, ce dernier morphisme a un inverse à droite. Donc $L_n(a\oplus  c\to a\oplus c\oplus b)$ puis $L_n(a\to a\oplus c)$ sont des isomorphismes.

Considérons à présent un morphisme $f : d\to c$ quelconque de $\C$ et montrons que c'est un $L_*$-isomorphisme. On choisit  des objets $a$ et $b$ de $\C$ et des entiers $i$ et $j$ tels qu'existent des isomorphismes $a\oplus d\simeq S(j)$ et $b\oplus c\simeq S(i)$. On peut ensuite trouver, par l'axiome (W), un entier $l$, qu'on peut supposer supérieur à $i$ et $j$, et un automorphisme $g$ de $S(l)$ tel que le diagramme suivant commute :
$$\xymatrix{d\ar[d]_f\ar[r] & a\oplus d\simeq S(j)\ar[r]^-{S(j\leq l)} & S(l)\ar[d]^g \\
c\ar[r] & b\oplus c\simeq S(i)\ar[r]^-{S(i\leq l)} & S(l).
}$$
Toutes les flèches horizontales sont des  $L_*$-isomorphismes par ce qui précède, de même que $g$ qui est un isomorphisme, donc c'est aussi le cas de $f$, ce qui achève la démonstration.
\end{proof}

\begin{rem}
C'est dans cette démonstration qu'apparaît l'intérêt de supposer symétrique la structure monoïdale de $\C$, hypothèse qui n'intervient nulle part ailleurs dans les constructions ou démonstrations.
\end{rem}

\begin{rem}\label{rqsco}
\begin{enumerate}
\item\label{sco2} Dans le cas de la catégorie $\mathbb{M}(A)$ des modules à gauche projectifs de type fini sur un anneau $A$ avec injections scindées (cf. exemple~\ref{exfond}.\,\ref{ex-sco}), l'hypothèse (C') est satisfaite mais pas (S) ; néanmoins, la conclusion de la proposition précédente a encore lieu. Cela provient de résultats établis par Betley dans~\cite{Bet2}. Nous en donnons une démonstration, fondée sur les travaux de Scorichenko, dans l'appendice~\ref{appBS}.
\item \label{bets2} Toutes les hypothèses de la proposition sont satisfaites par la catégorie $\Theta$ de l'exemple~\ref{exfond}.\,\ref{exbs}. Cela permet de retrouver les résultats de Betley (cf. \cite{Bet-sym}) sur l'homologie stable des groupes symétriques. C'est ce que nous ferons dans l'appendice~\ref{apsym}.
\item\label{qand} Toutes les hypothèses de la proposition~\ref{fdo} sont également satisfaites dans les catégories $\Eq(k)$ et $\E_{alt}(k)$ (où $k$ est un corps commutatif) munies de la somme orthogonale et de l'objet $\mathbf{H}$ (cf. exemples~\ref{exfond}.\,\ref{exquad} et~\ref{exfond}.\,\ref{exsymp}). Cependant, les résultats de ce corollaire sont peu maniables dans ce cas, ce qui nous amènera à travailler plutôt dans les catégories $\Eqd$ et $\E_{alt}^{deg}$, dans les sections suivantes.

On remarquera déjà que, le choix de l'espace quadratique $\mathbf{H}$ n'ayant pas réellement d'importance (cf. remarque~\ref{rq-goinf}), l'homologie de la catégorie $\Eq$ calcule l'homologie stable de {\em tous} les groupes orthogonaux sur $k$ : le morphisme canonique de la colimite des homologies des groupes orthogonaux $O_{n,n}$ (qui correspondent au choix de $\mathbf{H}$) vers la colimite des homologies de tous les groupes orthogonaux (associés à des formes quadratiques non dégénérées) est un isomorphisme. (On peut le voir directement très rapidement : c'est une simple conséquence de l'inaction des automorphismes intérieurs en homologie, comme dans la démonstration de la proposition~\ref{fdo}.)
\end{enumerate}
\end{rem}

\subsection{Arrêt de la suite spectrale et comparaison à la $K$-théorie stable}

Nous commençons par donner une condition suffisante pour que le morphisme
$H_*(G_{\infty} ; F_{\infty}) \to H_*(\C\times G_\infty ; \Pi^*F)$ du paragraphe \ref{section-2.1} soit un isomorphisme.

\begin{pr}\label{isoc1}
Supposons que la catégorie possède un objet initial, noté $0$, de sorte qu'il existe, pour tout objet $G$ de $\mathbf{Mod}-\C$, un morphisme naturel de $G$ vers le foncteur $G(0)$.

Soit $F\in {\rm Ob}\,\C-\mathbf{Mod}$ tel que le morphisme $L_q\to L_q(0)$ induise un isomorphisme ${\rm Tor}_*^\C(L_q,F)\xrightarrow{\simeq} {\rm Tor}_*^\C(L_q(0),F)$ en homologie pour tout $q\in\mathbb{N}$.

Alors le morphisme~(\ref{eq2}) : $H_*(G_\infty ; F_\infty)\to H_*(\C\times G_\infty ; \Pi^*F)$ est un isomorphisme.
\end{pr}

\begin{proof}
Par la proposition~\ref{homprod} et le corollaire~\ref{hsexp1}, $H_*(G_\infty ; F_\infty)\simeq {\rm Tor}^{\C\times G_\infty}_*(c\mapsto (P^\C_c)_\infty,\Pi^*F)$ et $H_*(\C\times G_\infty ; \Pi^*F)\simeq {\rm Tor}^{\C\times G_\infty}_*(\kk,\Pi^*F)$ sont respectivement l'aboutissement de suites spectrales convergentes de deuxièmes pages :
$${\rm I}^2_{p,q}={\rm Tor}^\C_p(L_q,F)\qquad\text{et}\qquad {\rm II}^2_{p,q}={\rm Tor}^\C_p(L_q(0),F)$$
(puisque $L_q(0)\simeq H_q(G_\infty;\kk)$). Le morphisme $(c\mapsto (P^\C_c)_\infty)\to (P^\C_0)_\infty\simeq\kk$ induit par hypothèse un isomorphisme ${\rm I}^2_{p,q}\to {\rm II}^2_{p,q}$, il induit donc un isomorphisme entre les aboutissements des suites spectrales, d'où la proposition.
\end{proof}

\`A partir de ce résultat, nous donnons un critère simple pour que la suite spectrale qui nous intéresse s'effondre à la deuxième page ; on obtient même mieux : non seulement le terme $E^\infty$, mais aussi l'aboutissement de la suite spectrale (qui en diffère par une filtration), sont isomorphes au terme $E^2$.

\begin{pr}\label{arret2}
Faisons les trois hypothèses suivantes :
\begin{enumerate}
\item l'anneau $\kk$ est de dimension homologique au plus $1$ ;
\item la catégorie $\C$ possède un objet initial $0$ ;
\item $F\in {\rm Ob}\,\C-\mathbf{Mod}$ est tel que pour tout entier $q$, le morphisme $L_q\to L_q(0)$ induise un isomorphisme ${\rm Tor}_*^\C(L_q,F)\xrightarrow{\simeq} {\rm Tor}_*^\C(L_q(0),F)$ en homologie.
\end{enumerate}

Alors il existe des isomorphismes
$$H_n(G_\infty ; F_\infty)\simeq H_n(\C\times G_\infty;\Pi^* F)\simeq\bigoplus_{p+q=n} {\rm Tor}^\C_p(H_q(G_\infty;\kk),F)$$
naturels en $F$.
\end{pr}

Les cas les plus importants sont ceux où $\kk$ est un corps ou égale $\mathbb{Z}$. On remarque que sous les hypoth\`eses de la proposition~\ref{fdo} les conditions $2$ et $3$ de la proposition pr\'ec\'edente sont satisfaites.

\begin{proof} Comme $\C$ a un objet initial, tout foncteur constant de $\mathbf{Mod}-\C$ en un $\kk$-module injectif est injectif et représente les morphismes de l'évaluation en $0$ vers ledit module. Par conséquent,
$${\rm Ext}^r_{\C^{op}}(M,N) \simeq {\rm Ext}^r_\kk(M,N)$$
lorsque $M$ et $N$ sont des $\kk$-modules, vus dans le terme de gauche comme foncteurs constants depuis $\C^{op}$. Ces groupes sont donc nuls pour $r\geq 2$. 

Soit  $C_\bullet$ le complexe de $\mathbf{Mod}-\C$ obtenu en prenant les coïnvariants par rapport à $G_{\infty}$ d'une r\'esolution projective de $\kk\in {\rm Ob}\,\mathbf{Mod}-(\C\times G_{\infty})$. 
La deuxième suite spectrale de la proposition~\ref{homprod} s'obtient en prenant le produit tensoriel au-dessus de $\C$ de $C_\bullet$ et d'une r\'esolution projective de $F$. Les propositions~\ref{obstruction}  (appendice~\ref{appA}) et~\ref{isoc1} donnent alors la conclusion.  
\end{proof}
 
On peut également utiliser l'autre suite spectrale pour l'homologie de la catégorie $\C\times G_\infty$ donnée par la proposition~\ref{homprod}, par l'intermédiaire de la propriété générale suivante :
\begin{pr}\label{eff-kun}
Soient $\A$ et $\B$ deux petites catégories, $\Pi : \A\times\B\to\A$ la projection et $F$ un objet de $\A-\mathbf{Mod}$. On suppose que l'anneau $\kk$ est de dimension homologique au plus~$1$.

Alors la suite spectrale
$$E^2_{p,q}(F)=H_p(\B; H_q(\A;F))\Rightarrow H_{p+q}(\A\times \B ; \Pi^*F)$$
donnée par la proposition~\ref{homprod} (où $H_q(\A;F)$ est vu comme objet constant de $\B-\mathbf{Mod}$) s'effondre à la deuxième page.

De surcroît, cette suite spectrale induit un scindement ({\em non nécessairement naturel} en $F$)
$$H_n(\A\times \B ; \Pi^*F)\simeq\bigoplus_{p+q=n}H_p(\B; H_q(\A;F)).$$
\end{pr}

\begin{proof}
On suit la démarche de la section~$5$ de l'article "Stable $K$-theory is bifunctor homology" de Franjou et Pirashvili, dans le volume \cite{FFPS}, qui elle-même s'inspire du théorème~$2$ de l'article \cite{BP} de Betley et Pirashvili.

Comme dans \cite{FFPS}, on remarque que cette suite spectrale est une {\em $\partial$-suite spectrale} (cf. $2.2$ de l'article de Franjou-Pirashvili), c'est-à-dire que pour tout suite exacte courte $0\to N\to P\to F\to 0$ de $\A-\mathbf{Mod}$ et tous entiers $r\geq 2$, $p$ et $q$, on dispose d'un morphisme $\partial^r : E^r_{p,q}(F)\to E^r_{p,q-1}(N)$ de sorte que $\partial^{r+1}$ est le morphisme induit par $\partial^r$ en homologie et que $\partial$ commute à la différentielle de la suite spectrale. Les morphismes $\partial^2$ sont induits par le morphisme de liaison $H_q(\A ; F)\to H_{q-1}(\A ; N)$ ; la structure de $\partial$-suite spectrale s'obtient en l'étendant aux complexes de chaînes sur $\A-\mathbf{Mod}$ (i.e. en passant à la catégorie dérivée).

On choisit ensuite une suite exacte courte $0\to N\to P\to F\to 0$ avec $P$ projectif : le mor\-phisme de liaison  $H_q(\A ; F)\to H_{q-1}(\A ; N)$ est un isomorphisme pour $q\geq 1$ et un monomorphisme pour $q=1$, lequel est {\em scindé} parce que son conoyau est un sous-module du module projectif $H_0(\A ; P)$ (si $P$ est un projectif standard $P^\A_a$, $H_0(\A ; P)\simeq\kk$), donc est lui-même projectif par l'hypothèse faite sur $\kk$. Par conséquent, $\partial^2 : E^2_{p,q}(F)\to E^2_{p,q-1}(N)$ est injectif, ce qui nous permet d'appliquer le lemme~$2.2$ de l'article de Franjou-Pirashvili susmentionné pour conclure à l'effondrement de la suite spectrale à la deuxième page.

Le scindement se démontre de façon analogue, mais indépendante, par récurrence sur $n$ (on s'inspire ici directement de \cite{BP}) : il est immédiat pour $n=0$, on suppose donc $n>0$ et l'on se donne comme avant une suite exacte $0\to N\to P\to F\to 0$ avec $P$ projectif, qui induit des isomorphismes $H_i(\A ; F)\simeq H_{i-1}(\A ; N)$ pour $i>1$ et un scindement (non naturel) $H_0(\A ; N)\simeq K\oplus H_1(\A ; F)$, où $K=Ker\,(H_0(\A ; P)\to H_0(\A ; F))$ est un $\kk$-module projectif. On considère alors le diagramme
$$\xymatrix{H_n(\A\times\B ; \Pi^*P)\ar[d]^\simeq_h\ar[r]^a & H_n(\A\times\B ; \Pi^*F)\ar[d]^f\ar[r]^b & H_{n-1}(\A\times\B ; \Pi^*N)\ar[d]^\simeq_g \\
H_n(\B ; H_0(\A ; P))\ar[r]^-c & \underset{p+q=n}{\bigoplus} H_p(\B ; H_q(\A ; F))\ar[r]^-d & \underset{p+q=n-1}{\bigoplus}H_p(\B ; H_q(\A ; N))
}$$
dans lequel :
\begin{enumerate}
\item la ligne supérieure est une partie de la suite exacte longue du foncteur homologique $H_*(\A\times\B ; -)$ ;
\item $h$ est le "coin" de la suite spectrale $(E^r_{p,q}(P))$, qui est un isomorphisme parce que $E^2_{p,q}(P)=0$ pour $q>0$, $P$ étant projectif ;
\item $g$ est un isomorphisme donné par l'hypothèse de récurrence ;
\item $f$ est le morphisme dont la composante $H_n(\A\times\B ; \Pi^*F)\to H_n(\B ; H_0(\A ; F))$ est le coin de la suite spectrale et la composante $H_n(\A\times\B ; \Pi^*F)\to H_{n-i}(\B ; H_i(\A ; F))$ est pour $i>0$ composée de
$$H_n(\A\times\B ; \Pi^*F)\xrightarrow{b}H_{n-1}(\A\times\B ; \Pi^*N)\xrightarrow{g}\underset{p+q=n-1}{\bigoplus}H_p(\B ; H_q(\A ; N)\twoheadrightarrow H_{n-i}(\B ; H_{i-1}(\A ; N))$$
et de la flèche induite par l'isomorphisme $H_i(\A ; F)\simeq H_{i-1}(\A ; N)$ pour $i>1$ et le monomorphisme scindé $H_1(\A ; F)\hookrightarrow K\oplus H_1(\A ; F)\simeq H_0(\A ; N)$ pour $i=1$ ;
\item le morphisme $d$ est défini de manière évidente pour que le carré de droite commute ;
\item le morphisme $c$ est la composée du morphisme $H_n(\B ; H_0(\A ; P))\to H_n(\B ; H_0(\A ; F))$ induit par $P\twoheadrightarrow F$ et de l'inclusion canonique. Cela assure que le carré de gauche commute, par naturalité du coin de la suite spectrale et par nullité de la composée $ba$.
\end{enumerate}
L'exactitude de la suite $H_n(\B ; H_0(\A ; P))\to H_n(\B ; H_0(\A ; F))\to H_{n-1}(\B ; K)$ déduite de la suite exacte courte $0\to K\to H_0(\A ; P)\to H_0(\A ; F)\to 0$ permet de voir que la ligne inférieure du diagramme est exacte.

En notant $H_i(A)$ pour $H_i(\A\times\B ; A)$ et $gr_iH(A)$ pour $\underset{r+s=i}{\bigoplus} H_r(\B ; H_s(\A ; A))$ pour alléger, le lemme des cinq appliqué au diagramme commutatif aux lignes exactes (déduit de l'hypothèse de récurrence)
$$\xymatrix{H_n(P)\ar[d]^\simeq_h\ar[r]^a & H_n(F)\ar[d]^f\ar[r]^b & H_{n-1}(N)\ar[d]^\simeq_g\ar[r] & H_{n-1}(P)\ar[d]^\simeq \\
gr_nH(P)\ar[r]^c & gr_nH(F)\ar[r]^-d  & gr_{n-1}H(N)\ar[r] & gr_{n-1}H(P)
}$$
montre que $f$ est surjectif.

On considère à présent un diagramme
$$\xymatrix{H_n(N)\ar[r]\ar[d]^{f'} & H_n(P)\ar[d]^\simeq_h\ar[r]^a & H_n(F)\ar[d]^f\ar[r]^b & H_{n-1}(N)\ar[d]^\simeq_g\\
gr_nH(N)\ar[r] & gr_nH(P)\ar[r]^c & gr_nH(F)\ar[r]^-d  & gr_{n-1}H(N)
}$$
dans lequel la flèche $f'$ est construite comme la flèche $f$, en remplaçant $F$ par $N$. Ce diagramme est commutatif et ses lignes sont exactes, pour des raisons similaires aux précédentes. Ce qu'on vient d'établir montre que la flèche $f'$ est surjective. Le lemme des cinq montre alors que $f$ est injective, ce qui achève la démonstration.
\end{proof}

En particulier, sous les hypothèses de la proposition~\ref{arret2}, on obtient une suite spectrale fonctorielle
$$E^2_{p,q}(F)=H_p(G_\infty ; H_q(\C ; F))\Rightarrow H_{p+q}(\C\times G_\infty ; \Pi^*F)\simeq H_{p+q}(G_\infty ; F_\infty),$$
où $G_\infty$ opère trivialement sur $H_q(\C ; F)$, qui s'effondre à la deuxième page et induit un scindement a priori non naturel
$$H_n(G_\infty ; F_\infty)\simeq\bigoplus_{p+q=n}H_p(G_\infty ; H_q(\C ; F)).$$

\begin{rem} L'hypothèse que $\C$ possède un objet initial n'est pas nécessaire pour cette suite spectrale. On notera par ailleurs que, même dans les cas usuels favorables, il n'existe généralement pas de scindement {\em naturel} de la filtration associée à cette suite spectrale (contrairement à ce que la rédaction de l'article \cite{BP} peut laisser penser).
\end{rem}

Les suites spectrales que nous venons d'étudier, qui convergent vers $H_*(G_\infty;F_\infty)$, qui s'arrêtent au terme $E^2$ dans les cas favorables, méritent d'être comparées à celle obtenue en $K$-théorie stable (et avec ses variantes faisant intervenir d'autres groupes que le groupe linéaire).

Rappelons sa construction. Soit $G$ un groupe dont le sous-groupe des commutateurs est parfait (c'est le cas du groupe linéaire infini sur un anneau arbitraire, mais aussi du groupe symétrique infini, du groupe orthogonal ou symplectique infini sur un corps commutatif). On applique la construction plus de Quillen au classifiant $BG$ de $G$ (comme groupe discret) et on forme la fibre homotopique $F_G$ de l'application canonique $BG\to BG^+$.

L'homologie de $F_G$ à coefficients dans un $G$-module $M$, vu comme $\pi_1(F_G)$-module via le morphisme canonique $\pi_1(F_G)\to\pi_1(BG)=G$, est par définition la $K$-théorie stable $G$-généralisée à coefficients dans $M$ ; on la note $K^{s}_*(G;M)$. Le point remarquable est le suivant : dans la suite spectrale de Serre
\begin{equation}\label{ssks}
E^2_{p,q}=H_p(BG^+;K^s_q(G;M))\simeq H_p(G;K^s_q(G;M))\Rightarrow H_{p+q}(G;M)
\end{equation}
l'action du groupe $G$ sur le groupe abélien $K^s_q(G;M)$ est {\em triviale} (cf. par exemple \cite{Ka}, théorème~$3.1$).

De plus, dans de nombreux cas, la suite spectrale s'effondre au terme $E^2$ et la filtration associée se scinde (de manière généralement {\em non naturelle}) : voir l'article \cite{BP} et l'article "Stable $K$-theory is bifunctor homology" de Franjou et Pirashvili, déjà évoqués, pour le cas classique du groupe linéaire (nous n'avons d'ailleurs fait que reprendre leur démonstration pour établir la proposition~\ref{eff-kun}) ; ce cas est adapté dans \cite{Bet-sym}, théorème~$1.3$, pour le groupe symétrique (qui comme l'article \cite{BP} peut suggérer un scindement naturel, incorrect, de la graduation). 

La $K$-théorie stable comme l'homologie de catégorie rendent donc généralement des services très analogues en terme de calcul de l'homologie du groupe $G$ à coefficients tordus : dans les cas où l'on sait les identifier naturellement, la suite spectrale~(\ref{ssks}) et celle de la proposition~\ref{eff-kun} sont isomorphes. Néanmoins, nous ne connaissons pas d'analogue à la suite spectrale de la proposition~\ref{fctder} (qui présente l'avantage non seulement de s'effondrer à la deuxième page, mais aussi de procurer une décomposition {\em fonctorielle} de l'homologie stable dans les cas favorables, contrairement à celle de la proposition~\ref{eff-kun}). 

Il convient d'ailleurs de noter que les arguments de Scorichenko pour identifier la $K$-théorie stable d'un anneau quelconque et l'homologie de la catégorie des modules projectifs de type fini pour des coefficients déduits de bifoncteurs polynomiaux n'utilise pas réellement la définition de la $K$-théorie stable, mais seulement l'existence de la suite spectrale naturelle (\ref{ssks}) et d'un morphisme naturel gradué $K^s_*(GL_\infty;M)\to H_*(GL_\infty;M)$ (induit par l'application canonique $F_G\to BG$). On déduit de ce morphisme un morphisme naturel de la $K$-théorie stable à coefficients provenant d'un bifoncteur vers l'homologie de ce bifoncteur, morphisme dont Scorichenko montre qu'il est bijectif si le bifoncteur est polynomial en exploitant uniquement la suite spectrale (\ref{ssks}). Dans le cas d'un corps fini, l'homologie du groupe linéaire infini à coefficients dans le corps est trivial, comme l'a montré Quillen dans \cite{Qui}, de sorte que la suite spectrale~(\ref{ssks}) se réduit à un isomorphisme naturel entre $K$-théorie stable et homologie de $GL_\infty$, c'est pourquoi le théorème de Scorichenko sur la $K$-théorie stable se réduit dans ce cas, traité auparavant par Betley (cf. \cite{Bet}) et Suslin (cf. l'appendice de \cite{FFSS}) indépendamment, à une identification entre homologie du groupe linéaire (à coefficients tordus) et homologie de catégorie.

Nous reviendrons dans les appendices~\ref{appBS} et~\ref{apsym} sur l'équivalence des méthodes de calcul d'homologie stable des groupes linéaires et symétriques respectivement, à coefficients venant d'un foncteur convenable, à l'aide de la $K$-théorie stable (généralisée au groupe symétrique infini dans le second cas) ou à l'aide de notre suite spectrale utilisant directement l'homologie d'une catégorie appropriée.

%% file: dv-s3.tex
Cette section fournit le r\'esultat principal de cet article, qui relie l'homologie stable des groupes orthogonaux (resp. symplectiques) \`a coefficients tordus et des groupes de torsion sur la cat\'egorie des espaces vectoriels. Ces derniers sont accessibles au calcul, dans les cas favorables, comme nous l'illustrons \`a la section~\ref{scal}.

On commence par identifier une certaine cat\'egorie de fractions de $\Eqd$ (resp. de $\E_{alt}^{deg}$) \`a une cat\'egorie de Burnside (cf. définition~\ref{span}). La décomposition des foncteurs de Mackey non additifs associés à cette catégorie (cf. théorème~\ref{equ-span}) intervient pour simplifier la deuxième page de la suite spectrale de la section précédente pour les groupes orthogonaux ou symplectiques. Combiné à des résultats d'annulation en homologie des foncteurs, cela permet d'obtenir notre résultat central, le théorème~\ref{thf-o}.

\subsection{Les catégories de fractions $\Eqd[(-\perp \mathbf{H})^{-1}]$ et $\E_{alt}^{deg}[(-\perp \mathbf{H})^{-1}]$}\label{par-fract}

Dans le cas o\`u $\C=\Eqd(k)$ (resp. $\C=\E_{alt}^{deg}(k)$), où $k$ est un corps commutatif fixé, le foncteur
 $\mathbb{L}_n Q_!(\kk)$ transforme l'inclusion canonique $D \to S_{\mathbf{H}}(i)\perp D= \mathbf{H}^{\perp i} \perp D$ en un isomorphisme d'apr\`es la proposition~\ref{fdo}. Cette observation justifie que l'on s'int\'eresse dans ce paragraphe \`a la cat\'egorie de fractions inversant les morphismes canoniques de la forme $D \xrightarrow{i} D \perp \mathbf{H}^{\perp n}$.
 
Soient $(-\perp \mathbf{H})$ l'ensemble de fl\`eches de $\Eqd$ suivant:
$$(-\perp \mathbf{H})=\{ D \xrightarrow{i} D \perp H\ \mid H \mathrm{ \ espace\  hyperbolique\ }; i \mathrm{ \ inclusion\ canonique} \}$$
 et $\Eqd[(-\perp \mathbf{H})^{-1}]$ la cat\'egorie de fractions correspondante (qui existe par des r\'esultats g\'en\'eraux, voir par exemple \cite{Gabriel-Zisman}). On note  $\phi: \Eqd \to \Eqd[(-\perp \mathbf{H})^{-1}]$ le foncteur canonique. Le but de ce paragraphe est de d\'emontrer qu'il existe une \'equivalence de cat\'egories 
$$\Psi: \Eqd[(-\perp \mathbf{H})^{-1}] \xrightarrow{\simeq} Sp(\E^f_{inj})$$
o\`u $\E^f_{inj}$ est la sous-catégorie des injections de $\E^f(k)$ (cf. définition~\ref{def-cataux}) et $Sp(.)$ d\'esigne la cat\'egorie de Burnside dont on rappelle la construction dans la définition~\ref{span} de l'appendice~\ref{amm}.

Apr\`es avoir donn\'e quelques r\'esultats essentiels sur les morphismes de $\Eqd[(-\perp \mathbf{H})^{-1}]$, on d\'efinira le foncteur $\Psi$, puis on montrera que $\Psi$ est essentiellement surjectif, plein et fid\`ele.

\begin{nota}
Dans les diagrammes dans la cat\'egorie de fractions $\Eqd[(-\perp \mathbf{H})^{-1}]$ on notera par 
\xymatrix{\ar@{~>}[r] & } les morphismes de $\Eqd$ qui sont \'el\'ements de $(-\perp \mathbf{H})$ et par $\to$ tout  morphisme de $\Eqd$.
\end{nota}

Remarquons que $(-\perp \mathbf{H})$ est stable par composition et que pour tout objet $D$ de $\Eqd$, $Id_{D} \in (-\perp \mathbf{H})$. On a \'egalement la propri\'et\'e suivante:
\begin{lm}
Soient $H$ un espace hyperbolique et $f: D \to D'$ un morphisme de $\Eqd$, il existe $g: D \perp H \to D' \perp H$ rendant le diagramme suivant commutatif:
$$\xymatrix{
D \ar@{~>}[r] \ar[d]_f& D \perp H\ar[d]^g\\
D' \ar@{~>}[r] &D' \perp H.
}$$
\end{lm}

\begin{proof}
Ce diagramme est  commutatif pour  $g=f \perp Id_{H}$.
\end{proof}

Ce lemme nous permet de donner la description suivante des morphismes de $\Eqd[(-\perp \mathbf{H})^{-1}]$.
\begin{lm} \label{decompo}
Tout morphisme de $\Eqd[(-\perp \mathbf{H})^{-1}]$ s'\'ecrit sous la forme $g^{-1} f$ o\`u $g$ est un \'el\'ement de $(-\perp \mathbf{H})$ et $f$ est un morphisme de $\Eqd$.
\end{lm}

\begin{rem}
L'ensemble $(-\perp \mathbf{H})$ n'admet pas un "calcul \`a gauche des fractions" au sens de Gabriel et Zisman (\cite{Gabriel-Zisman}). En effet, 
pour $\xymatrix{ D \ar@{~>}[r]^-i& D \perp H}$ un \'el\'ement de $(-\perp \mathbf{H})$, $f: D \perp H \to D'$ et $g: D \perp H \to D'$ deux morphismes de $\Eqd$ tels que $fi=gi$, il n'existe pas de morphisme de $(-\perp \mathbf{H})$ $\xymatrix{ D' \ar@{~>}[r]^-{i'}& D' \perp K}$ tel que $i'f=i'g$. Nous n'avons donc pas une description de $\Eqd[(-\perp \mathbf{H})^{-1}]$ aussi simple que dans le cas consid\'er\'e par Gabriel-Zisman, il est donc a priori difficile de savoir quand deux morphismes sont égaux dans cette cat\'egorie, l'\'ecriture d'un morphisme de $\Eqd[(-\perp \mathbf{H})^{-1}]$ sous la forme $g^{-1} f$ n'\'etant pas unique. N\'eanmoins nous avons les r\'esultats suivants qui seront essentiels dans la preuve de la fid\'elit\'e du foncteur $F$ du th\'eor\`eme~\ref{thm-fractions}.
\end{rem}

\begin{lm} \label{auto}
Soient $D$ un objet de $\Eqd$,  $H$ un espace hyperbolique et $f \in {\rm Aut}_{\Eqd}(D \perp H)$ tel que $f_{\mid D}=Id_D$ alors $f=Id_{D \perp H}$ dans $\Eqd[(-\perp \mathbf{H})^{-1}]$.
\end{lm}
\begin{proof}
Soient $D,\ H$ et $f$ comme dans l'\'enonc\'e. L'inclusion canonique $D \xrightarrow{i} D \perp H$ est un \'el\'ement de $(-\perp \mathbf{H})$ et est donc un isomorphisme dans $\Eqd[(-\perp \mathbf{H})^{-1}]$. Or on a: $fi=i$ dont on d\'eduit que $f=Id_{D \perp H}$ dans $\Eqd[(-\perp \mathbf{H})^{-1}]$.
\end{proof}

Le cas particulier o\`u $D$ est l'espace vectoriel nul fournit le lemme suivant:
\begin{lm} \label{lm1}
Soient $H$ un espace hyperbolique et $f \in {\rm Aut}_{\Eqd}(H)$. On a alors $f=Id_H$ dans $\Eqd[(-\perp \mathbf{H})^{-1}]$.
\end{lm}

On d\'eduit de ce lemme le r\'esultat suivant:
\begin{pr} \label{prop1}
Soient $D$ un objet de $\Eqd$,  $H$ un espace hyperbolique et $f,g \in Hom_{\Eqd}(D,H)$ alors $f=g$ dans $\Eqd[(-\perp \mathbf{H})^{-1}]$.
\end{pr}
\begin{proof}
Par le th\'eor\`eme de Witt, il existe $h \in  {\rm Aut}_{\Eqd}(H)$ tel que $hf=g$. Le lemme~\ref{lm1} permet d'en d\'eduire que $f=g$ dans $\Eqd[(-\perp \mathbf{H})^{-1}]$.
\end{proof}

Nous aurons \'egalement besoin dans la suite du lemme tr\`es facile suivant:
\begin{lm} \label{memeH}
Soient $\alpha, \beta: V \to W$ dans $\Eqd[(-\perp \mathbf{H})^{-1}]$ alors il existe un espace hyperbolique $H$  tel que $\alpha=i^{-1} f$ et $\beta = i^{-1} g$ avec $f,g: V \to W \perp H$ et $i: W \to W \perp H$ l'inclusion canonique.
\end{lm}

\begin{proof}
D'apr\`es le lemme~\ref{decompo} on a: $\alpha=i_1^{-1} f_1$ et $\beta = i_2^{-1} g_2$ o\`u $f_1: V \to W \perp H_1$, $g_2: V \to W \perp H_2$ et $i_1: W \to W \perp H_1$, $i_2: W \to W \perp H_2$ sont  les inclusions canoniques. Il suffit alors de v\'erifier que $\alpha=i_2^{-1} i_2i_1^{-1} f_1= i^{-1} f$ o\`u $f: V \to W \perp H_1 \perp H_2$ est la compos\'ee de $f_1$ et de l'inclusion canonique $W \perp H_1\to W \perp H_1 \perp H_2$. De m\^eme pour $\beta$.
\end{proof}

\begin{nota}
Dans la suite, un morphisme de $\Eqd[(-\perp \mathbf{H})^{-1}]$ se d\'ecomposant sous la forme $\alpha=i^{-1} f$  o\`u $f: V \to W \perp H$ et $i: W \to W \perp H$ l'inclusion canonique, sera not\'e:
$$\xymatrix{  V \ar[r]^-f & W \perp H & W  \ar@{~>}[l]}.$$

\end{nota}

\begin{rem}
La cat\'egorie $\Eqd[(-\perp \mathbf{H})^{-1}] $ est \'equivalente \`a la cat\'egorie de fractions de $\Eqd$ o\`u l'on inverse l'ensemble des inclusions canoniques $D \rightarrow D \perp K$ o\`u $K$ est un espace quadratique non d\'eg\'en\'er\'e. En effet, comme tout espace quadratique non d\'eg\'en\'er\'e $K$ se plonge dans un espace hyperbolique $K \perp K'$, la compos\'ee suivante:
$$ V \xrightarrow{a} V \perp K \xrightarrow{b} V \perp K \perp K'$$
est inversible dans $\Eqd[(-\perp \mathbf{H})^{-1}] $, ce qui implique que $a$ est inversible \`a gauche et $b$ est inversible \`a droite dans $\Eqd[(-\perp \mathbf{H})^{-1}] $. De plus, comme la compos\'ee
$$ V \perp K \xrightarrow{b} V \perp K \perp K'  \xrightarrow{c} V \perp K \perp K' \perp K$$
est inversible dans  $\Eqd[(-\perp \mathbf{H})^{-1}] $, $b$ est inversible \`a gauche dans  $\Eqd[(-\perp \mathbf{H})^{-1}] $. On en d\'eduit que $b$ est inversible dans  $\Eqd[(-\perp \mathbf{H})^{-1}] $  et  donc que $a$ l'est \'egalement.
\end{rem}

Afin de d\'efinir le foncteur de $\Eqd$ dans $Sp(\E^f_{inj})$ dont d\'ecoule l'\'equivalence de cat\'egories annonc\'ee en d\'ebut de section nous avons besoin du lemme facile suivant.

\begin{lm} \label{lm}
Soit $f \in {\rm Hom}_{\Eqd}(D,D')$, on a $f^{-1}(Rad(D')) \subset Rad(D)$ o\`u $Rad$ d\'esigne le radical.
\end{lm}

\begin{pr} \label{tildeF}
Il existe un foncteur $\widetilde{\Psi}: \Eqd \to Sp(\E^f_{inj})$ d\'efini sur les objets par $\widetilde{\Psi}(D)=Rad(D)$ et sur les morphismes par:
$$\widetilde{\Psi}(D \xrightarrow{f} D')=Rad(D) \xleftarrow{i} f^{-1}(Rad(D')) \xrightarrow{\widetilde f} Rad(D')$$
o\`u $i$ est l'inclusion obtenue dans le lemme~\ref{lm}  et $\widetilde{f}$ est la restriction de $f$ \`a $f^{-1}(Rad(D'))$.
\end{pr}
\begin{proof}
V\'erifions que $\widetilde{\Psi}$ d\'efinit bien un foncteur.

Pour $D \in \Eqd$:
$\widetilde{\Psi}(D \xrightarrow{Id} D)= Rad(D) \xleftarrow{Id} Rad(D) \xrightarrow{Id} Rad(D)$.

Pour une paire de fl\`eches composables de $\Eqd$: $D \xrightarrow{f} D' \xrightarrow{g} D''$ on a:

$$\widetilde{\Psi}(g \circ f)=Rad(D) \hookleftarrow (g \circ f)^{-1} (Rad(D'')) \xrightarrow{\widetilde{g \circ f}} Rad(D'')$$

et $\widetilde{\Psi}(g) \circ \widetilde{\Psi}(f)$ est donn\'e par le diagramme en escalier suivant:
$$\xymatrix{
f^{-1}( g^{-1}(Rad(D'')) \ar@{^{(}->}[d]  \ar[r]^-f & g^{-1}(Rad(D'')) \ar@{^{(}->}[d]  \ar[r]^-{\widetilde g} & Rad(D'')\\
f^{-1}(Rad(D'))  \ar@{^{(}->}[d] \ar[r]_{\widetilde f} & Rad(D')\\
Rad(D).
}$$
D'o\`u 
$$\widetilde{\Psi}(g \circ f)=\widetilde{\Psi}(g) \circ \widetilde{\Psi}(f).$$
\end{proof}

\begin{pr}
Le  foncteur $\widetilde{\Psi}: \Eqd \to Sp(\E^f_{inj})$ induit un unique foncteur $\Psi:  \Eqd[(-\perp \mathbf{H})^{-1}] \to Sp(\E^f_{inj})$  rendant le diagramme suivant commutatif :

$$\xymatrix{
 \Eqd \ar[r]^{\widetilde{\Psi}}  \ar[d]_{\phi}& Sp(\E^f_{inj})\\
  \Eqd[(-\perp \mathbf{H})^{-1}]. \ar@{-->}[ur]_{\Psi}
}$$
\end{pr}
\begin{proof}
Soit $ D \xrightarrow{i} D \perp H$ un \'el\'ement de $(-\perp \mathbf{H})$, on a 
$$\widetilde{\Psi}(i)=Rad(D) \xleftarrow{Id} Rad(D) \xrightarrow{Id} Rad(D)$$
(qui est bien un isomorphisme de $Sp(\E^f_{inj})$!). 

Par la propri\'et\'e universelle de la cat\'egorie des fractions, il existe un unique foncteur $\Psi: \Eqd[(-\perp \mathbf{H})^{-1}] \to Sp(\E^f_{inj})$ tel que $\Psi \circ \phi=\widetilde{\Psi}$.
\end{proof}

Afin de montrer que le foncteur $\Psi$ fournit une \'equivalence de cat\'egories, on prouve dans la suite qu'il est essentiellement surjectif et pleinement fid\`ele.

\begin{pr}
Le  foncteur $\Psi: \Eqd[(-\perp \mathbf{H})^{-1}] \to Sp(\E^f_{inj})$ est essentiellement surjectif.
\end{pr}
\begin{proof}
Soit $V \in \E^f_{inj}$, on munit $V$ de la forme quadratique nulle. L'espace quadratique $D$ ainsi obtenu \'etant totalement isotrope on a $\Psi(D)=V$.  
\end{proof}

\begin{pr}
Le  foncteur $\Psi: \Eqd[(-\perp \mathbf{H})^{-1}] \to Sp(\E^f_{inj})$ est plein.
\end{pr}
\begin{proof}
Soient $V=Rad(V) \perp H$ et $W=Rad(W) \perp K$ deux objets de $\Eqd[(-\perp \mathbf{H})^{-1}] $ et $f \in {\rm Hom}_{Sp(\E^f_{inj})}(Rad(V), Rad(W))$. Par d\'efinition des morphismes de $Sp(\E^f_{inj})$ il existe un sous-espace vectoriel $X$ de $Rad(V)$ tel que $f=Rad(V) \xleftarrow{i} X \xrightarrow{\beta} Rad(W)$ o\`u $i$ est l'inclusion  et $\beta$ est une injection lin\'eaire. 

Soit $X'$ un espace vectoriel suppl\'ementaire de $X$ dans $Rad(V)$ et $X''$ un suppl\'ementaire de $\beta(X)$ dans $Rad(W)$. Comme les formes quadratiques sur les radicaux sont nulles on a:
$$Rad(V) \simeq X \perp X' \qquad Rad(W) \simeq \beta(X) \perp X''.$$
Par les propri\'et\'es g\'en\'erales de la cat\'egorie $\Eqd$ il existe un espace hyperbolique $L$  et $g \in {\rm Hom}_{\Eqd}(X' ,L)$.

Soit $k: V \simeq X \perp X' \perp H \to (\beta(X) \perp X''  \perp K )\perp H \perp L \simeq W \perp H \perp L $ l'application donn\'ee par la matrice:
$$\left(\begin{array}{ccc} 
\beta & 0 &0 \\
0 &0 &Id_H \\
0 &g & 0
\end{array}\right).$$

L'application $u= \xymatrix{V  \ar[r]^-k & W \perp H \perp L & W \ar@{~>}[l] } $ de $\Eqd[(-\perp \mathbf{H})^{-1}] $ v\'erifie $\Psi(u)=f$.
\end{proof}

\begin{pr}
Le  foncteur $\Psi: \Eqd[(-\perp \mathbf{H})^{-1}] \to Sp(\E^f_{inj})$ est fid\`ele.
\end{pr}
\begin{proof}
Soient $\alpha, \beta \in Hom_{\Eqd[(-\perp \mathbf{H})^{-1}] }(V,W)$ tels que $\Psi(\alpha)=\Psi(\beta)$. D'apr\`es les lemmes~\ref{decompo} et~\ref{memeH} on peut \'ecrire $\alpha$ et $\beta$ sous la forme:
$$\alpha=\xymatrix{V \ar[r]^-f  & W \perp H & W  \ar@{~>}[l] }; \quad \beta=\xymatrix{V \ar[r]^-g & W \perp H & W  \ar@{~>}[l]   }.$$
Par la d\'efinition de la cat\'egorie de Burnside $Sp(\E^f_{inj})$ rappel\'ee en~\ref{span}, l'\'egalit\'e $\Psi(\alpha)=\Psi(\beta)$ implique que $f^{-1}(Rad(W))=g^{-1}(Rad(W))$ et $\widetilde{f}=\widetilde{g}$.

On a les d\'ecompositions suivantes des espaces $V$ et $W \perp H$: $V=f^{-1}(Rad(W)) \perp D $ o\`u $D$ est un objet de $\Eqd$; $W \perp H = Rad(W) \perp L$ o\`u $L$ est un espace non d\'eg\'en\'er\'e. Dans des bases de $V$ et $W \perp H$ obtenues en juxtaposant des bases de $f^{-1}(Rad(W))$ et $\ D $ pour $V$ et de $Rad(W)$ et $L$ pour $W \perp H$ la matrice de l'application $f$ s'\'ecrit:

$$\left(\begin{array}{cc} 
\widetilde f & \epsilon \\
0& i 
\end{array} \right)$$
o\`u $i: D \to L$ est une injection (pr\'eservant la forme quadratique) puisque $D \cap f^{-1}(Rad(W))=0$,  et $\epsilon: D \to Rad(W)$ est une application lin\'eaire. La matrice de $g$ dans la m\^eme base est de la forme 
$$\left(\begin{array}{cc} 
\widetilde f& \epsilon' \\
0& i' 
\end{array} \right)$$
avec les m\^emes conditions sur $i'$ et $ \epsilon'$ que pr\'ec\'edemment.

Comme $i$ et $i'$ sont des injections pr\'eservant les formes quadratiques, \`a valeurs dans un espace quadratique non  d\'eg\'en\'er\'e, par le th\'eor\`eme de Witt, il existe $u \in O(L)$ tel que $ui=i'$. De plus, comme $i$ est injective, il existe $\alpha: L \to Rad(W)$ tel que $\epsilon+ \alpha i =\epsilon'$.

Consid\'erons l'automorphisme $l$ de $Rad(W) \perp L$ de matrice 

$$
\left(\begin{array}{cc} 
Id & \alpha\\
0& u
\end{array} \right)
$$

l'application $lf$ a pour matrice

$$
\left(\begin{array}{cc} 
f_{\mid f^{-1}(Rad(W))} & \epsilon+ \alpha i\\
0& ui
\end{array} \right)
=
\left(\begin{array}{cc} 
f_{\mid f^{-1}(Rad(W))} & \epsilon' \\
0& i' 
\end{array} \right)$$

dont on d\'eduit que $lf=g$. Or, par le lemme~\ref{auto}, on a $l=Id_{Rad(W) \perp L}$ d'o\`u: $f=g$ dans $ \Eqd[(-\perp \mathbf{H})^{-1}]$.
\end{proof}

Nous avons donc d\'emontr\'e le r\'esultat suivant:

 \begin{thm} \label{thm-fractions}
 Le foncteur $\widetilde{\Psi}: \Eqd \to Sp(\E^f_{inj})$ d\'efini \`a la proposition~\ref{tildeF} induit 
 une \'equivalence de cat\'egories 
$$\Psi: \Eqd[(-\perp \mathbf{H})^{-1}] \xrightarrow{\simeq} Sp(\E^f_{inj})$$
o\`u $\E^f_{inj}$ est la cat\'egorie ayant pour objets les $k$-espaces vectoriels de dimension finie et pour morphismes les applications lin\'eaires injectives. 
\end{thm}

Pour la cat\'egorie $\E_{alt}^{deg}$, on d\'emontre de mani\`ere similaire le th\'eor\`eme suivant.

\begin{thm} \label{thm-fractions-symp}
 Le foncteur $\widetilde{\Psi}: \E_{alt}^{deg} \to Sp(\E^f_{inj})$ d\'efini comme \`a la proposition~\ref{tildeF} mutatis mutandis induit 
 une \'equivalence de cat\'egories 
$$\Psi: \E_{alt}^{deg}[(-\perp \mathbf{H})^{-1}] \xrightarrow{\simeq} Sp(\E^f_{inj})$$
o\`u $\mathbf{H}$ d\'esigne l'espace symplectique d\'efini en~\ref{exfond}.\ref{exsymp}. 
\end{thm}

\subsection{Théorème fondamental} \label{thm-fond}

On commence par donner quelques notations utilisées couramment dans la suite :
\begin{nota}\label{not-gi} Soit $k$ un corps.
\begin{enumerate}
\item On désigne par $O_\infty(k)$ (resp. $Sp_\infty(k)$) le groupe $\underset{n\in\mathbb{N}}{\col}O_{n,n}(k)$ (resp. $\underset{n\in\mathbb{N}}{\col}Sp_{2n}(k)$), la colimite étant prise conformément aux conventions de la section~\ref{s-un}.
\item On désigne par $S^i$ (resp. $\Gamma^i$, $\Lambda^i$) l'endofoncteur $i$-ème puissance symétrique (resp. divisée, extérieure) des $k$-espaces vectoriels.
\item On note $\F(k)$, voire $\F$, la catégorie $\E^f(k)-\mathbf{Mod}$ des foncteurs depuis les $k$-espaces vectoriels de dimension finie vers les $\kk$-modules.
\item Lorsque $F$ est un foncteur depuis une catégorie d'espaces vectoriels, on note $F^\vee$ la précomposition de $F$ par le foncteur de dualité $(-)^*$ : $F^\vee(V)=F(V^*)$.
\end{enumerate}
\end{nota}
(On rappelle que $F^*$ désigne pour sa part, conformément à la notation de l'appendice~\ref{appA}, la {\em post}composition de $F$ par la dualité, lorsque $\kk$ est un corps.)

La notation $O_\infty$ est justifiée, pour les considérations homologiques qui sont les nôtres, par le fait que la colimite sur tous les groupes orthogonaux donne des résultats canoniquement isomorphes (même si le groupe orthogonal infini obtenu n'est pas isomorphe à celui qu'on considère) --- cf. remarque~\ref{rqsco}.\,\ref{qand}.

On aura également besoin de la notion classique suivante :
\begin{defi}[Cf. \cite{FFPS}] \label{poly}Soit $k$ un corps. 
\begin{enumerate}
\item Le {\em foncteur différence} de $\F(k)$ est l'endofoncteur $\Delta$ de cette catégorie d\'efini comme \'etant le noyau de l'épimorphisme scindé évident $(-\oplus k)^*\to Id$.
\item Un foncteur $F$ de $\F(k)$ est dit {\em polynomial} s'il existe $n\in\mathbb{N}$ tel que $\Delta^n(F)=0$ ; son {\em degré} est alors le plus grand entier $d$ tel que $\Delta^d(F)\neq 0$.
\item Un foncteur {\em analytique} est une colimite (qu'on peut supposer filtrante) de foncteurs polynomiaux.
\end{enumerate}
\end{defi}

Le théorème ci-dessous constitue le résultat principal de cet article.

\begin{thm}\label{thf-o}
Soient $k$ un corps fini et $F$ un foncteur analytique de $\F(k)$.

Il existe des isomorphismes naturels gradués
$$H_*(O_\infty(k) ; F_\infty)\simeq {\rm Tor}_*^{\E^f(k)\times O_\infty(k)}(\kk[S^2]^\vee,F)$$
et
$$H_*(Sp_\infty(k) ; F_\infty)\simeq {\rm Tor}_*^{\E^f(k)\times Sp_\infty(k)}(\kk[\Lambda^2]^\vee,F)$$
(où les groupes $O_\infty(k)$ et $Sp_\infty(k)$ agissent trivialement).
\end{thm}

\begin{cor}\label{thfvdec}
Sous les hypothèses précédentes, il existe des suites spectrales naturelles en $F$
données par
$$E^2_{p,q}={\rm Tor}^{\E^f(k)}_p(H_q(O_\infty(k) ;\kk)\underset{\kk}{\otimes}\kk[S^2]^\vee,F)\Rightarrow H_{p+q}(O_\infty(k) ; F_\infty)$$
et
$$E^2_{p,q}={\rm Tor}^{\E^f(k)}_p(H_q(Sp_\infty(k) ;\kk)\underset{\kk}{\otimes}\kk[\Lambda^2]^\vee,F)\Rightarrow H_{p+q}(Sp_\infty(k) ; F_\infty).$$

Lorsque $\kk$ est un anneau de dimension homologique au plus $1$, elles s'effondrent au terme $E^2$ et procurent des isomorphismes canoniques
$$H_n(O_\infty(k) ; F_\infty)\simeq\bigoplus_{p+q=n}{\rm Tor}^{\E^f(k)}_p(H_q(O_\infty(k);\kk)\underset{\kk}{\otimes}\kk[S^2]^\vee,F)$$
et
$$H_n(Sp_\infty(k) ; F_\infty)\simeq\bigoplus_{p+q=n}{\rm Tor}^{\E^f(k)}_p(H_q(Sp_\infty(k);\kk)\underset{\kk}{\otimes}\kk[\Lambda^2]^\vee,F).$$
\end{cor}

\begin{cor}\label{thfssv2}
Sous les hypothèses du théorème, il existe des suites spectrales naturelles
$$E^2_{p,q}=H_p(O_\infty(k) ; {\rm Tor}^{\E^f(k)}_q(\kk[S^2]^\vee,F))\Rightarrow H_{p+q}(O_\infty(k) ; F_\infty)$$
et
$$E^2_{p,q}=H_p(Sp_\infty(k) ; {\rm Tor}^{\E^f(k)}_q(\kk[\Lambda^2]^\vee,F))\Rightarrow H_{p+q}(Sp_\infty(k) ; F_\infty)$$
qui s'effondrent à la deuxième page lorsque l'anneau $\kk$ est de dimension au plus~$1$ et procurent alors des isomorphismes non fonctoriels
$$H_n(O_\infty(k) ; F_\infty)\simeq\bigoplus_{p+q=n}H_p(O_\infty(k) ; {\rm Tor}^{\E^f(k)}_q(\kk[S^2]^\vee,F))$$
et
$$H_n(Sp_\infty(k) ; F_\infty)\simeq\bigoplus_{p+q=n}H_p(Sp_\infty(k) ; {\rm Tor}^{\E^f(k)}_q(\kk[\Lambda^2]^\vee,F)).$$
\end{cor}

Avant de démontrer ces résultats, on donne quelques d\'efinitions puis un lemme reliant les différents foncteurs introduits ainsi que ceux de l'appendice~\ref{sect-cof}. Les mentions au corps $k$ sont sous-entendues.

\begin{defi} On désigne par $\E^f_{iso}$ la sous-catégorie de $\E^f$ ayant les mêmes objets et les isomorphismes pour morphismes.
\end{defi}

Dans ce qui suit, on utilise également les catégories introduites dans la définition~\ref{def-cataux}.

\begin{defi}\begin{itemize}
 \item Le foncteur $\beta: \mathbf{Mod}-\E^f_{iso}\to\mathbf{Mod}-\E^f_{surj}$ envoie un foncteur $F$ sur le foncteur prenant les mêmes valeurs sur les objets et les isomorphismes et envoyant les surjections strictes sur $0$.
\item Le foncteur  $\gamma : \mathbf{Mod}-\Eqd[(-\oplus \mathbf{H})^{-1}] \to\mathbf{Mod}-\Eqd$  est la précomposition par le foncteur canonique  $\phi: \Eqd\to\Eqd[(-\oplus \mathbf{H})^{-1}]$.
\item Le foncteur $\delta : \mathbf{Mod}-\E^f\to\mathbf{Mod}-\E^f_{inj}$  est la précomposition par le foncteur d'inclusion $\E^f_{inj}\to\E^f$.
\end{itemize}
\end{defi}

\begin{rem}
On rappelle qu'une  forme quadratique sur un espace vectoriel $V$ est un \'el\'ement de $S^2(V^*)$. Par  cons\'equent, pour  $k$ un corps commutatif (éventuellement de caractéristique $2$), on a l'\'egalit\'e: $$\E_q^{deg}=(\E^f_{inj})_{(S^2)^{\vee}}$$ o\`u $(\E^f_{inj})_{(S^2)^\vee}$ est la cat\'egorie dont les objets sont les couples $(V,\alpha)$ o\`u $V$ est un objet de $\E^f_{inj}$ et $\alpha$ est un \'el\'ement de $S^2(V^{*})$. On note  $\Omega_{(S^2)^\vee}: \mathbf{Mod}-\E_q^{deg}=\mathbf{Mod}-(\E^f_{inj})_{(S^2 )^{\vee}}\to \mathbf{Mod}-\E^f_{inj}$ le foncteur  d\'efini \`a la proposition~\ref{torcomod} de l'appendice~\ref{appA}.
\end{rem}

La preuve du th\'eor\`eme~\ref{thf-o} repose sur le lemme suivant, qui combine les résultats du §\,\ref{par-fract} et la décomposition des foncteurs de Mackey. On rappelle que les foncteurs $\rho$, $\kappa$ et $\omega$ sont introduits dans la définition~\ref{def-cataux}.

\begin{lm}\label{dcm}
Le diagramme suivant :
$$\xymatrix{\mathbf{Mod}-\E^f_{iso}\ar[d]_{\beta} \ar[rr]^-{\xi}_-{\simeq} & & \mathbf{Mod}-Sp(\E^f_{inj}) \ar[r]^-\eta_-\simeq & \mathbf{Mod}-\Eqd[(-\oplus \mathbf{H})^{-1}]\ar[r]^-{\gamma} & \mathbf{Mod}-\Eqd\ar[d]^{\Omega_{(S^2)^\vee}} \\
 \mathbf{Mod}-\E^f_{surj}\ar[rr]^-{\rho\otimes\kappa(\kk[S^2]^\vee)} & & \mathbf{Mod}-\E^f_\Gr\ar[r]^-\omega & \mathbf{Mod}-\E^f\ar[r]^-{\delta} &  \mathbf{Mod}-\E^f_{inj}
}$$ 
o\`u $\eta$ est induit par l'\'equivalence de cat\'egories du th\'eor\`eme~\ref{thm-fractions} et $\xi$ est l'\'equivalence donn\'ee dans le th\'eor\`eme~\ref{equ-span}, commute à isomorphisme canonique près.
\end{lm}

\begin{proof}
Soit $\alpha$ un objet de $\mathbf{Mod}-\E^f_{iso}$. Son image dans $\mathbf{Mod}-\E^f_{inj}$ en suivant le chemin supérieur est donnée par
$$V\mapsto\bigoplus_{q\in S^2(V^*)}\bigoplus_{W\subset Rad(V,q)}\alpha(W).$$
Sur les morphismes, une flèche $f : V\to V'$ de $\E^f_{inj}$ est transformée en l'application linéaire dont la composante $\alpha(W')\to\alpha(W)$, pour $q\in S^2(V^*)$, $q'\in S^2(V'^*)$, $W\subset Rad(V,q)$ et $W'\subset Rad(V',q')$, est $\alpha(\bar{f})$ lorsque $S^2(^t f)(q')=q$ et que $f$ induit un isomorphisme $\bar{f} : W\to W'$ (i.e. $f(W)=W'$) et est nulle sinon. 

On peut simplifier cette écriture en notant que la donnée d'une forme quadratique $q$ sur $V$ et d'un sous-espace $W$ de $V$ inclus dans son radical est équivalente à la donnée d'un sous-espace $W$ de $V$ et d'une forme quadratique $\bar{q}$ sur $V/W$. Notre foncteur devient alors
$$V\mapsto\bigoplus_{W\subset V}\alpha(W)\otimes \kk[S^2((V/W)^*)].$$
Sur les morphismes, la composante $\alpha(W')\otimes \kk[S^2((V'/W')^*)]\to\alpha(W)\otimes \kk[S^2((V/W)^*)]$ induite par $f$ est nulle si $f(W)\neq W'$ et est sinon égale au produit tensoriel de $\alpha(\bar{f})$ et du morphisme $\kk[S^2((V'/W')^*)]\to \kk[S^2((V/W)^*)]$ associant $[\overline{S^2(^t f)(q')}]$ à $[\bar{q'}]$ pour $q'$ forme quadratique sur $V'$ nulle sur $W'$. Ce dernier morphisme n'étant autre que $\kappa(\kk[S^2]^\vee)(f)$, cela établit le lemme.
\end{proof}

\begin{proof}[D\'emonstration du th\'eor\`eme~\ref{thf-o}]
 On se contente d'établir l'assertion relative au groupe orthogonal, l'autre étant tout-à-fait analogue.
 
 On utilise la notation $L_q$ de~\ref{not-fdl}, le triplet sous-jacent étant $(\Eqd,S_\mathbf{H},G_\mathbf{H})$ (cf. exemple~\ref{exfond}.\,\ref{exquad}).
 
 Compte-tenu de la proposition~\ref{torcomod}~, la suite spectrale de la proposition~\ref{fctder} prend la forme :
 $$E^2_{p,q}={\rm Tor}_p^{\Eqd}(L_q,U^*F)\simeq {\rm Tor}_p^{\E^f_{inj}}(\Omega_{(S^2)^\vee}(L_q),F)\Rightarrow H_{p+q}(O_\infty ; F_\infty)$$
 o\`u $U: \Eqd=(\E^f_{inj})_{(S^2)^{\vee}} \to \E^f_{inj}$ est le foncteur d'oubli (pour alléger, on a noté encore $F$ la restriction de ce foncteur à $\E^f_{inj}$).
 
  La proposition~\ref{fdo} montre que les foncteurs $L_q$ inversent les inclusions $V\hookrightarrow V\perp H$, où l'espace quadratique $H$ est non dégénéré ; ils appartiennent par conséquent à l'image essentielle du foncteur $\gamma$. Il existe donc, par le lemme précédent, des objets $\alpha_q$ de $\mathbf{Mod}-\E^f_{iso}$ tels que $\Omega_{(S^2)^\vee}(L_q)\simeq\delta\omega(\rho\beta(\alpha_q)\otimes\kappa(\kk[S^2]^\vee))$. Autrement dit, $\Omega_{(S^2)^\vee}(L_q)$ est isomorphe à la restriction à $\E^f_{inj}$ du foncteur $\omega(\rho\beta(\alpha_q)\otimes\kappa(\kk[S^2]^\vee))$ de $\mathbf{Mod}-\E^f$.
 
 Comme $F$ est par hypothèse analytique, on déduit du théorème~\ref{th-susl} (dû à Suslin) que le morphisme canonique ${\rm Tor}_p^{\E^f_{inj}}(\Omega_{(S^2)^\vee}(L_q),F)\to {\rm Tor}_p^{\E^f}(\omega(\rho\beta(\alpha_q)\otimes\kappa(\kk[S^2]^\vee)),F)$ est un isomorphisme.
 
Le caractère analytique de $F$ implique également que le monomorphisme canonique $\lambda(\rho\beta(\alpha_q)\otimes\kappa(\kk[S^2]^\vee))\to\omega(\rho\beta(\alpha_q)\otimes\kappa(\kk[S^2]^\vee))$ (le foncteur $\lambda$ est défini en~\ref{def-cataux}), induit un isomorphisme ${\rm Tor}_p^{\E^f}(\lambda(\rho\beta(\alpha_q)\otimes\kappa(\kk[S^2]^\vee)),F)\to {\rm Tor}_p^{\E^f}(\omega(\rho\beta(\alpha_q)\otimes\kappa(\kk[S^2]^\vee)),F)$, grâce au théorème~\ref{th-dja}.

On note à présent que $\lambda(\rho\beta(\alpha_q)\otimes\kappa(\kk[S^2]^\vee))$ est isomorphe à $\alpha_q(0)\otimes\kk[S^2]^\vee$. Quant au $\kk$-module $\alpha_q(0)$, c'est la valeur en $0$ du foncteur $L_q$, c'est-à-dire $H_q(O_\infty(k) ; \kk)$.
 
Par la proposition~\ref{isoc1}, on en déduit que le morphisme naturel $H_*(O_\infty(k);F_\infty)\to H_*(\Eqd(k)\times O_\infty(k);U^*F)$ est un isomorphisme. La proposition~\ref{torcomod} fournit par ailleurs un isomorphisme entre $H_*(\Eqd(k)\times O_\infty(k);U^*F)$ et ${\rm Tor}^{\E^f_{inj}\times O_\infty(k)}_*(\kk[S^2]^\vee,F)$, lui-même isomorphe à ${\rm Tor}^{\E^f\times O_\infty(k)}_*(\kk[S^2]^\vee,F)$ grâce au théorème de Suslin (\ref{th-susl}) déjà invoqué (comparer les suites spectrales de Künneth associées à ces deux derniers groupes de torsion). Cela termine la démonstration.
\end{proof}

\begin{proof}[Démonstration des corollaires~\ref{thfvdec} et~\ref{thfssv2}] Ils se déduisent du théorème~\ref{thf-o} et des propositions~\ref{arret2} et~\ref{eff-kun} respectivement, en utilisant comme précédemment la proposition~\ref{torcomod} et le théorème~\ref{th-susl} pour l'identification des deuxièmes pages.
\end{proof}

\begin{rem}
On vient de montrer que le morphisme naturel ${\rm Tor}_p^{\Eq}(M,F)\to {\rm Tor}_p^{\Eqd}(M,F)$, où $M$ est un foncteur constant et $F$ un foncteur analytique de $\F(k)$ (par abus on a omis la notation de précomposition par différents foncteurs d'oubli), est un isomorphisme.
\end{rem}

\begin{rem}
L'hypothèse de finitude du corps $k$ ne sert essentiellement que dans le lemme~\ref{dcm}, pour pouvoir affirmer que le foncteur $\xi$ du théorème~\ref{equ-span} est une équivalence de catégories (l'application dudit théorème à la catégorie $\E^f_{inj}$ suppose en effet que l'ensemble des sous-espaces d'un $k$-espace vectoriel de dimension finie soit fini).

Nous pensons cependant que le théorème~\ref{thf-o} reste valable lorsque $k$ est un corps commutatif infini. 
\end{rem}

\begin{nota}\label{not-inji}
Dans la catégorie $\F(k)$, on désigne par $I$ le foncteur $(\kk[-]^\vee)^* : V\mapsto\kk^{V^*}$.
\end{nota}
(C'est un objet injectif de $\F$.)

En utilisant la trivialité de l'homologie stable des groupes classiques sur un corps fini à coefficients dans le même corps (voir \cite{FPH}, chapitre~3, §\,4), la caractéristique $2$ étant exclue pour les groupes orthogonaux, on déduit du théorème~\ref{thf-o} le résultat suivant. Dans le cas où $k$
est de caractéristique $2$, l'homologie $H_*(O_\infty(k);\FF)$ n'est pas triviale mais isomorphe à $H_*(\mathbb{Z}/2;\FF)$ (donc de dimension $1$ en chaque degré), car le sous-groupe d'indice $2$ de $O_\infty(k)$, noté $DO$ dans \cite{FPH} (défini au chapitre 2, §\,7, via l'invariant de Dickson), analogue en caractéristique $2$ du groupe spécial orthogonal, a une homologie triviale (cf. \cite{FPH}, chapitre 3, §\,4).

\begin{cor}\label{crf-o}
Supposons que $\kk=k$ est un corps fini de caractéristique $p$. Si $F$ est un foncteur analytique de $\F(k)$, il existe des isomorphismes naturels
$$H_*(O_\infty(k) ; F_\infty)\simeq {\rm Tor}^{\E^f}_*(k[S^2]^\vee,F)\quad\text{si }p\text{ est impair,}$$
$$H_*(O_\infty(k) ; F_\infty)\simeq {\rm Tor}^{\E^f}_*(k[S^2]^\vee,F)\otimes H_*(\mathbb{Z}/2;k)\quad\text{si }p=2$$
$$\text{et}\qquad H_*(Sp_\infty(k) ; F_\infty)\simeq {\rm Tor}^{\E^f}_*(k[\Lambda^2]^\vee,F).$$

Les duaux de ces espaces vectoriels s'identifient canoniquement à
$${\rm Ext}^*_\F(F,I\circ\Gamma^2)\quad\text{si }p\text{ est impair,}\qquad{\rm Ext}^*_\F(F,I\circ\Gamma^2)\otimes H^*(\mathbb{Z}/2;k)\quad\text{si }p=2$$
$$\text{et}\qquad\qquad {\rm Ext}^*_\F(F,I\circ\Lambda^2)\qquad\qquad\text{respectivement.}$$
\end{cor}

En utilisant les résultats de stabilité établis par Charney dans \cite{Charney}, on en déduit le corollaire suivant.

\begin{cor}\label{crf-os}
Sous les hypothèses précédentes, supposons de plus que $F$ est polynomial de degré au plus $d$. Alors pour tous entiers $i$ et $n$ tels que $n\geq 2i+d+6$, on a des isomorphismes naturels
$$H_i(O_{n,n}(k) ; F(k^{2n}))^*\simeq {\rm Ext}^i_\F(F,I\circ\Gamma^2)\quad\text{si }p\text{ est impair},$$
$$H_i(O_{n,n}(k) ; F(k^{2n}))^*\simeq\bigoplus_{j=0}^i {\rm Ext}^j_\F(F,I\circ\Gamma^2)\quad\text{si }p=2$$
$$\text{et}\qquad H_i(Sp_{2n}(k) ; F(k^{2n}))^*\simeq {\rm Ext}^i_\F(F,I\circ\Lambda^2).$$
\end{cor}

Par dualité entre homologie et cohomologie d'un groupe fini, on obtient la variante suivante :
\begin{cor}\label{corco} Sous les mêmes hypothèses, il existe des isomorphismes naturels
$$H^i(O_{n,n}(k) ; F(k^{2n}))\simeq {\rm Ext}^i_\F(k[S^2],F)\quad\text{si }p\text{ est impair},$$
$$H^i(O_{n,n}(k) ; F(k^{2n}))\simeq\bigoplus_{j=0}^i {\rm Ext}^i_\F(k[S^2],F)\quad\text{si }p=2$$
$$\text{et}\qquad H^i(Sp_{2n}(k) ; F(k^{2n}))\simeq {\rm Ext}^i_\F(k[\Lambda^2],F).$$
\end{cor}

%% file: appc.tex
Cet appendice a pour but de montrer que l'argument donn\'e par le second auteur  pour obtenir l'\'equivalence de cat\'egories du th\'eor\`eme $4.2$ de \cite{CV}: 
$$Sp(\Eqd(\mathbb{F}_2))-\mathbf{Mod}\simeq\prod_{V \in \mathcal{S}} \mathbb{F}_2[O(V)]-\textbf{Mod}$$
o\`u $Sp(\Eqd(\FF))$ est la cat\'egorie de Burnside associ\'ee \`a $\Eqd(\FF)$ dont on rappelle la d\'efinition en~\ref{span} et $\mathcal{S}$ est un ensemble de rep\'esentants des classes d'isom\'etrie des objets de $\Eqd(\FF)$, s'adapte facilement \`a d'autres situations pour donner des \'equivalences de cat\'egories non triviales. L'une d'entre elles intervient dans la preuve du th\'eor\`eme~\ref{thf-o}.

Cet argument repose sur la combinaison de l'inversion de Möbius et du th\'eor\`eme de Morita-Freyd dont on rappelle ici les \'enonc\'es.

La fonction de Möbius $\mu_X$ d'un ensemble fini partiellement ordonn\'e $(X,\leq)$ est d\'efinie par:
$$\mu_X(x,x)=1 \quad \mathrm{pour \ tout\ }  x \mathrm{\ dans\ } X$$
$$\sum_{x \leq z \leq y} \mu_X(x,z)=0 \quad \mathrm{pour \ tous\ } x <y \  \mathrm{dans}\ X.$$

\begin{thm}[Th\'eor\`eme 3.9.2 de \cite{Stanley}] \label{stanley}
Soit $(X, \leq)$ un ensemble fini partiellement ordonn\'e, dans lequel toute paire $\{x,y\}$ a une borne inférieure $x \wedge y$. Supposons que $X$ a un plus grand \'el\'ement $M$. Soit  $R$ un anneau (avec unit\'e $1_R$), et supposons que $\alpha \mapsto e_{\alpha}$ est une application de $X$ dans  $R$ v\'erifiant les propri\'et\'es suivantes: $e_{\alpha}e_{\beta}=e_{\alpha \wedge \beta}$ pour tout $\alpha, \beta \in X$, et $e_M=1_R$. Pour $\alpha \in X$, on d\'efinit:
$$f_{\alpha}=\sum_{\beta \leq \alpha}\mu_X(\beta, \alpha) e_{\beta},$$
o\`u $\mu_X$ est la fonction de  M\"{o}bius de $X$. Alors les \'el\'ements $f_{\alpha}$, pour $\alpha \in X$, sont des idempotents orthogonaux de $R$, dont la somme est \'egale \`a $1_R$.
\end{thm}

Le th\'eor\`eme suivant, dû à Freyd, doit \^etre vu comme une forme générale, "à plusieurs objets", du théorème classique de Morita sur l'équivalence des catégories de modules.

\begin{thm}[Morita-Freyd] \label{Morita-Freyd}
Soient $\A$ une catégorie de Grothendieck $\kk$-linéaire (i.e. enrichie sur les $\kk$-modules) et $\mathcal{G}$ une sous-catégorie pleine petite de $\A$ dont les objets sont projectifs de type fini et engendrent $\A$. Alors $\A$ est équivalente à la catégorie des foncteurs $\kk$-linéaires de $\mathcal{G}^{op}$ dans $\mathbf{Mod}_\kk$.
\end{thm}

On commence par appliquer ces outils \`a la cat\'egorie $\Gamma$ des ensembles finis point\'es et on montre que cela permet de retrouver le th\'eor\`eme \`a la  Dold-Kan d\'emontr\'e par Pirashvili dans l'article \cite{PDK}, qu'on utilisera dans l'appendice~\ref{apsym}. On \'etudie ensuite le cas, tr\`es proche de celui consid\'er\'e dans \cite{CV},  de la cat\'egorie de Burnside associ\'ee \`a une "bonne cat\'egorie" et qui fournit une \'equivalence de cat\'egories qu'on emploie dans le paragraphe~\ref{thm-fond}.

On désigne par  $\Omega$ la catégorie des ensembles finis avec surjections. On note $(-)_+$ l'adjoint à gauche au foncteur d'oubli de $\Gamma$ vers la catégorie des ensembles finis :  pour $E$ un ensemble fini,  $E_+$ s'obtient en adjoignant un point de base externe à $E$.

\begin{thm}[Pirashvili]\label{prdk}
Il existe une équivalence de catégories $cr : \mathbf{Mod}-\Gamma\to\mathbf{Mod}-\Omega$ telle que
$$cr(F)(E)=Coker\,\big(\bigoplus_{e\in E}F(E_e)\to F(E_+)\big)$$
morphisme induit par les surjections $E_+\twoheadrightarrow E_e$ égales à l'identité hors du point de base, où $E_e$ désigne l'ensemble $E$ pointé par $e$.

De plus, le foncteur $i_! : \Omega-\mathbf{Mod}\to\Gamma-\mathbf{Mod}$ tel que 
$$G\underset{\Gamma}{\otimes}i_!(F)\simeq cr(G)\underset{\Omega}{\otimes} F$$
est donné par
$$i_!(F)(E)=\bigoplus_{E'\subset E\setminus\{*\}}F(E').$$
\end{thm}

\begin{proof}
Soit $E$ un objet de $\Gamma$ de point de base $*$. On consid\`ere l'ensemble $\mathfrak{p}(E)$ de ses sous-objets (i.e. de ses sous-ensembles contenant le point de base) ordonn\'e par l'inclusion;  $\mathfrak{p}(E)$ admet pour plus grand \'el\'ement $E$ et toute paire d'\'el\'ement de $\mathfrak{p}(E)$, $\{A, B\}$ admet pour plus grand minorant l'intersection $A \cap B$. Soit $R=\kk[\mathrm{End}_{\Gamma}(E)]$ la $\kk$-algèbre du monoïde $\mathrm{End}_{\Gamma}(E)$. À un \'el\'ement $A$ de $\mathfrak{p}(E)$ on associe l'\'el\'ement $e_A \in \mathrm{End}_{\Gamma}(E)$  donné par $e_A(x)=x$ si $x\in A$ et $e_A(x)=*$ sinon. On a $e_E=Id_E$ et $e_Ae_B=e_{A\cap B}$. Par cons\'equent, par le th\'eor\`eme~\ref{stanley} les \'el\'ements $f_A$ d\'efinis par:
$$f_A=\sum_{B\subset A}\mu_{\mathfrak{p}(E)}(B,A)e_B$$
sont des idempotents orthogonaux de $R$ de somme \'egale \`a $Id_E$.
On en déduit la décomposition 
\begin{equation}
P^\Gamma_E\simeq \bigoplus_{A\in\mathfrak{p}(E)} f_A P^\Gamma_E,
\end{equation}
où les $P^\Gamma_E$ désignent les projectifs standard de $\Gamma-\mathbf{Mod}$. Les foncteurs $P^\Gamma_E$ sont de type fini et donc leurs facteurs directs $f_A P^\Gamma_E$ forment un ensemble de g\'en\'erateurs projectifs de type fini de $\Gamma-\mathbf{Mod}$.

Afin d'appliquer le th\'eor\`eme~\ref{Morita-Freyd}, nous avons besoin d'identifier les modules suivants
$${\rm Hom}\,(f_A P^\Gamma_E,f_B P^\Gamma_E)\simeq f_B \kk[{\rm End}_\Gamma (E)] f_A.$$

 Pour cela, on fait les observations suivantes, pour tous $A\in\mathfrak{p}(E)$ et $t\in {\rm End}_\Gamma (E)$ :
\begin{enumerate}
\item $e_A t=t e_{t^{-1}(A)}$ ;
\item $t=t e_B$, où $B$ désigne l'élément de $\mathfrak{p}(E)$ réunion du point de base et de l'ensemble complémentaire de $t^{-1}(*)$ dans $E$.
\end{enumerate}

Du premier point on déduit
$$f_A t=t\sum_{B\subset A}\mu_{\mathfrak{p}(E)}(B,A)e_{t^{-1}(B)}$$
puis, compte-tenu de ce que $e_C f_D=f_D$ si $D\subset C$, $0$ sinon,
$$f_A t f_{A'}=t\sum_{t(A')\subset B\subset A}\mu_{\mathfrak{p}(E)}(B,A) f_{A'}=\left\lbrace\begin{array}{l}
t f_{A'}\text{ si }t(A')=A, \\
0\text{ sinon.}
\end{array}
\right.$$

Du second vient que $t f_{A'}$ est nul sauf si $A'\subset (E\setminus t^{-1}(*))\cup\{*\}$. Cette condition signifie que $A'$ induit une fonction de $A'\setminus\{*\}$ dans $E\setminus\{*\}$, ou une surjection de $A'\setminus\{*\}$ vers $A\setminus\{*\}$ si $t(A')=A$.

Des deux observations précédentes, et du fait que $t f_{A'}$, comme les $t e_B$ pour $B\subset A'$, ne dépend que de la restriction de $t$ à $A'$, on déduit une application linéaire surjective
\begin{equation}\label{grabi}
\kk[{\rm Surj}\,(A'\setminus\{*\},A\setminus\{*\})]\twoheadrightarrow f_A\kk[{\rm End}_\Gamma (E)] f_{A'},
\end{equation}
où Surj désigne l'ensemble des fonctions surjectives entre deux ensembles.

Un argument de rang montre que cette application est en fait bijective : en effet, la somme directe des $\kk$-modules libres  $f_A\kk[{\rm End}_\Gamma (E)] f_{A'}$ lorsque $A$ et $A'$ parcourent $\mathfrak{p}(E)$ est isomorphe à $\kk[{\rm End}_\Gamma (E)]$. On conclut par la bijection
$${\rm End}_\Gamma (E)\simeq\bigsqcup_{(A,A')\in\mathfrak{p}(E)^2}{\rm Surj}\,(A'\setminus\{*\},A\setminus\{*\})$$ 
obtenue en associant à un endomorphisme $t$ de $E$ les éléments $A=t(E)$ et $A'=(E\setminus t^{-1}(*))\cup\{*\}$ de $\mathfrak{p}(E)$ et la surjection $A'\setminus\{*\}\twoheadrightarrow A\setminus\{*\}$ induite par $t$.

On constate d'autre part que l'isomorphisme~(\ref{grabi}) est compatible à la composition en un sens évident. Il nous reste \`a expliciter les \'equivalences de cat\'egories donn\'ees par le théorème de Freyd-Morita. La d\'ecomposition (\ref{grabi}) fournit un isomorphisme
$$(i_!(X))(E)\simeq\bigoplus_{A\in\mathfrak{p}(E)}X(A\setminus\{*\})=\underset{*\not\in E'}{\bigoplus_{E'\subset E}}X(E')$$
pour tout objet $X$ de $\Omega-\mathbf{Mod}$ et tout objet $E$ de $\Gamma$. On vérifie facilement que la fonctorialité s'obtient en écrivant $i_!(X)(E)$ comme conoyau de l'injection naturelle
$$\bigoplus_{A\in\mathfrak{p}(E)}X(A)\hookrightarrow\bigoplus_{A\subset E}X(A).$$
(Autrement dit, la composante $X(E')\to X(F')$ de $X(f)$, pour $f : E\to F$ morphisme de $\Gamma$, est induite par $f$ si $f(E')=F'$ et nulle sinon.)

La formule pour $cr$ s'obtient de mani\`ere analogue.
\end{proof}

Avant d'appliquer cette m\'ethode au cas dont on se sert au paragraphe~\ref{thm-fond}, on rappelle la d\'efinition de la {\em cat\'egorie de Burnside} d'une cat\'egorie admettant des produits fibr\'es. Cette terminologie est issue de la théorie des représentations; la notation $Sp$ provient du terme anglais \textit{span}. 

\begin{defi}[\cite{Benabou}]  \label{span}
Soit $\C$ une cat\'egorie admettant des produits fibr\'es, la cat\'egorie de Burnside de $\C$, not\'ee $Sp(\C)$, est d\'efinie de la mani\`ere suivante:
\begin{enumerate}
\item les objets de $Sp(\C)$ sont ceux de $\C$;
\item pour $A$ et $B$ deux objets de $Sp(\C)$, ${\rm Hom}_{Sp(\C)}(A,B)$ est l'ensemble des classes d'\'equivalence de diagrammes dans $\C$ de la forme $A \xleftarrow{f} D \xrightarrow{g} B$, pour la relation d'\'equivalence qui identifie les deux diagrammes $A \xleftarrow{f} D \xrightarrow{g} B$ et $A \xleftarrow{u} D' \xrightarrow{v} B$ s'il existe un isomorphisme $\alpha: D \to D'$ rendant le diagramme suivant commutatif:
$$\xymatrix{
D \ar[rrd]^-g \ar[rd]^(.65){\alpha} \ar[ddr]_-f\\
 &D' \ar[r]_v \ar[d]^u & B\\
 & A.
}$$
Le morphisme de ${\rm Hom}_{Sp(\C)}(A,B)$  repr\'esent\'e par le diagramme $A \xleftarrow{f} D \xrightarrow{g} B$ sera not\'e $[A \xleftarrow{f} D \xrightarrow{g} B]$.
\item Pour deux morphismes $T_1=[A \leftarrow D \rightarrow B]$ et $T_2=[B \leftarrow D' \rightarrow C]$ la composition est donn\'ee par: 
$$T_2 \circ T_1=[A \leftarrow D \underset{B}{\times} D' \rightarrow C ].$$
\end{enumerate}

On appelle foncteur de Mackey non additif depuis $\C$ un foncteur $Sp(\C)\to\mathbf{Mod}_\kk$.
\end{defi}

Donnons quelques notations et définitions supplémentaires dans le cas où la catégorie $\C$ est petite et où tous ses morphismes sont des monomorphismes. Un {\em sous-objet} d'un objet $i$ de $\C$ est une classe d'équivalence de morphisme de but $i$ pour la relation identifiant $f : j\to i$ à $f' : j'\to i$ lorsqu'il existe un isomorphisme $g : j\xrightarrow{\simeq} j'$ tel que $f=f'g$. On notera $\mathfrak{p}(i)$ l'ensemble des sous-objets de $i$ ; on supposera de plus effectué le choix d'un représentant dans chaque classe (ces représentants seront spécifiés par la notation $\hookrightarrow$), auquel on identifiera celle-ci. L'ensemble $\mathfrak{p}(i)$ est muni d'une relation d'ordre notée $\subset$ définie par $k\subset j$ lorsque le morphisme $k\hookrightarrow i$ se factorise (de manière unique puisque les flèches de $\C$ sont des monomorphismes) par $j\hookrightarrow i$ ; l'existence de produits fibrés dans $\C$ assure que deux éléments $j$ et $k$ de cet ensemble ordonné possèdent une borne inférieure notée $j\cap k$.
 
L'image d'un morphisme de $\C$ est par définition sa classe dans l'ensemble des sous-objets de son but. Si $t=[i\xleftarrow{u} k\xrightarrow{v}j]$ est un élément de ${\rm Hom}_{Sp(\C)}(i,j)$, l'image de $u$ (resp. $v$), qui ne dépend pas du choix des représentants, est notée $coim(t)$ (resp. $im(t)$) ; $t$ possède une unique écriture $t=[i\hookleftarrow coim(t)\xrightarrow{\bar{t}}j]$.

Pour $f\in {\rm Hom}_\C(i,i')$, $j\in\mathfrak{p}(i)$ et $j'\in\mathfrak{p}(i')$, on note $f(j)=im(j\hookrightarrow i\xrightarrow{f}i')$ et $f^{-1}(j')$ l'élément de $\mathfrak{p}(i)$ donné par le produit fibré de $f$ et de $j'\hookrightarrow j$.

\begin{thm}[Cf. \cite{CV}, théorème~$4.2$] \label{equ-span}
Soit $\C$ une petite cat\'egorie admettant des produits fibr\'es, dont tous les morphismes sont des monomorphismes et telle que l'ensemble $\mathfrak{p}(i)$ est fini pour tout $i\in {\rm Ob}\,\C$. Notons $\C_{iso}$ la sous-cat\'egorie de $\C$ ayant les m\^emes objets et pour morphismes les isomorphismes de $\C$. Il existe une \'equivalence de cat\'egories
$$\xi : \C_{iso}-\mathbf{Mod}\xrightarrow{\simeq} Sp(\C)-\mathbf{Mod}$$
telle que $\xi(M)(i)=\underset{j\in\mathfrak{p}(i)}{\bigoplus} M(j)$ et qu'un morphisme $t\in {\rm Hom}_{Sp(\C)}(i,i')$ induit l'application $\xi(M)(i)\to\xi(M)(i')$ dont la composante $M(j)\to M(j')$ (pour $j\in\mathfrak{p}(i)$ et $j'\in\mathfrak{p}(i')$) est $M(\tilde{t})$ si $j\subset coim(t)$ et que $\bar{t}(j)=j'$, de sorte que $\bar{t}$ induit un isomorphisme $\tilde{t} : j\xrightarrow{\simeq}j'$, et est nulle sinon.
\end{thm}

\begin{proof} Reprenant les arguments de \cite{CV}, on procède de manière très analogue à la démonstration du théorème~\ref{prdk}. 

Soient $c$ un objet de $\C$ et $i$ un élément de $\mathfrak{p}(c)$ ; on note $e^c_i$ (ou simplement $e_i$) l'idempotent $[c\hookleftarrow i\hookrightarrow c]$ de ${\rm End}_{Sp(\C)}(c)$. On a $e^c_c=Id_c$ et $e_i^c.e_j^c=e^c_{i\cap j}$, ce qui permet d'appliquer le théorème~\ref{stanley} pour obtenir une famille complète d'idempotents orthogonaux
$$f_a^c=\sum_{i\subset a}\mu_{\mathfrak{p}(c)}(i,a)e_i^c$$
de l'anneau $\kk[{\rm End}_{Sp(\C)}(c)]$.

On remarque maintenant que pour tout morphisme $t\in {\rm Hom}_{Sp(\C)}(c,d)$, on a :
\begin{enumerate}
\item $t e^c_i=e^d_{\bar{t}(i)}t$ pour $i\subset coim(t)$ ;
\item $t e^c_{coim(t)}=t$.
\end{enumerate}
Le second point montre que $t f^c_a=t e^c_{coim(t)} f^c_a=0$ sauf si $a\subset coim(t)$.

On déduit par ailleurs du premier point que, sous l'hypothèse  $a\subset coim(t)$ :
$$f^d_b t f^c_a=\sum_{i\subset a}\mu_{\mathfrak{p}(c)}(i,a)f^d_b e^d_{\bar{t}(i)} t=\underset{b\subset\bar{t}(i)}{\sum_{i\subset a}}\mu_{\mathfrak{p}(c)}(i,a) f^d_b t=\left\lbrace\begin{array}{ll}
f^d_b t\text{ si }\bar{t}(a)=b, \\
0\text{ sinon.}
\end{array}
\right.$$

Par conséquent, $f^d_b t f^c_a$ est nul sauf si $a\subset coim(t)$ et que $\bar{t}$ induit un isomorphisme $\tilde{t}$ de $a$ vers $b$. De surcroît, comme $f^d_b=f^d_b e^d_b$ et $f^c_a=e^c_a f^c_a$, $f^d_b t f^c_a$ ne dépend alors que de $e^d_b t e^c_a$, c'est-à-dire de $\tilde{t}$. Il s'ensuit que l'application linéaire
$$\kk[{\rm Hom}_{\C_{iso}}(a,b)]\to f^d_b\kk[{\rm Hom}_{Sp(\C)}(c,d)]f^c_a\qquad [u]\mapsto f^d_b[c\hookleftarrow a\xrightarrow{u}b\hookrightarrow d]f^c_a$$
est surjective. On montre que c'est en fait un isomorphisme par un argument de rang analogue à celui utilisé pour le théorème~\ref{prdk}, et l'on termine aussi la démonstration de la même façon.
\end{proof}

%% file: apcof.tex
Cet appendice a pour objet de rappeler brièvement deux résultats d'annulation homologique utilisés de manière cruciale dans la démonstration du théorème~\ref{thf-o}. On s'y donne un corps {\em fini} $k$, qui sera souvent sous-entendu dans les notations.

\begin{defi}\label{def-cataux}
 \begin{itemize}
  \item On désigne par  $\E^f_{surj}$ (resp. $\E^f_{inj}$) la sous-catégorie de la catégorie $\E^f$ des $k$-espaces vectoriels de dimension finie ayant les mêmes objets et les surjections (resp. les injections) pour morphismes.
\item On désigne par $\E^f_\Gr$ est la catégorie des couples $(V,W)$ formés d'un espace vectoriel de dimension finie $V$ et d'un sous-espace vectoriel $W$, avec pour morphismes $(V,W)\to (V',W')$ les applications linéaires $f : V\to V'$ telles que $f(W)=W'$.
\item On note $\iota : \mathbf{Mod}-\E^f\to \mathbf{Mod}-\E^f_\Gr$ le foncteur de précomposition par $\E^f_\Gr\to\E^f\quad (V,W)\mapsto V$.
\item On note $\rho : \mathbf{Mod}-\E^f_{surj}\to \mathbf{Mod}-\E^f_\Gr$ la précomposition par $\E^f_\Gr\to\E^f_{surj}\quad (V,W)\mapsto W$.
 \item On note  $\kappa : \mathbf{Mod}-\E^f\to\mathbf{Mod}-\E^f_\Gr$ la précomposition par $\E^f_\Gr\to\E^f\quad (V,W)\mapsto V/W$. 
 \item On note $\lambda : \mathbf{Mod}-\E^f_\Gr\to\mathbf{Mod}-\E^f$ la précomposition par le foncteur $\E^f\to\E^f_\Gr\quad V\mapsto (V,0)$
 \item On désigne par $\omega :  \mathbf{Mod}-\E^f_\Gr\to\mathbf{Mod}-\E^f$ le foncteur défini par $\omega(X)(V)=\bigoplus_{W\subset V} X(V,W)$ et le fait que $\omega(X)(f)$, pour $f : V\to V'$ morphisme de $\E^f$, a pour composante $X(f) : X(V',W')\to X(V,W)$ si $f(W)=W'$ et $0$ sinon (cf. proposition~\ref{torcomod}).
 \end{itemize}
\end{defi}

On remarque que l'on dispose d'un isomorphisme canonique $\omega(X\otimes\iota(F))\simeq\omega(X)\otimes F$, où $X\in {\rm Ob}\,\mathbf{Mod}-\E^f_\Gr$ et $F\in {\rm Ob}\,\mathbf{Mod}-\E^f$. Par ailleurs, l'inclusion évidente de $X(V,0)$ dans $\bigoplus_{W\subset V} X(V,W)$ définit une transformation naturelle $\lambda\hookrightarrow\omega$.

\smallskip

Le premier résultat d'homologie des foncteurs (sa démonstration établit clairement qu'il s'agit d'un résultat d'annulation) dont nous avons besoin pour établir le théorème principal de cet article est le suivant. Il s'agit de la variante en termes de groupes de torsion d'un cas particulier fondamental du théorème~10.2.1 de \cite{Dja} (cf. aussi son corollaire~10.2.2).

\begin{thm}[Djament]\label{th-dja}
  Soient $X\in\mathbf{Mod}-\E^f_\Gr$ et $F$ un foncteur analytique de $\F$ (cf. définition~\ref{poly}). L'inclusion naturelle $\lambda(X)\hookrightarrow\omega(X)$ induit un isomorphisme
  $${\rm Tor}^{\E^f}_*(\lambda(X),F)\simeq {\rm Tor}^{\E^f}_*(\omega(X),F).$$
\end{thm}

Il revient au même de démontrer l'annulation de ${\rm Tor}^{\E^f}_*(\omega(X)/\lambda(X),F)$. Quitte à remplacer $X$ par son quotient $X'$ défini par $X'(V,W)=X(V,W)$ si $W\neq 0$ et $X'(V,0)=0$, on est donc ramené à établir l'annulation de ${\rm Tor}^{\E^f}_*(\omega(X),F)$ lorsque $\lambda(X)=0$. On s'appuie pour cela sur les lemmes suivants.

On commence par noter $\bar{P}\in\mathbf{Mod}-\E^f$ le noyau du morphisme d'augmentation $P^{(\E^f)^{op}}_k\twoheadrightarrow P^{(\E^f)^{op}}_0=\kk$. Ainsi, les éléments $[l]- [0]$, où $l$ parcourt les formes linéaires non nulles sur $V$, forment une base de $\bar{P}(V)$. On a classiquement :
\begin{lm}\label{lmad1}
 Il existe un isomorphisme naturel ${\rm Tor}^{\E^f}_*(F\otimes\bar{P},G)\simeq {\rm Tor}^{\E^f}_*(F, \Delta(G))$.
\end{lm}
(On rappelle que le foncteur différence $\Delta$ est défini en~\ref{poly}.)

Ce lemme se déduit aussitôt du cas du degré $0$, par exactitude des foncteurs en jeu ; la propriété provient alors de l'isomorphisme canonique $P^{(\E^f)^{op}}_V\otimes P^{(\E^f)^{op}}_W\simeq P^{(\E^f)^{op}}_{V\oplus W}$ en considérant une présentation projective de $F$.

\begin{lm}\label{lmad2}
 Si $X\in\mathbf{Mod}-\E^f_\Gr$ vérifie $\lambda(X)=0$, il existe une suite exacte courte $0\to Y\to X\otimes\iota(\bar{P})\to X\to 0$ où $Y\in\mathbf{Mod}-\E^f_\Gr$ vérifie $\lambda(Y)=0$.
\end{lm}

\begin{proof}
 On définit un morphisme $\iota(\bar{P})\to\kk$ en associant à $[l]-[0]\in \iota(\bar{P})(V,W)=Ker\,(\kk[V^*]\to\kk)$ $0$ si la restriction de $l$ à $W$ est nulle et $1$ dans le cas contraire. Ce morphisme est surjectif si $W$ est non nul. En le tensorisant par $X$, on obtient le résultat souhaité. 
\end{proof}

\begin{proof}[Démonstration du théorème~\ref{th-dja}] Par un argument de colimite, on peut supposer $F$ polynomial : il existe $d$ tel que $\Delta^d(F)=0$. On a vu qu'on peut également supposer $\lambda(X)=0$. En itérant $d$ fois l'épimorphisme du lemme~\ref{lmad2}, on obtient une suite exacte $0\to Y\to X\otimes\iota(\bar{P})^{\otimes d}\to X\to 0$ où $Y\in\mathbf{Mod}-\E^f_\Gr$ vérifie $\lambda(Y)=0$. Le foncteur $\omega$ étant exact, on en déduit une suite exacte $0\to\omega(Y)\to\omega(X)\otimes\bar{P}^{\otimes d}\to\omega(X)\to 0$. La suite exacte longue d'homologie associée fournit, compte-tenu de ce que ${\rm Tor}^{\E^f}_*(\omega(X)\otimes\bar{P}^{\otimes d},F)\simeq {\rm Tor}^{\E^f}_*(\omega(X),\Delta^d(F))=0$ (par le lemme~\ref{lmad1} appliqué $d$ fois), l'annulation de ${\rm Tor}^{\E^f}_0(\omega(X),F)$ et des isomorphismes ${\rm Tor}^{\E^f}_i(\omega(X),F)\simeq {\rm Tor}^{\E^f}_{i-1}(\omega(Y),F)$. Par récurrence sur le degré homologique, le théorème s'ensuit. 
\end{proof}

Le deuxième résultat d'homologie des foncteurs dont nous avons besoin dans cet article est le suivant. C'est une variante en termes de groupes de torsion du théorème A.8 de l'appendice de \cite{FFSS} ; nous esquissons une démonstration se fondant sur le théorème~\ref{th-dja}, en suivant \cite{Dja}, §\,13.2.

\begin{thm}[Suslin]\label{th-susl} Soient $A\in {\rm Ob}\,\mathbf{Mod}-\E^f$ et $F\in {\rm Ob}\,\F$ analytique. Le morphisme canonique
 $${\rm Tor}_*^{\E^f_{inj}}(\delta(A),\delta(F))\xrightarrow{\simeq} {\rm Tor}_*^{\E^f}(A,F)$$
 est un isomorphisme, où $\delta$ désigne le foncteur de restriction à la sous-catégorie $\E^f_{inj}$ des injections de $\E^f$.
\end{thm}

\begin{lm}\label{adel2}
 Soient $X\in {\rm Ob}\,\mathbf{Mod}-\E^f_{inj}$ et $F\in {\rm Ob}\,\F$. Il existe un isomorphisme naturel ${\rm Tor}_*^{\E^f_{inj}}(X,\delta(F))\simeq {\rm Tor}_*^{\E^f}(\varpi(X),F)$, où $\varpi(X)$ est défini par $\varpi(X)(V)=\bigoplus_{W\subset V}X(V/W)$ et le fait que, pour toute application linéaire $f : V\to V'$, $\varpi(X)(f) : \varpi(X)(V')\to\varpi(X)(V)$ a pour composante $X(V'/W')\to X(V/W)$ le morphisme induit par le monomorphisme $V/W\hookrightarrow V'/W'$ induit par $f$ si $f^{-1}(W')=W$, $0$ sinon. 
\end{lm}

\begin{proof}
 L'isomorphisme en degré $0$ se déduit facilement de la décomposition
 $${\rm Hom}_{\E^f}(V,-)\simeq\coprod_{W\subset V} {\rm Hom}_{\E^f_{inj}}(V/W,-)$$
 (qui traite le cas $F=P^{\E^f}_V$).
 Le cas général s'en déduit par exactitude de $\varpi$ et $\delta$.
\end{proof}

\begin{lm}\label{laux2} Les endofoncteurs $\omega\kappa$ et $\varpi\delta$ de $\mathbf{Mod}-\E^f$ sont isomorphes.
\end{lm}

\begin{proof} On a $\omega\kappa(A)(V)=\varpi\delta(A)(V)=\bigoplus_{W\subset V}A(V/W)$, mais l'effet sur les morphismes n'est pas le même (l'un utilise l'image directe, l'autre l'image inverse d'un sous-espace vectoriel par une application linéaire). Néanmoins, on vérifie facilement (cf. \cite{Dja}, proposition~13.2.1 pour les détails) que l'application linéaire $\omega\kappa(A)(V)\to\varpi\delta(A)(V)$ ayant pour composante $A(V/W)\to A(V/W')$ le morphisme induit par la projection $V/W\twoheadrightarrow V/W'$ lorsque $W\subset W'$ et $0$ sinon est un isomorphisme et définit une transformation naturelle $\omega\kappa\to\varpi\delta$.
\end{proof}

\begin{proof}[Démonstration du théorème~\ref{th-susl}] En utilisant successivement les lemmes~\ref{adel2} et~\ref{laux2} puis le théorème~\ref{th-dja}, on obtient des isomorphismes naturels
$${\rm Tor}_*^{\E^f_{inj}}(\delta(A),\delta(F))\simeq {\rm Tor}_*^{\E^f}(\varpi\delta(A),F)\simeq {\rm Tor}_*^{\E^f}(\omega\kappa(A),F)\simeq {\rm Tor}_*^{\E^f}(\lambda\kappa(A),F);$$
 on conclut en remarquant que $\lambda\kappa\simeq id$.
\end{proof}

En fait, le théorème~\ref{th-susl} est valable dans un cadre beaucoup plus général :
\begin{thm}[Scorichenko]\label{thscogl}
 Soient $\A$ une catégorie additive (essentiellement) petite et $\A_{inj}$ la sous-catégorie des monomorphismes scindés de $\A$. Soient $F\in {\rm Ob}\,\mathbf{Mod}-\A$ et $G\in {\rm Ob}\,\A-\mathbf{Mod}$, $G$ étant supposé analytique\,\footnote{La définition est analogue à celle donnée dans le cadre de la catégorie $\F$. Le seul point à noter dans la définition de foncteur polynomial est qu'il faut imposer la nilpotence pour {\em tous} les foncteurs différences imaginables (avec un degré commun).}. Alors le morphisme canonique
 $${\rm Tor}^{\A_{inj}}_*(F,G)\to {\rm Tor}^{\A}_*(F,G)$$
 est un isomorphisme (où l'on a, par abus, omis les foncteurs de restriction à $\A_{inj}$ dans le membre de gauche).
\end{thm}

La principale difficulté réside dans le fait que l'analogue du foncteur $\varpi$ qui apparaît dans le contexte général n'est ni explicite ni exact. Pour la démonstration, nous renvoyons à \cite{Sco} (non publié), ou aux notes \cite{Dja-notes} qui en rendent disponibles les arguments.

Dans l'appendice~\ref{appBS}, nous utiliserons le cas particulier de ce théorème dans lequel $\A$ est la catégorie des modules projectifs de type fini sur un anneau fixé et $F$ un foncteur constant pour appliquer notre formalisme général à l'homologie des groupes linéaires.

%% file: dv.bbl
\providecommand{\bysame}{\leavevmode ---\ }
\providecommand{\og}{``}
\providecommand{\fg}{''}
\providecommand{\smfandname}{\&}
\providecommand{\smfedsname}{\'eds.}
\providecommand{\smfedname}{\'ed.}
\providecommand{\smfmastersthesisname}{M\'emoire}
\providecommand{\smfphdthesisname}{Th\`ese}
\begin{thebibliography}{TvdK08}

\bibitem[B{\'e}n67]{Benabou}
{\scshape J.~B{\'e}nabou} -- {\og Introduction to bicategories\fg}, in
  \emph{Reports of the Midwest Category Seminar}, Springer, Berlin, 1967,
  p.~1--77.

\bibitem[Bet89]{Bet1}
{\scshape S.~Betley} -- {\og Vanishing theorems for homology of {${\rm Gl}\sb
  nR$}\fg}, \emph{J. Pure Appl. Algebra} \textbf{58} (1989), no.~3,
  p.~213--226.

\bibitem[Bet92]{Bet2}
\bysame , {\og Homology of {${\rm Gl}(R)$} with coefficients in a functor of
  finite degree\fg}, \emph{J. Algebra} \textbf{150} (1992), no.~1, p.~73--86.

\bibitem[Bet99]{Bet}
\bysame , {\og Stable {$K$}-theory of finite fields\fg}, \emph{$K$-Theory}
  \textbf{17} (1999), no.~2, p.~103--111.

\bibitem[Bet02]{Bet-sym}
\bysame , {\og Twisted homology of symmetric groups\fg}, \emph{Proc. Amer.
  Math. Soc.} \textbf{130} (2002), no.~12, p.~3439--3445 (electronic).

\bibitem[BP94]{BP}
{\scshape S.~Betley {\normalfont \smfandname} T.~Pirashvili} -- {\og Stable
  {$K$}-theory as a derived functor\fg}, \emph{J. Pure Appl. Algebra}
  \textbf{96} (1994), no.~3, p.~245--258.

\bibitem[Cha87]{Charney}
{\scshape R.~Charney} -- {\og A generalization of a theorem of {V}ogtmann\fg},
  in \emph{Proceedings of the Northwestern conference on cohomology of groups
  (Evanston, Ill., 1985)}, vol.~44, 1987, p.~107--125.

\bibitem[Cha08]{Chal}
{\scshape M.~Cha{\l}upnik} -- {\og Koszul duality and extensions of exponential
  functors\fg}, \emph{Adv. Math.} \textbf{218} (2008), no.~3, p.~969--982.

\bibitem[Dja07]{Dja}
{\scshape A.~Djament} -- {\og Foncteurs en grassmanniennes, filtration de
  {K}rull et cohomologie des foncteurs\fg}, \emph{M\'em. Soc. Math. Fr. (N.S.)}
  (2007), no.~111, p.~xx+213.

\bibitem[Dja09]{Dja-notes}
\bysame , {\og Les r\'esultats d'annulation homologique de {S}corichenko\fg},
  disponible sur http://hal.archives-ouvertes.fr/hal-00411929/, 2009.

\bibitem[Dol60]{Dold}
{\scshape A.~Dold} -- {\og Zur {H}omotopietheorie der {K}ettenkomplexe\fg},
  \emph{Math. Ann.} \textbf{140} (1960), p.~278--298.

\bibitem[FFPS03]{FFPS}
{\scshape V.~Franjou, E.~M. Friedlander, T.~Pirashvili {\normalfont
  \smfandname} L.~Schwartz} -- \emph{Rational representations, the {S}teenrod
  algebra and functor homology}, Panoramas et Synth\`eses [Panoramas and
  Syntheses], vol.~16, Soci\'et\'e Math\'ematique de France, Paris, 2003.

\bibitem[FFSS99]{FFSS}
{\scshape V.~Franjou, E.~M. Friedlander, A.~Scorichenko {\normalfont
  \smfandname} A.~Suslin} -- {\og General linear and functor cohomology over
  finite fields\fg}, \emph{Ann. of Math. (2)} \textbf{150} (1999), no.~2,
  p.~663--728.

\bibitem[FLS94]{FLS}
{\scshape V.~Franjou, J.~Lannes {\normalfont \smfandname} L.~Schwartz} -- {\og
  Autour de la cohomologie de {M}ac {L}ane des corps finis\fg}, \emph{Invent.
  Math.} \textbf{115} (1994), no.~3, p.~513--538.

\bibitem[FP78]{FPH}
{\scshape Z.~Fiedorowicz {\normalfont \smfandname} S.~Priddy} -- \emph{Homology
  of classical groups over finite fields and their associated infinite loop
  spaces}, Lecture Notes in Mathematics, vol. 674, Springer, Berlin, 1978.

\bibitem[FS97]{FS}
{\scshape E.~M. Friedlander {\normalfont \smfandname} A.~Suslin} -- {\og
  Cohomology of finite group schemes over a field\fg}, \emph{Invent. Math.}
  \textbf{127} (1997), no.~2, p.~209--270.

\bibitem[Gab62]{Gab}
{\scshape P.~Gabriel} -- {\og Des cat\'egories ab\'eliennes\fg}, \emph{Bull.
  Soc. Math. France} \textbf{90} (1962), p.~323--448.

\bibitem[GZ67]{Gabriel-Zisman}
{\scshape P.~Gabriel {\normalfont \smfandname} M.~Zisman} -- \emph{Calculus of
  fractions and homotopy theory}, Ergebnisse der Mathematik und ihrer
  Grenzgebiete, Band 35, Springer-Verlag New York, Inc., New York, 1967.

\bibitem[Kas82]{Ka}
{\scshape C.~Kassel} -- {\og La {$K$}-th\'eorie stable\fg}, \emph{Bull. Soc.
  Math. France} \textbf{110} (1982), no.~4, p.~381--416.

\bibitem[Kuh98]{K-comp}
{\scshape N.~J. Kuhn} -- {\og Computations in generic representation theory:
  maps from symmetric powers to composite functors\fg}, \emph{Trans. Amer.
  Math. Soc.} \textbf{350} (1998), no.~10, p.~4221--4233.

\bibitem[Lod98]{Loday}
{\scshape J.-L. Loday} -- \emph{Cyclic homology}, second \smfedname,
  Grundlehren der Mathematischen Wissenschaften [Fundamental Principles of
  Mathematical Sciences], vol. 301, Springer-Verlag, Berlin, 1998, Appendix E
  by Mar\'\i a O. Ronco, Chapter 13 by the author in collaboration with
  Teimuraz Pirashvili.

\bibitem[Pfi95]{Pfister}
{\scshape A.~Pfister} -- \emph{Quadratic forms with applications to algebraic
  geometry and topology}, London Mathematical Society Lecture Note Series, vol.
  217, Cambridge University Press, Cambridge, 1995.

\bibitem[Pir00a]{PDK}
{\scshape T.~Pirashvili} -- {\og Dold-{K}an type theorem for
  {$\Gamma$}-groups\fg}, \emph{Math. Ann.} \textbf{318} (2000), no.~2,
  p.~277--298.

\bibitem[Pir00b]{P-hodge}
\bysame , {\og Hodge decomposition for higher order {H}ochschild homology\fg},
  \emph{Ann. Sci. \'Ecole Norm. Sup. (4)} \textbf{33} (2000), no.~2,
  p.~151--179.

\bibitem[Qui72]{Qui}
{\scshape D.~Quillen} -- {\og On the cohomology and {$K$}-theory of the general
  linear groups over a finite field\fg}, \emph{Ann. of Math. (2)} \textbf{96}
  (1972), p.~552--586.

\bibitem[Sch85]{Scharlau}
{\scshape W.~Scharlau} -- \emph{Quadratic and {H}ermitian forms}, Grundlehren
  der Mathematischen Wissenschaften [Fundamental Principles of Mathematical
  Sciences], vol. 270, Springer-Verlag, Berlin, 1985.

\bibitem[Sco00]{Sco}
{\scshape A.~Scorichenko} -- {\og Stable {K}-theory and functor homology over a
  ring\fg}, \smfphdthesisname, Evanston, 2000.

\bibitem[Sta97]{Stanley}
{\scshape R.~P. Stanley} -- \emph{Enumerative combinatorics. {V}ol. 1},
  Cambridge Studies in Advanced Mathematics, vol.~49, Cambridge University
  Press, Cambridge, 1997, With a foreword by Gian-Carlo Rota, Corrected reprint
  of the 1986 original.

\bibitem[Tou08]{Tcu}
{\scshape A.~Touz{\'e}} -- {\og Universal classes for algebraic groups\fg}, \`a
  para\^itre au {\em Duke Math. J.}, 2008.

\bibitem[Tou09]{Touze}
\bysame , {\og Cohomology of classical algebraic groups from the functorial
  point of view\fg}, arXiv:0902.4459, 2009.

\bibitem[Tro02]{T1}
{\scshape A.~Troesch} -- {\og Quelques calculs de cohomologie de composition de
  puissances sym\'etriques\fg}, \emph{Communications in Algebra} \textbf{30(7)}
  (2002).

\bibitem[TvdK08]{TvdK}
{\scshape A.~Touz{\'e} {\normalfont \smfandname} W.~van~der Kallen} -- {\og
  Bifunctor cohomology and cohomological finite generation for reductive
  groups\fg}, \`a para\^itre au {\em Duke Math. J.}, 2008.

\bibitem[Ves08]{CV}
{\scshape C.~Vespa} -- {\og Generic representations of orthogonal groups: the
  functor category {$\mathcal{F}_{quad}$}\fg}, \emph{J. Pure Appl. Algebra}
  \textbf{212} (2008), no.~6, p.~1472--1499.

\end{thebibliography}
